\documentclass[ssy]{amsart}
\usepackage{amssymb}
\usepackage{algorithm}
\usepackage{algorithmic}
\usepackage{verbatim} 
\newtheorem{thm}{Theorem}[section]
\newtheorem{prop}[thm]{Proposition}
\newtheorem{lem}[thm]{Lemma}
\newtheorem{cor}[thm]{Corollary}
\newtheorem{defin}[thm]{Definition}
\newtheorem{example}[thm]{Example}

\newcommand{\mysection}{\setcounter{equation}{0}\setcounter{secnumdepth}{3} 
\section}
\makeatletter
      \def\@setcopyright{}
      \def\serieslogo@{}
      \makeatother
\begin{document}
\author{Michael L. Wenocur}
\address{Palo Alto, California, USA}
\email{mwenocur@acm.org}
\title{New Representations for G/G/1 Waiting Times}

\begin{abstract}
We obtain an explicit representation for the Laplace transform of the waiting time for a wide class of distributions
by solving the Wiener-Hopf factorization problem via the Hadamard product theorem.
Under broad conditions it is shown that this representation
is invertible by an infinite partial fraction expansion. Computational schema illustrated by
a variety of examples demonstrate
the feasibility of numerically solving a rich class of G/G/1 queue possessing either bounded service or arrival times.
In particular very tractable computations are derived for the M/M/1 queue with gated arrivals.

Much of the analysis hinges on the notion of {\it exponential order} given in Definition \ref{expOrdDef} and referenced frequently.
\end{abstract}
\subjclass{Primary 60K25; Secondary 45E10,44A10}
\keywords{G/G/1, Wiener-Hopf Factorization, Waiting Times}

\thanks{Thanks to J. Michael Harrison for teaching me how to conduct research, to Gaurav Khanna whose encouragement helped me cross the finish line, and to my wife Deborah for her manifold kindnesses.}

\dedicatory{Dedicated to Austin J. Lemoine who long ago provided time and encouragement to begin this research.}
\date{\today}
\maketitle

\mysection {Unfinished Workload Computation}
The steady-state unfinished workload $V$
of the queue prior to an arrival
satisfies the equation
\begin{equation}
V \stackrel{\mathcal L}{=} (V+U)^{+}
\end{equation}
Under mild regularity conditions the Laplace transform of $U=Y-X$ may be written as
\begin{equation}
B(\theta)A(-\theta)-1 = \theta\chi(\theta)/\psi(\theta)
\end{equation}
where $A(\theta)={\mathcal E} [\exp(-\theta X)] $ is the Laplace transform of the interarrival times,
$B(\theta)={\mathcal E}[\exp(-\theta Y)]$ is the Laplace transform of the service times, and
$\psi(\theta) ={\mathcal E} [\exp(-\theta V)] $.
For a variety of interesting models we are
able to characterize the asymptotic behavior of zeroes of $\chi$ and the poles of $\psi$. Moreover, in special cases we give explicit
representations for $\chi$ and $\psi$ that are computationally tractable,
with detailed computations for gated M/M/1 and uniform/D/1 processes.

\mysection {Wiener-Hopf Factorization Results}

In this section we present  summary results for representing the unfinished workload's Laplace transform for a rich class of G/G/1 queues. Results are also given for inverting these transforms to obtain the tail probabilities.
\begin{thm}
\label{smithLemma}
Let $A(\theta) ={\mathcal E}[\exp(-\theta X)]$ and
$B(\theta)={\mathcal E}[\exp(-\theta Y)]$ be
analytic in the non-null strip $\{\theta: -\epsilon < \Re{\theta} <\epsilon \}$ with either $X$ or $Y$ non-lattice,
and set
$F(\theta)=B(\theta)A(-\theta)-1$ then $F$ may be factored
as $\theta\chi(\theta)/\psi(\theta)$ where $\chi$ is analytic and non-vanishing on  $\{\theta:  \Re{\theta} <\epsilon_1 \}$ and $\psi$ is analytic and non-vanishing on  $\{\theta: -\epsilon_1 < \Re{\theta} \}$ for some $\epsilon_1 > 0$. 

Moreover, the zeroes of $\chi$ are co-extensive with all the zeroes of
$F$ that have positive real part and the poles of $\psi$ are co-extensive with the zeroes of $F$ that have negative real part.

If $A$ has no essential singularities and its poles are nowhere dense then the poles of $\chi$ are co-extensive with the poles
of $F$ having positive real part. If $B$ has no essential singularities and its poles are nowhere dense then the zeroes of $\psi(\theta)$ are co-extensive with all the poles of $F$ having negative real part.
\end{thm}
{\bf Proof} This theorem is proved in Smith\cite{smith}.

\noindent
{\bf Remark}
The gated M/M/1 queue model analyzed in Section \ref{ComputationalExamples} provides an example of a service time transform which has an essential singularity at $-\mu$. In this case the singularity appears in $\psi$ as a infinite product of terms of the form $(\theta +\mu)/(\theta - s_n)$ where 
$s_n = -\mu + cn^{-1} + O(n^{-2})$.

\begin{thm}
	\begin{flushleft}
	Let $F, A,$ and $B$ satisfy the conditions of Theorem \ref{smithLemma}, with
	 $A$ an entire function such that $|A(\theta)| < C\exp(\lambda |\theta|)$.
	Moreover suppose the right-hand zeroes of $F$
	have the enumeration
	$\{z_n ,\overline{z}_n, 0 < n < \infty \}$
        and $\Re{z_n} = o(n^{1/2})$
	and $\Im{z_n} = 2\pi\alpha^{-1}n + \beta + o(1)$. Then
	\end{flushleft}

	\begin{equation}
		\chi (\theta)= {\chi}_{0} \exp\{\alpha\theta/2\}
		\prod_{k=1}^{\infty} (1- \frac{2 \Re{z_k}\theta}
		{z_k\overline{z}_k}
		+\frac{\theta^{2}}{z_k\overline{z}_k} )
	\end{equation}
\label{chiExpThm}
\end{thm}
	{\bf Proof}
	Theorem \ref{smithLemma}
	implies $\chi$ is an entire function whose
	zeroes occur with the same position and multiplicity
	as those of
	$F$ in the strict right plane.
	The rest of the proof consists of showing $\chi(\theta)$
	satisfies the exponential bounds of Theorem \ref{inf_prod_characterization}.
	Smith\cite{smith} proves that $\chi(\theta)$ is bounded and analytic
	on $\Re{\theta} \leq 0$.
   	Since
	$\chi= F \psi/\theta$ holds in the right plane, the exponential bound on $A$
	implies the existence
	of a constants $C$ and $M$ for which
	$ | \chi (\theta)| \leq C \exp(M\mid\theta\mid)$.
	This completes the proof.
\begin{thm}
		\begin{flushleft}
			Let $F, A,$ and $B$ be as defined above, with
			$A$ analytic around a neighborhood zero, and $B$ an entire function, and where $F$ is 
			of exponential order $-\alpha$ at $-\infty$ where $\alpha < 0$.
			Moreover suppose the left-hand zeroes
			have the enumeration
			$\{z_{n},  n= -\infty , \infty\}$
			such that
			$z_{-n} = \overline{z}_n$,
			$\Re z_n = o(n^{1/2})$ and
			$\Im{z_n} = -2\pi n/ \alpha + \beta + o(1)$,
			and that there exists an $\eta$-sublinear (cf. Definition \ref{subLinDef})
                        sequence of radii $\{r_n, n \geq n_1\}$
			for which
			$|F(\theta)| > \epsilon > 0$ whenever $\Re \theta < 0$ and $|\theta| = r_n$
			then
		\end{flushleft}
\begin{equation}
\label{psiProdEquation}
		\psi (\theta) =\exp\{-\alpha\theta/2\}(1-\theta/z_{0})^{-1}
		\prod_{k=1}^{\infty} (1-\frac{2\Re z_k\theta}{z_k\overline{z}_k}
		+\frac{\theta^{2}}{z_k\overline{z}_k})^{-1}
\end{equation}
		
\label{prodWaitingThm}
\end{thm}
{\bf Proof}
        The proof proceeds by showing that $\psi$ satisfies the conditions of
        Theorem \ref{onThemeOfWeierstrass}.
        Marshall \cite{marshall} shows that $\psi$ may be written in terms of $I$, 
        the transform of the idleness distribution, as
        \begin{equation}
	\psi(\theta)= c_0 (1 - I(-\theta))F(\theta)^{-1}
	\label{idlenessIdentity}
        \end{equation}
    where  $c_0$ is a normalization constant.
        
    The conditions on $r_n$ together with (\ref{idlenessIdentity}) imply
	$|\psi(\theta)| \leq \max(2 c_0 \epsilon^{-1}, 1)$
	for $|\theta| = r_n$.

	Since $F$ is of exponential order $|\alpha|$ at $-\infty$, and $1 - I(-\theta)$ is of exponential order
	$0$ at $-\infty$, it follows from Proposition \ref{prod_exp_order} that $\psi$
	is of exponential order $-|\alpha|$
	at $-\infty$. This shows that $\psi$ satisfies the conditions of Theorem \ref{onThemeOfWeierstrass}.
	
\vspace{0.20in}
\noindent
{\bf Remark}
With slight adjustment to Theorem \ref{onThemeOfWeierstrass} one may loosen the conditions on $B$ to permit its possessing a finite number of simple poles. In this case the representation of $\psi$ would be modified by adding a polynomial factor whose zeroes would be the poles of $B$. An example of such a transform $B$, would be the transform of shifted exponential distribution, ie,
$B(\theta) = \exp(-b\theta)(1 + \theta/\mu)^{-1}$.
\begin{thm}
		\begin{flushleft}
			Let $F, A,$ and $B$ be as in Theorem \ref{prodWaitingThm}
			then
		\end{flushleft}
		\begin{equation}
		\label{psiPartFracExpansionEqn}
		\psi (\theta) =
		     1 + \sum_{n= -\infty}^{\infty} \sum_{j=1}^{k_n}  
		                       A_{n,j}((\theta - z_n)^{-j} - (-z_n)^{-j})
    \end{equation}
     where $A_{n,j}$ is the coefficient of $(\theta - z_n)^{-j}$ in the Laurent expansion of
                $\psi$ around $z_n$, and $k_n$ is the multiplicity of $z_n$.
	
\end{thm}

{\bf Proof} Theorem \ref{prodWaitingThm} implies $\psi$ is the ratio of the
two analytic functions, $\exp\{\alpha\theta/2\}$ and $\gamma$ respectively, with
$\gamma$ not identically zero. Thus $\psi$ is globally meromorphic. Mittag-Leffler
integrability follows from the boundness of $\psi$ on $|\theta| = r_n$. Thus $\psi$
satisfies the conditions of Theorem \ref{mittagTheorem}.
\begin{cor}
\label{corLaplaceTranExpans}
		\begin{flushleft}
			Let $F, A,$ and $B$ be as in Theorem \ref{prodWaitingThm}
			then
		\end{flushleft}
		\begin{equation}
		\psi(\theta) =
		     1 + \sum_{n= -\infty}^{\infty} \sum_{j=1}^{k_n}  
		                       a_{n,j}((1 - \theta/z_n)^{-j} - 1)
                 \label{partFracResult}
                \end{equation}
                
		\begin{equation}
		\psi ^{(\nu)}(\theta) =
		      \sum_{n= -\infty}^{\infty} \sum_{j=1}^{k_n}  
		           a_{n,j}(1 - \theta/z_n)^{-j - \nu} z_n^{-\nu} (\nu + j - 1)! / (j - 1)!
                 \label{partFracDerivResult}
                \end{equation}
\end{cor}
where $a_{n,j} = A_{n,j}(-z_n)^{-j}$ and $\nu > 0$

\begin{thm}
		\begin{flushleft}
			Let $F, A,$ and $B$ be as in Corollary \ref{corLaplaceTranExpans} 
			and suppose that there exist $M_1$ finite and $\gamma >0$ such that $|a_{n,j}| \leq M_1 |n|^{-1}$ and 
			$\Re{z_n} = -\gamma \log|n| + \gamma_0 + o(1)$
			then the tail waiting probability for $t >0 $ is given by
		\end{flushleft}

\begin{equation}
P(W > t) = \lim_{N \rightarrow \infty}
\sum_{n= -N}^N \sum_{j=1}^{k_n} a_{n,j}\sum_{l = 0}^{j - 1}
		                        \exp(z_nt)(-tz_n)^l / l! 
\label{waitingTailExpansion}
\end{equation} 
\end{thm}
{\bf Proof}
Observe by Lemma \ref{WaitingTimeExpFromPsi} the expansion holds for positive $t$ that are points of continuity for the distribution. The l.h.s. of (\ref{waitingTailExpansion}) is a monotone decreasing right continuous function, so it has at most a countable number of discontinuities. By Lemma 
\ref{monotoneContinuityLemma} it suffices to show that the r.h.s. of (\ref{waitingTailExpansion}) is continuous on $[u, \infty]$
for all $u > 0$.
 
Denote the right hand side by $Q(t)$. We proceed by showing $Q(t)$ belongs to $C_0[u, \infty]$ for $u > 0$.  The hypothesis $\Re{z_n} < 0$ for all $n$ implies that each summand is a continuous function converging to zero as $t \rightarrow \infty$ .
Define $Q_0(t)$ as follows
\begin{equation}
Q_0(t) = \sum_{|n| > 0}  a_{n,1}\exp(z_nt) 
\end{equation}
and $Q_1(t) = Q(t) - Q_0(t)$. 

The hypotheses on $z_n$ show that at most a finite number of them can repeat. Thus $Q_1$ is a finite sum of $C_0[u, \infty]$ functions, implying $Q_1$ belongs to $C_0[u, \infty]$. By the hypotheses on $z_n$ for $t \geq u$ we may choose $M$ sufficiently large so that 
$\Re{z_n} < -\gamma\log|n| + M$. Thus
\begin{equation}
|a_{n,1}\exp(\Re{z_n}t) <  |a_{n,1}|\exp(\Re{z_n}u) 
< M_1|n|^{-1 -\gamma u}\exp(Mu)
\end{equation}
\begin{equation}
|Q_0(t)| < \sum_{|n| > 0} |a_{n,1}|\exp(\Re{z_n}u) 
< \sum_{|n| > 0}^{\infty} M_1|n|^{-1 -\gamma u}\exp(Mu) < \infty
\end{equation}

The last equation implies the partials sums of $Q_0(t)$ form a Cauchy sequence in 
$C_0[u, \infty]$. Since $C_0[u, \infty]$ is a complete space it follows that $Q$ is a member of $C_0[u, \infty]$. 

\begin{lem}
\label{monotoneContinuityLemma}
Suppose $f$ is right continuous on $[a, b)$ and $g$ is continuous on $[a, b]$. 
Define $U = \{t \in [a, b]: f(t) \neq g(t)\}$. Suppose $U$ has at most a countable set of points, then $f$ is continuous on $[a, b)$.
\end{lem}
{\bf Proof} Let $W = \{w \in [a,b): f(w) = g(w)\}$. The set $W$ is by hypothesis dense in $[a, b)$, hence for
$t \in [a, b)$ there exists a strictly monotone decreasing sequence in $W$ that converges to $t$. By right continuity of $f$ and the continuity of $g$ it follows that $f(t) = g(t)$. Since $t$ is arbitrary, it follows that $f = g$ on $[a,b)$.
\subsection{Formulas for computing the $A_{n, j}$}
We conclude this section by deriving a product representation of
$A_{n,j}$. We begin 
by noting that $A_{n,j}$ may be computed as
\begin{equation}
\label{a_n_j_Formula}
A_{n,j} (k_n - j)! =  \lim_{\theta \rightarrow z_n}D^{k_n - j} (\psi(\theta) (\theta - z_n)^{k_n}) 
\end{equation}
where $D$ is the derivative operator with respect to $\theta$.
If we introduce the helping functions
$\psi(\theta, n) = \psi(\theta) (\theta - z_n)^{k_n}$, and $\Upsilon(\theta, n) = \log(\psi(\theta, n))$, then
\begin{equation}
A_{n,j} = \psi^{(k_n - j)}(z_n, n)/ (k_n -j)! 
\end{equation}
Equations \ref{psiProdEquation} and \ref{a_n_j_Formula} jointly imply that $A_{n, k_n}$ may be computed as
\begin{equation}
A_{n, k_n} = (-z_n)^{k_n}\exp(\alpha z_n / 2) \prod_{l != n} (1 - z_n / z_l)^{-k_l}
\end{equation}
$A_{n, k_n - 1} = \Upsilon^{(1)}(z_n, n) A_{n, k_n}$ where 
\begin{equation}
\Upsilon^{(1)}(z_n, n) = \alpha/2 + \sum_{m \neq n} k_m (z_m - z_n)^{-1}
\end{equation}
$A_{n, k_n - 2} = (\Upsilon^{(2)}(z_n, n) + \Upsilon^{(1)}(z_n, n)^2)A_{n, k_n}$ where 
\begin{equation}
\Upsilon^{(2)}(z_n, n) = \sum_{m \neq n} k_m (z_m - z_n)^{-2}
\end{equation}
It is a straight forward task to to derive formulas for $k_n \geq 3$ akin to deriving moments from
the cumulant function. Computational experience suggests that even a single repeated root is
unusual.
\mysection{Examples}
In this part we give five examples that illustrate the results of the previous section. In particular we consider the time-gated M/M/1 queue where Poisson arrivals are only admitted only at the end of each unit time interval.

The next two examples study the $E_m$/D/1 queue culminating in an alternative representation of Takacs' formula for the M/D/1 unfinished workload tail distribution. In the final two examples we consider U/D/1 and U/U/1 queues.
\label{ComputationalExamples}
\begin{example}
Gated M/Cox/1 Queue
\end{example}
We study the unfinished work process of the gated M/Cox/1 queue. This is equivalent to D/G/1 system
where 
$A(\theta) =\exp(\theta)$,
$B(\theta) = \exp\{-\lambda + \lambda R(\theta)\}$, and
$F(\theta) = \exp\{-\lambda + \lambda R(\theta) + \theta\} - 1$
and $R(\theta)$ is the Laplace transform of the Cox service time distribution.

It is helpful to recall that
\begin{equation}
\exp(u) - 1 =  u \exp(u/2)\prod_{n=1}^{\infty} (1 + u^2 /(2n\pi)^{2})
\end{equation}

We then see that
\begin{equation}
F(\theta) =  H(\theta) \exp(H(\theta)/2)\prod_{n \neq 0} (1 - iH(\theta) /2n\pi)
\label{gMM1Basic}
\end{equation}

where $H(\theta) = -\lambda + \lambda R(\theta) + \theta$

It follows from the previous equation that $\theta$ is a zero of $F$ if and only
there exists an integer $n$ such that
\begin{equation}
H(\theta) = i2n \pi
\label{gatedRootEqn}
\end{equation}

Note that $\overline{H(\theta)} = H(\overline{\theta})$ so that we may restrict attention to
nonnegative $n$ in \ref{gatedRootEqn}.

Using Rouche's Theorem one can show that for all but a finite number of $n$'s
there is only one root $\theta$ to \ref{gatedRootEqn} with $\Re \theta > 0$.

Specializing to the case of exponential service times with rate $\mu$ we find
that \ref{gatedRootEqn} simplifies to
\begin{equation}
\theta^2 + \theta(\mu - \lambda -i2\pi n) -i2\mu n\pi = 0
\end{equation}

Solving for $\theta$ yields
\begin{equation}
\theta = a + ib_n  \pm(c_n - id_n)
\end{equation}

where
$a = (\lambda - \mu)/2$, $b_n = n\pi$, $x_n = (\lambda - \mu)^2/4 - n^2 \pi^2$,
$y_n = (\mu + \lambda)n\pi$, $c_n = \Re((x_n + iy_n)^{1/2})$, and
$d_n = \Im((x_n + iy_n)^{1/2})$

The formula for the principal square root of $x_n + iy_n$ for $x_n < 0$ is given
\begin{equation}
2^{1/2}c_n = ((x_n^2 + y_n^2)^{1/2} + x_n)^{1/2} = |x_n|^{1/2}((1 + x_n^{-2}y_n^2)^{1/2} -1)^{1/2}
\end{equation}
\begin{equation}
{\rm sgn}(n)2^{1/2}d_n = ((x_n^2 + y_n^2)^{1/2} - x_n)^{1/2}
           = |x_n|^{1/2}((1 + x_n^{-2}y_n^2)^{1/2} +1)^{1/2}
\end{equation}
and by $x_n \geq 0$ is given
\begin{equation}
2^{1/2}c_n = ((x_n^2 + y_n^2)^{1/2} + x_n)^{1/2} = |x_n|^{1/2}((1 + x_n^{-2}y_n^2)^{1/2} +1)^{1/2}
\end{equation}
\begin{equation}
{\rm sgn}(n)2^{1/2}d_n = ((x_n^2 + y_n^2)^{1/2} - x_n)^{1/2}
           = |x_n|^{1/2}((1 + x_n^{-2}y_n^2)^{1/2} -1)^{1/2}
\end{equation}
For the remainder of this discussion we limit attention to  
$n > (\mu -\lambda)2^{-1}\pi^{-1}$.

By employing the local expansions
\begin{equation}
2^{-1/2}((1 + u^2)^{1/2} - 1)^{1/2} = |u|2^{-1} -|u|^3 16^{-1} + O(u^5)
\end{equation}
\begin{equation}
2^{-1/2}((1 + u^2)^{1/2} + 1)^{1/2} = 1 + |u|^2 8^{-1} + O(u^4)
\end{equation}
we obtain
\begin{eqnarray}
c_n & = & y_n |x_n|^{-1/2} 2^{-1} - y_n^3 |x_n|^{-5/2} 16^{-1} + O(n^{-3})\nonumber \\
    & = & (\lambda + \mu)0.5  -
    \mu\lambda(\mu + \lambda) (2n\pi)^{-2} +O(n^{-3})
\end{eqnarray}
and
\begin{eqnarray}
d_n &=& |x_n|^{1/2} + y_n^2 |x_n|^{-3/2}8^{-1} + O(n^{-3}) \nonumber \\
    & = & n\pi + \lambda \mu (2n\pi)^{-1}   + O(n^{-3})
\end{eqnarray}

The zeroes $r_n$ with nonnegative real parts are given by
\begin{eqnarray}
\label{rn_gatedArrivals}
r_n &=& a + ib_n + c_n + id_n \nonumber \\
	  &=& \lambda + i2n\pi +  i\lambda\mu (2n\pi)^{-1} - 
	  \mu\lambda(\mu + \lambda) (2n\pi)^{-2}+ O(n^{-3})
\end{eqnarray}
The zeroes $s_n$ with negative real parts are given by
\begin{eqnarray}
s_n &=& a + ib_n - c_n - id_n \nonumber \\
    & =& -\mu - i\lambda \mu (2n\pi)^{-1}
      + \mu\lambda(\mu + \lambda) (2n\pi)^{-2}+ O(n^{-3})
\end{eqnarray}
Simple algebra shows that $c_n > |a|$ for $n \neq 0$ yielding
$\Re{r_n} > 0 > \Re{s_n}$ for $n \neq 0$.
Moreover $r_0 = 0$ and $s_0 = \lambda - \mu$.

We may now rewrite \ref{gMM1Basic} as
\begin{equation}
F(\theta) = H(\theta)\exp(H(\theta)/2)
             \prod_{n \neq 0} ((\theta - r_n)(\theta -s_n)(\theta + \mu)^{-1}(i2n\pi)^{-1}
\end{equation}
where $H(\theta) = (\theta + \mu -\lambda)\theta(\theta + \mu)^{-1}$.

The power series substitution theorem (cf Apostol \cite{apostol}, Theorem 13-27)
implies that a function formed from the exponentiation of a
function analytic around zero is in turn a function analytic around zero.
Therefore since $R$ is analytic in a neigborhood of zero so must $B$. 
In light of \ref{rn_gatedArrivals} we may apply Theorem $\ref{chiExpThm}$ to show
\begin{equation}
F(\theta) = \theta\chi(\theta)/\psi(\theta)
\label{fChiPsiGated}
\end{equation}
where
\begin{equation}
		\chi (\theta)= (1 - \lambda\mu^{-1} ) \exp\{\theta/2\}
		\prod_{n = 1}^{\infty} (\theta r_n^{-1}-1)(\theta r_{-n}^{-1} -1)
\end{equation}
We may divide $F(\theta)$ by $\theta\chi(\theta)$ to get
\begin{equation}
\psi(\theta) = H_0(\theta) \exp(H_1(\theta))\prod_{n = 1}^{\infty} \Upsilon_n(\theta) \Upsilon_{-n}(\theta)
\label{mm1GatedProdForm}
\end{equation}
where $H_0(\theta) = (1 -\rho)(1 - \lambda (\mu + \theta)^{-1} )^{-1} $,
$H_1(\theta) = \lambda 2^{-1}(1 - \mu(\mu + \theta)^{-1})$,
and $\Upsilon_n(\theta) = ((\theta + \mu)(\theta -s_n)^{-1} (i2n\pi)r_n^{-1})$.
Using $s_n r_n = -i2n\mu\pi$ it is easy to show that $\psi(0) = 1$.

Letting $\theta \rightarrow +\infty$ shows the idle probability is given by 
\begin{equation}
W_0 = (1 -\rho) \exp(\lambda/2) \prod_{n = 1}^{\infty}  4n^2\pi^2 (r_n r_{-n})^{-1}
\end{equation}
The latter equation may be efficiently computed using the approximation given in \ref{eulerProdApprox} and by noting that
\begin{equation}
4\pi^2|r_n|^{-2} = 1 + (\lambda^2 + 2\lambda\mu)(4 \pi^2 n ^2)^{-1} + O(n^{-4})
\end{equation}
Now observe that asymptotically for $\theta$ fixed
\begin{equation}
\Upsilon_n(\theta) = (1 + (s_n + \mu)(\theta - s_n)^{-1})(1 - i\lambda(2n\pi)^{-1} + O(n^{-2}))
\end{equation}
or that
\begin{equation}
\Upsilon_n(\theta) = (1 + i\lambda\mu(2n\pi)^{-1}(\theta - s_n)^{-1})(1 - i\lambda(2n\pi)^{-1}) + O(n^{-2})
\end{equation}
It follows from the previous equation that the product appearing in \ref{mm1GatedProdForm} is absolutely
convergent. Surprisingly this transform is invertible as will be demonstrated henceforth.

We proceed by reparameterizing  \ref{fChiPsiGated} and \ref{mm1GatedProdForm} in terms of
$w = (\theta +\mu)^{-1}$ and $w_n = (s_n + \mu)^{-1}$. Employ
$(\theta + \mu)(\theta - s)^{-1}  = (1 - (s + \mu)w)^{-1}$ in  \ref{mm1GatedProdForm} to get
\begin{equation}
\Psi(w) = (1-\rho) G_0(w) \exp(G_1(w))
          \prod_{n = 1}^{\infty}  g_n(w) g_{-n}(w)
\label{psiWGated}
\end{equation}
where
$G_0(w) = (1-(\mu -\lambda)w)^{-1}$,
$G_1(w) = -\lambda\mu 2^{-1} w + \lambda 2^{-1}$ and $g_n(w) = (1 -(s_n + \mu)w)^{-1}i2n\pi r_n^{-1}$.

Reparameterizing \ref{fChiPsiGated} in terms of $w$ yields
\begin{equation}
\Psi(w) = \phi_0(w)/ \phi_1(w)
\end{equation}
where $\phi_0(w) = (1-\rho)(w^{-1} -\mu)\chi(w^{-1} -\mu)$ and   
$\phi_1(w) = (\exp\{\lambda\mu w - \lambda + w^{-1} -\mu\} - 1)$

Observe that as $|w| \rightarrow \infty$ $\phi_0(w) \rightarrow (\lambda - \mu)\chi(-\mu)$.
One can show that $\phi_1(w)$ remains bounded below by from 0.5 on
$R_n = \{w: |w| = ((2n + 1)\pi + \lambda + \mu)/\lambda\mu\}$ for $n$ sufficiently large. 
Equation \ref{psiWGated} may be used to prove
$\Psi$ is analytic at zero.

It thus follows that $\Psi(w)$ satisfies the conditions of Theorem \ref{mittagTheorem}.
Thus
\begin{equation}
\Psi(w)  = \Psi(0) + \sum_{-\infty < j < \infty} q_j [(w -w_j)^{-1} + w_j^{-1}]
\end{equation}
or reparameterizing in terms of $\theta$ and simplifying yields
\begin{equation}
\psi(\theta)  = \psi_0 + \sum_{-\infty < j < \infty} p_j  (1 - \theta/s_j)^{-1}
\label{gatedMM1Psi}
\end{equation}
where
\begin{equation}
p_j s_j = -{\bf res}(\psi)\left|_{\theta = s_j} = 
-s_j\chi(s_j)/F^\prime(s_j)  = O(j^{-2}) \right.
\end{equation}
which may be reduced to 
\begin{equation}
p_j = \chi(s_j)/ (1 - \lambda(1 +s_j\mu^{-1})^{-2}) = O(j^{-2})
\end{equation}
One may invert \ref{gatedMM1Psi} term by term to get
\begin{equation}
P(V > t) =  \sum_{-\infty < j < \infty} p_j  \exp(s_j t)
\label{gatedMM1Distr}
\end{equation}
\[
\begin{tabular}{|r|r|r|r|r|r|}\hline
\multicolumn{6}{|c|}{Tail Distribution of Unfinished Workload }\\ \hline
$\lambda$&$\mu$&$\rho$&Eqn $\ref{gatedMM1Distr}$&Markov&Time\\ \hline
3.0&4.0&0.750&0.510817&0.509864& 0.0\\ \hline
3.0&4.0&0.750&0.200318&0.200301& 1.0\\ \hline
3.0&4.0&0.750&0.074312&0.074312& 2.0\\ \hline
3.5&4.0&0.875&0.728580&0.727993& 0.0\\ \hline
3.5&4.0&0.875&0.454497&0.454487& 1.0\\ \hline
3.5&4.0&0.875&0.276359&0.276359& 2.0\\ \hline

\end{tabular}
\]
These computations were implemented in Python. The sums appearing
in \ref{gatedMM1Distr} were truncated after the first = 60 terms with the residues approximated by
computing the first 2000 factors.
The finite queue Markov chain approximation was computed iteratively starting from a
queue of 0 and with
a maximum queue length of 200, where the iterative threshold was $10 ^{-10}$
with a maximum of 10000 iterative steps.
\subsection*{Computing The Mean Unfinished Workload}
By differentiating the log of \ref{mm1GatedProdForm} we get the mean unfinished workload

\begin{equation}
-\psi^\prime(0) = -2^{-1}\rho + \rho (\mu -\lambda)^{-1} -
                      \sum_{n = 1}^{\infty} 2 \Re(\mu^{-1} + s_n^{-1})
\label{aveWorkFromPhiGated}
\end{equation}
We may extrapolate the infinite sum in \ref{aveWorkFromPhiGated} by leveraging
\begin{equation}
\Re(\mu^{-1} + s_n^{-1})= - 4^{-1}\lambda (\pi n)^{-2} + O(n^{-4})
\end{equation}
and the Eulerian
identity
\begin{equation}
\sum_{n = 1}^{\infty} n^{-2} = 6^{-1}\pi^{2}
\end{equation}
to get the following computationally tractable formula for the unfinished workload
\begin{eqnarray}
-\psi^\prime(0) + 2^{-1}\rho - \rho (\mu -\lambda)^{-1}  =
                   -\sum_{n = 1}^{\infty} 2 \Re(\mu^{-1} + s_n^{-1}) \nonumber \\
                   =
                   -\sum_{n = 1}^m \left(2 \Re(\mu^{-1} + s_n^{-1}\right) +
                   2^{-1}\lambda (\pi n)^{-2}) + 12 ^{-1} \lambda  + O(m^{-3})
\label{gatedSnMean}
\end{eqnarray}
By differentiating the equation $\theta \chi(\theta) = F(\theta)\psi(\theta)$ twice we find
that $-\psi^\prime(0) =  2^{-1}(1 -\rho) +\rho (\mu - \lambda)^{-1} - \chi^\prime(0)(1 -\rho)^{-1}$
which simplifies to
\begin{equation}
-\psi^\prime(0) = -2^{-1}\rho + \rho (\mu -\lambda)^{-1} -
                      \sum_{n = 1}^{\infty} 2 \Re{r_n} |r_n|^{-2}
\end{equation}
Using \ref{eulerSum} we may extrapolate the previous equation by
\begin{equation}
-\psi^\prime(0) = -2^{-1}\rho + \rho (\mu -\lambda)^{-1} -
                      \sum_{n = 1}^{m} ( 2 \Re{r_n} |r_n|^{-2} - 2\lambda q_n)
                      +  \lambda Q_n  + O(m^{-3})
\label{gatedRnMean}
\end{equation}
where $q_n = \lambda / (4 (n\pi)^2 + a)$, $a  = \lambda^2 + 2\lambda \mu$.
and where
\begin{equation}
Q_n = \sum_{n=1}^\infty q_n = 1 + a^{1/2}\pi * \coth(a^{1/2}\pi) / (2 * a) - 1 / a 
\end{equation}

\vspace{0.20in}
The following tables give an indication of accuracy and the computational tractability
of the above formulas.
\[
\begin{tabular}{|r|r|r|r|r|r|}\hline
\multicolumn{6}{|c|}{Computed Mean Unfinished Workload}\\ \hline
$\lambda$&$\mu$&$\rho$&Equation $\ref{gatedRnMean}$&Equation $\ref{gatedSnMean}$&Finite Queue\\ \hline
3.0&4.0&0.750&0.53620286355&0.53620286365&0.53620286352\\ \hline
3.5&4.0&0.875&1.49447474664&1.49447474664&1.49447473470\\ \hline
\end{tabular}
\]
\[
\begin{tabular}{|r|r|r|r|r|r|}\hline
\multicolumn{6}{|c|}{Average Execution Times in Milliseconds}\\
\multicolumn{6}{|c|}{on a 450 MegaHertz Pentium II}\\ \hline
$\lambda$&$\mu$&$\rho$&Equation $\ref{gatedRnMean}$&Equation $\ref{gatedSnMean}$&Finite Queue\\ \hline
3.0&4.0&0.750&21&21&22900\\ \hline
3.5&4.0&0.875&21&21&80000\\ \hline
\end{tabular}
\]
These computations were implemented in Python. The sums appearing
in \ref{gatedRnMean} and in \ref{gatedSnMean} were truncated at $m = 1000$. 
The finite queue approximation was computed iteratively starting from a queue of 0 and with 
a maximum queue length of 200, where the iterative threshold was $10 ^{-12}$
with a maximum of 800 iterative steps.

\begin{example}
Em/D/1
\end{example}
By rescaling time, so that the service time is one, we may take
$A(-\theta) = (1 - \theta/\lambda)^{-m}$ and $B(\theta) = \exp(-\theta)$ where $m  > \lambda$.
Let $u_1, u_2, ... u_{m-1}$ such that $F(u_i) = 0$ and $\Re {u_i} > 0$.
We now can write
\begin{equation}
\psi(\theta) = (m/\lambda - 1)\theta(1 - \theta/\lambda)^{-m}
F(\theta)^{-1}\prod_{0< i < m}(1 - \theta/u_i)
\label{EnDLapT}
\end{equation}
The equation $F(\theta) = 0$ is equivalent to $\sigma(z) = 0$ where
 $z = \lambda  -\theta$, $\beta = m\ln(\lambda) - \lambda$ and $\sigma$ is defined by \ref{tDefEqn}.
Since $\lambda < m$ it follows that $\sigma(m) \neq 0$,
and therefore by Corollary \ref{tizeroes} it follows that $F$ has only simple zeroes.

For the remainder of this example we will consider the case that $m = 2$ and $\lambda = 1$
Equation \ref{EnDLapT} simplifies to
\begin{equation}
\psi(\theta) = (\theta(1 - \theta/u_1)/(\exp(-\theta) - (1 - \theta)^2)
\end{equation}
where $u_1 = 1.477670$.

The mean, mean square and mean cube waiting time are given by $1/u_1 - 1/2$, 
$5/6 - 1/ u_1$,  and $(5/u_1 - 3)/2$ respectively. Computing with the first thousand terms we find the approximate mean value is

\[
\begin{tabular}{|r|r|r|r|}\hline
\multicolumn{4}{|c|}{Exact versus approximate mean waiting time statistics}\\ \hline
Computational Method & Mean Wait & Mean Square Wait & Mean Cube Wait\\ \hline
Exact              &0.176775&0.156592276220&0.1918526427820\\ \hline
Eqn \ref{spectralMomentFormula} with N = 1000 &0.176741&0.156592276251&0.1918526427803\\ \hline
\end{tabular}
\]
\begin{example}
M/D/1
\end{example}
By rescaling we can take the service time to be 1 with the arrival rate $\lambda < 1$.
The exact formula for waiting time distribution is given by (cf. Takacs\cite{takacs})
\begin{equation}
P\{V \leq t\} = (1 - \rho) \sum_{n = 0}^{\lfloor t \rfloor} \exp(-\lambda(n - t))
                (\lambda(n - t))^n /n!
\end{equation}

\[
\begin{tabular}{|r|r|r|r|r|}\hline
\multicolumn{5}{|c|}{Exact versus approximate tail probabilities}\\
\multicolumn{5}{|c|}{for $\lambda = 1/3$ and $D = 1$}\\ \hline
\multicolumn{1}{|r|}{Terms} & \multicolumn{1}{r|}{Equation \ref{waitingTailExpansion}} & \multicolumn{1}{r|} {Exact Tail} &
               \multicolumn{1}{c|}{Absolute Error} & \multicolumn{1}{r|}{Time}\\ \hline
10  & 0.271886491 & 0.275397300 &  0.003510809  &  0.25\\ \hline
10  & 0.212919003 & 0.212426391 &  0.000492611  &  0.50\\ \hline
10  & 0.070737664 & 0.069591717 &  0.001145947  &  1.00\\ \hline
10  & 0.011647294 & 0.011646734 &  0.000000560  &  2.00\\ \hline
100  & 0.275606488 & 0.275397300 &  0.000209187  &  0.25\\ \hline
100  & 0.212443434 & 0.212426391 &  0.000017042  &  0.50\\ \hline
100  & 0.069704561 & 0.069591717 &  0.000112845  &  1.00\\ \hline
100  & 0.011646735 & 0.011646734 &  0.000000001  &  2.00\\ \hline 
1000  & 0.275409123 & 0.275397300 &  0.000011823  &  0.25\\ \hline
1000  & 0.212426937 & 0.212426391 &  0.000000545  &  0.50\\ \hline
1000  & 0.069602977 & 0.069591717 &  0.000011261  &  1.00\\ \hline
1000  & 0.011646734 & 0.011646734 &  0.000000000  &  2.00\\ \hline
\end{tabular}
\]

\begin{example}
U/D/1
\end{example}
In this example we consider an arrival process governed by a uniform[$a$, $b$] and
a constant time service time of $D$, ie,
$A(-\theta) = (\exp(b\theta) - \exp(a\theta))/((b - a)\theta)$ and $B(\theta) =\exp(-D\theta)$ where $b + a  > 2 D$.

$F(\theta) = (\exp((b - D)\theta) - \exp((a - D)\theta))/((b - a)\theta) - 1$. The helper function is given 
$H(\theta) =  - \exp((a - D)\theta))/((b - a)\theta) - 1$.

We now restrict attention to the case where $a, b, D$ respectively equal $0, 6, 1$  with a traffic intensity of $1/3$.
\[
\begin{tabular}{|l|r|r|r|r|}\hline
\multicolumn{5}{|c|}{Comparison of Cumulant Computational Methods}\\
\multicolumn{5}{|c|}{Telescoped Summation vs Truncated Sum vs Finite Expansion}\\ \hline
Method& Number &  Cumulant &  Cumulant &  Cumulant \\
& of Terms &$\kappa_1$&$\kappa_2$ & $\kappa_3$ \\  \hline

Telescoped &5       &0.1095808& 0.0838003& 0.0795173\\ \hline
Equation \ref{FirstCumulantForm}, \ref{JthCumulantForm}  &1000    &0.1089962& 0.0838510& 0.0795173\\ \hline
Equation \ref{FirstCumulantForm}, \ref{JthCumulantForm}  &10000   &0.1095107& 0.0838054& 0.0795173\\ \hline
Equations \ref{spectralMomentFormula} &1000     &0.1096221&	0.0837913&	0.0795084 \\ \hline
Equations \ref{spectralMomentFormula} &4000     &0.1095911&	0.0837981&	0.0795150 \\ \hline
\end{tabular}
\]

\[
\begin{tabular}{|c|r|}\hline
\multicolumn{2}{|c|}{2000 Term Spectral}\\
\multicolumn{2}{|c|}{Expansion for Tail Distribution}\\ \hline
{\bf t}	&${\bf P\{W > t\}}$\\ \hline
0.000000 & 0.184930\\ \hline
0.250000 & 0.143236\\ \hline
0.500000 & 0.101570\\ \hline
0.750000 & 0.059903\\ \hline
1.000000 & 0.018440\\ \hline
1.250000 & 0.011422\\ \hline
1.500000 & 0.006322\\ \hline
1.750000 & 0.002958\\ \hline
2.000000 & 0.001330\\ \hline
2.250000 & 0.000718\\ \hline
\end{tabular}
\]

We computed the value for $t = 0$ using equation \ref{helperBasedComputationOfIdleProb} to compute the probability of zero waiting and then subtracting that from one instead of using \ref{waitingTailExpansion} because the convergence of the spectral expansion at zero is at most conditional and in any case would be very slow.
\begin{example}
U/U/1
\label{U_U_Example}
\end{example}
In this example we consider an arrival process governed by a unif[$a_0$, $a_1$] distribution and
a service time distributed as unif[$b_0$, $b_1$], ie,
$A(-\theta) = (\exp(a_1\theta) - \exp(a_0\theta))/((a_1 - a_0)\theta)$ and $B(\theta) =(\exp(-b_0\theta) - \exp(-b_1\theta))/((b_1 - b_0)\theta)$ where $a_1 + a_0  > b_0 + b_1$.

$F(\theta) = \gamma (\exp((a_1 - b_0)\theta) -\exp((a_1 - b_1)\theta)
										 -\exp((a_0 - b_0)\theta) +\exp((a_0 - b_1)\theta))\theta^{-2} - 1$
where $\gamma^{-1} = (a_1-a_1)(b_1 - b_0)$. 
We now will introduce the following auxiliary functions to $F$, $H$ the component of $F$ that does not uniformly vanish on the left-hand side of the complex plane, $T$ the dominant term of $H$ needed for proving that $\psi$ has a partial fraction expansion.

Setting $\alpha_3 = a_1 - b_0$, 
$\alpha_2 = \max(a_1 - b_1, a_0 - b_0)$, $\alpha_1 = \min(a_1 - b_1, a_0 - b_0)$, and
$\alpha_0 = a_0 - b_1$. Observe $\alpha_2 > 0$ to satisfy the stability equation. Similarly $\alpha_0 \geq 0$ would imply waiting time is always zero. If $\alpha_1 < 0$ the helper function is given 
$H(\theta) = \gamma(\exp(\alpha_0\theta) - \exp(\alpha_1\theta))\theta^{-2} - 1$
otherwise the helper function is given by
$H(\theta) = \gamma\exp(\alpha_0\theta)\theta^{-2} - 1$.

$T(\theta) = \gamma\exp(\alpha_0\theta)\theta^{-2} - 1$.

We now restrict attention to the case where $a_0, a_1, b_0, b_1$ respectively equal $0, 5, 1, 2$  with a traffic intensity of $0.6$.
\[
\begin{tabular}{|l|r|r|r|r|}\hline
\multicolumn{5}{|c|}{Comparison of Cumulant Computational Methods}\\
\multicolumn{5}{|c|}{Telescoped Summation vs Truncated Sum vs Finite Expansion}\\ \hline
Method& Number &  Cumulant &  Cumulant &  Cumulant \\
& of Terms &$\kappa_1$&$\kappa_2$ & $\kappa_3$ \\  \hline

Telescoped &4       &0.4575838 &0.6797302 &1.4058925\\ \hline
Equations \ref{FirstCumulantForm}, \ref{JthCumulantForm} &1000    &0.4555881 &0.6799327 &1.4058925\\ \hline
Equations \ref{FirstCumulantForm}, \ref{JthCumulantForm}  &100000  &0.4575545 &0.6797323 &1.4058925\\ \hline
Equation \ref{spectralMomentFormula}  &1000    &0.4576456 &0.6796736 &1.4058054\\ \hline
Equation \ref{spectralMomentFormula}  &10000   &0.4575899 &0.6797245 &1.4058840\\ \hline
\end{tabular}
\]

\[
\begin{tabular}{|c|r|}\hline
\multicolumn{2}{|c|}{5000 Term Spectral}\\
\multicolumn{2}{|c|}{Expansion for Tail Distribution}\\ \hline
{\bf t}	&${\bf P\{W > t\}}$\\ \hline
0.000000 &0.389364\\ \hline
0.250000 &0.339889\\ \hline
0.500000 &0.290281\\ \hline
0.750000 &0.240581\\ \hline
1.000000 &0.190809\\ \hline
1.250000 &0.144900\\ \hline
1.500000 &0.107201\\ \hline
1.750000 &0.078343\\ \hline
2.000000 &0.058953\\ \hline
2.250000 &0.045736\\ \hline
\end{tabular}
\]

Note that the value for $t = 0$ was computed using equation \ref{helperBasedComputationOfIdleProb} because the convergence of the spectral expansion at zero is at best conditional.
\mysection{Computational Schema}
\label{computationSchemaSection}
The computation of the waiting time probabilities, idleness probability, and cumulants requires evaluating infinite products and sums involving the poles of $\psi$. The key to telescoping these computations is to find closed form functions whose zeroes grow asymptotically close to the poles of $\psi$. Fortunately we may obtain such functions
from $F$, by removing those parts of $F$ which vanish as $\theta  \rightarrow \infty$ with $\rm{re}(\theta) < 0 $. The generation of such helper functions from F is analyzed in detail in Section \ref{transformAsympSection}.

To motivate the discussion below we reintroduce Example \ref{U_U_Example} of Section \ref{ComputationalExamples}.

\begin{eqnarray}
A(-\theta) & = & (\exp(5\theta) - 1)(5\theta)^{-1} \nonumber \\
B(\theta) & = &  (\exp(-\theta) - \exp(-2\theta)) \theta^{-1}  \nonumber \\F(\theta) & = & (\exp(4\theta) - \exp(3\theta) + \exp(-2\theta) - \exp(-\theta)) (5 \theta^2)^{-1} - 1  \nonumber \\
f(\theta) & = & F(\theta) \theta^2 \nonumber \\
H(\theta) & = &  (\exp(-2\theta) - \exp(-\theta))(5 \theta^2)^{-1}  - 1\nonumber \\
h(\theta) & = & H(\theta)\theta^2 \nonumber \\
p(\theta) &= &\theta (1 - \theta / w_1) ( 1- \theta/w_2) \nonumber \\
q(\theta) &=& p(\theta)\theta^{-2} \nonumber \\
H(\theta) &=& q(\theta)\exp\{-\theta\}
	\prod_{k=1}^{\infty} (1- \frac{2\Re w_k\theta}{w_k\overline{w}_k}
	+\frac{\theta^{2}}{w_k\overline{w}_k}) \nonumber \\
T(\theta) & = &  \exp(-2\theta)(5 \theta^2)^{-1}  - 1\nonumber \\
t(\theta) & = &  \exp(-2\theta)  - 5 \theta^2 \nonumber 
\end{eqnarray}
\begin{verbatim}
FLeftZeroes =  -1.112636162915984,
               -2.362945135569372  +  4.24463938127872j,
               -2.894232391480601  +  7.457728791764555j,
               -3.214439899719532  + 10.723473915095289j
HZeroes     =  -0.329175737104105, 
               -1.107062156887795,
               -2.3629486802831456 +  4.244641429104492j,
               -2.8942321368738169 +  7.457728708110446j,
               -3.2144399507792216 + 10.723473898300561j
TZeroes     =   0.3235824310955529,
               -0.8692442923882397 +  0.618192635996526j,
               -2.3782928911558301 +  4.196823183937124j,
               -2.8870718848316352 +  7.485890763202157j,
               -3.2185697792830807 + 10.703472361378987j
\end{verbatim}

where the constants for $p$ are given by $w_1 =  -0.32917573710410553$ and $w_2 = -1.107062156887795$.

Notice that for $F$ and $H$ we have defined computationally convenient stand-ins $f$ and $h$ that are more easily differentiated and are well behaved around zero. It is these versions which we use to compute the zeroes necessary for computing the waiting time moments and distribution.
\subsection{Finding the Roots of $F$ and $H$}
The availability of interesting quantities of the waiting time distribution is predicated on the ability to extract the zeroes of $F$ which lie in the left half of the complex plane. Moreover, computationally efficient formulas require extracting a large number of $H$'s zeroes. Fortunately, asymptotic considerations show that except for a very small number of zeroes of $F$ and $H$ lying near the origin all the  zeroes of $H$ are very close to all the left hand zeroes of $F$ as illustrated by the table of zeroes given above. 

We first address the finding of zeroes $z_{n+1}$ and $w_{n+1}$ given zeroes $z_n$ and $w_n$ of $F$ and $H$ respectively. We assume that the sequences $\{z_k\}$ and $\{w_k\}$ satisfy the following typical asymptotic equations
$\Im{z_k} = (2k + m/2)\pi\alpha^{-1} + o(1)$,
$\Re {z_k} = m\ln(2k\pi)\alpha^{-1} + O(1)$, and
$|z_k - w_k| = O(k^{-c})$ where $c > 0$.
\newpage
\subsubsection{Computing Zeroes Sequentially With Newton-Raphson Method}
Computational experience shows that one may very efficiently compute the helper functions zeroes conjointly with the left hand zeroes of $F$ by a Newton-Raphson iteration in the complex plane. As the zeroes grow in size the computation of a root z of $F$ from a root $w$ of $H$ requires only 2 or 3 iterative steps. Moreover the computation of next set roots can be very accurately initialized by adding a small constant to the current root $w$.

The iterative root computation is given in Algorithm \ref{alg1} using a Pascal like syntax.
\begin{algorithm}
\caption{Estimate $z_{n+1}$ and $w_{n+1}$ given $ w_n$ to within $\epsilon$ using at most $N$ Newton-Raphson iterations }
\label{alg1}
\begin{algorithmic}[1]
\REQUIRE $|h(w_n)| < \epsilon$
\ENSURE $|f(z_{n+1})| < \epsilon$ and $|h(w_{n+1})| < \epsilon$
\STATE $w = w_n + i2\pi/\alpha$
\FOR{$0 \leq j < N$}
\STATE $w = w -h(w)/h^\prime(w)$
\IF{ $|h(w)| < \epsilon$}
\STATE break
\ENDIF
\ENDFOR
\IF{j== N}
\STATE return error
\ENDIF
\STATE $w_{n+1}= w$
\STATE $z = w$
\FOR{$0 \leq j < N$}
\STATE $z = z -f(z)/f^\prime(z)$
\IF{ $|f(z)| < \epsilon$}
\STATE break
\ENDIF
\ENDFOR
\IF{j== N}
\STATE return error
\ENDIF
\STATE $z_{n+1}= z$
\STATE return $(z_{n+1}, w_{n+1})$
\end{algorithmic}
\end{algorithm}

\subsubsection{Approaches to Finding Zeroes Near the Origin}
Finding the zeroes of $H$ near to the origin and the left hand zeroes of $F$ is somewhat problematic but  Schaefer and Bubeck  \cite{rootFinder} provide a systematic approach to computing all required zeroes near the origin. The author has also tried using a homotopic approach with indifferent results.
\subsection{Computing Waiting Time Probabilities}

The exploitation of formula \ref{waitingTailExpansion} requires the efficient computation of
$a_{n,j}$. But since  $a_{n,j} = A_{n,j}(-z_n)^{-j}$ it suffices to focus
attention on computing $A_{n,j}$ through the use of equation \ref{a_n_j_Formula}. Computational experience suggests that poles are rarely repeated, so we simplify the presentation by stipulating that all $k_n = 1$, but we note that this same framework can be easily extended to the case of repeated poles, and in any case there are only a finite number of repeated poles in the general case as well.

\begin{equation}
\label{seriesCoefProdForm}
a_n = \exp(\alpha z_n / 2) \prod_{k \neq n} (1 - z_n / z_k)^{-1}
\end{equation}

The key to computing \ref{seriesCoefProdForm}
is estimating the bilateral tail product ${\bf a}_n^{(k)}$ defined by
\begin{equation}
{\bf a}_n^{(k)} = \prod_{j = n + k + 1}^\infty
(1 - z_n / z_j)^{-1}(1 - z_n / \overline{z}_j)^{-1}
\end{equation}

Computing the partial product $a_n^{(k)} = a_n / {\bf a}_n^{(k)}$ is a straightforward task.

To compute the tail product ${\bf a}_n^{(k)}$ we posit the existence of a helper function $H$ (cf. equation \ref{H_representation}) of closed form having an expansion
\begin{equation}
\label{helperFunctionEquation}
H(\theta) = \exp(-\alpha \theta/2) q(\theta) \prod_{j = n_0}^\infty (1 - \theta/w_j)
(1 - \theta/\overline{w}_j)
\end{equation}
where $w_k = z_k + O(k^{-l})$ for $l \geq 1$ and  
$q$ is a rational function such that
\begin{equation}
\label{helperFunctionEquation2}
q(\theta) = q_0  \theta^{-m_q}  \prod_{j = 0}^{n_q - 1} (1 - \theta/u_j)
\end{equation}
and $n_q - m_q = 2 n_0 - 1$

Defining ${\bf h}_n^{(k)}$ as
\begin{equation}
{\bf h}_n^{(k)} = \prod_{j = n + k + 1}^\infty
(1 - z_n / w_j)(1 - z_n / \overline{w}_j)
\end{equation}
We may compute ${\bf h}_n^{(k)}  = H(z_n) / ((h_n^{(k)} (1 - z_n/w_n)) $  where 
\begin{equation}
h_n^{(k)} = (1 - z_n/w_n)^{-1}\exp(-\alpha z_n/2) q(z_n) 
\prod_{j = n_0}^{n + k} (1 - z_n/w_j) (1 - z_n/\overline{w}_j)
\end{equation}
One may robustly compute $H(z_n) / (1 - z_n/w_n)$ using a Taylor expansion of $H$ around
$w_n$ by exploiting the identity
\begin{equation}
H(z_n)(1 - z_n/w_n)^{-1}  = -w_n(H(z_n) - H(w_n))(z_n - w_n)^{-1}
\end{equation}
i.e., 
\begin{equation}
(H(z_n) - H(w_n))(z_n - w_n)^{-1} = H^{(1)}(w_n) + H^{(2)}(w_n)(z_n - w_n)/2 + ...
\end{equation}
Using ${\bf h}_n^{(k)}$
We estimate $a_n$ by $\hat{a}_n$ where
\begin{equation}
\hat{a}_n = a_n^{(k)}/{\bf h}_n^{(k)}
\end{equation}
We now estimate the relative error that results from replacing
${\bf a}_n^{(k)}$ by $1/{\bf h}_n^{(k)}$

\begin{equation}
{\bf a}_n^{(k)}{\bf h}_n^{(k)} = \prod_{j = n + k + 1}^\infty
\frac{z_j\overline{z}_j (w_j - z_n)(\overline{w}_j -z_n)}
{w_j\overline{w}_j (z_j - z_n)(\overline{z}_j -z_n) }
\end{equation}

\begin{equation}
\label{transformed_a_h_estimate}
{\bf a}_n^{(k)}{\bf h}_n^{(k)} = \prod_{j = n + k + 1}^\infty
	(\eta_j^{(0)} \overline{\eta}_j^{(0)} \eta_{j, n}^{(1)} \eta_{j, n}^{(2)})
\end{equation}
where $\eta_j^{(0)} = z_j/w_j$,
$\eta_{j, n}^{(1)}= (w_j - z_n)/(z_j - z_n)$ and $\eta_{j, n}^{(2)}= (\overline{w}_j -z_n)/(\overline{z}_j -z_n)$.

By virtue of the asymptotic hypotheses on $z_j$ and $w_j$ we may deduce

\begin{equation}
\eta_j{(0)} = 1 + (z_j - w_j)/w_j = 1 + O(j^{-l - 1})
\end{equation}

\begin{equation}
\eta_{j, n}^{(1)} = 1 + (w_j - z_j)/(z_j - z_n) = 1 + O(j^{-l}(j-n)^{-1})
\end{equation}

\begin{equation}
\eta_{j, n}^{(2)} = 1 + (\overline{w}_j - \overline{z}_j)/(\overline{z}_j - z_n) = 1 + O(j^{-l}(j+n)^{-1})
\end{equation}

Combining the three above equations we get
\begin{equation}
\label{combinedEtaFor_a_h_estimate}
\eta_j^{(0)} \overline{\eta}_j^{(0)} \eta_{j, n}^{(1)} \eta_{j, n}^{(2)} =
1 + O(j^{-l}(j-n)^{-1})
\end{equation}
Setting $j = n + k $ in equations \ref{transformed_a_h_estimate} and  \ref{combinedEtaFor_a_h_estimate} 
yields
\begin{equation}
{\bf a}_n^{(k)}{\bf h}_n^{(k)} = 1 + O(k^{-l})
\end{equation}

\subsection{Computing The Idle Probability}
Motivated by the inequality \ref{wk_xToZk_xRatioEst} we sketch how the idle probability can be efficiently computed given the existence of closed form helping function $H(x)$.
where 
\begin{equation}
\label{hProdFormForIdleProb}
H(\theta) = \exp(-\alpha\theta/2)q(\theta) 
						\prod_{k = n_0}^\infty(1-\theta/w_k)(1 - \theta/\overline{w}_k)
\end{equation}
where $q$ is the ratio of the polynomials $d_0$ to the polynomial $d_1$ such that 
deg($d_0$) - deg($q_1$) $= 2 n_0 - 1$.
\begin{equation}
P\{W = 0\} = \lim_{x \rightarrow \infty} \psi(x)
\end{equation}

Assuming the existence of a function $H$ as given by equation \ref{H_representation} such that
$\lim_{x \rightarrow \infty} H(x) = -1$ we may rewrite the idleness probability as
\begin{equation}
P\{W = 0\} = -\lim_{x \rightarrow \infty} \psi(x) H(x)
\end{equation}

Exploiting the representations of $H$ and $\psi$ respectively given by 
\ref{hProdFormForIdleProb} and \ref{psiProdEquation} we may deduce
\begin{equation}
\label{emptyQueueProbRepr}
P\{W = 0\} = -\lim_{x \rightarrow \infty}r(x) \prod_{k = n_0}^\infty \frac{z_k\overline{z}_k (w_k - x)(\overline{w}_k -x)}
{w_k\overline{w}_k (z_k - x)(\overline{z}_k - x) }
\end{equation}

where $r(x)$ is defined by
\begin{equation}
r(x) = q(x) \prod_{k = -n_0 + 1}^{n_0 - 1} (1 - x/z_k)^{-1}
\end{equation}

Lemma \ref{limitOfMobiusRatios} can now be applied to equation \ref{emptyQueueProbRepr} to show
\begin{equation}
\label{helperBasedComputationOfIdleProb}
P\{W = 0\} = -\lim_{x \rightarrow \infty}r(x) \prod_{k = n_0}^\infty \frac{z_k\overline{z}_k }
{w_k\overline{w}_k }
\end{equation}

\subsection {Moment Computation Considerations}

We now study the computation of waiting time moments via the product expansion given
in Theorem \ref{prodWaitingThm}. It is very convenient to use the cumulant function of
$K(\theta) = \log(\psi(-\theta)$ to
derive the waiting time moments.   The cumulants
defined by $\kappa_j = d^j(K(t))/(dt)^j |_{t=0}$ satisfy the following relationship to the moments
(Cf. Cramer\cite{cramer_cit}, p 186).
\begin{eqnarray}
m_1 &=& \kappa_1 \\
m_2 &=& \kappa_2 + {\kappa_1}^2 \\
m_3 &=& \kappa_3 + 3\kappa_1 \kappa_2 + {\kappa_1}^3
\end{eqnarray}

The equation for $K$ is given by

\begin{equation}
K(\theta) = -\alpha\theta/2 - \sum_{k=-\infty}^{\infty} \log(1+ \theta/z_k)
\end{equation}

Differentiating the previous equation once and setting $\theta$ to zero yields
\begin{equation}
\label{FirstCumulantForm}
\kappa_1 = -\alpha/2 - \sum_{k=-\infty}^{\infty} 1/z_k  = 
         -\alpha/2 -z_0^{-1} - 2\sum_{k=1}^{\infty}  \Re z_k|z_k|^{-2}
\end{equation}

and in for $j > 1$

\begin{equation}
\label{JthCumulantForm}
\kappa_j = \sum_{k=-\infty}^{\infty}(j-1)!(-z_k)^{-j} =
         (j-1)!(-1)^j \left ({z_0}^{-j} + \sum_{k=1}^{\infty}  
         2\Re (z_k^j)|z_k|^{-2j}\right )
\end{equation}

Typically $\Re z_k = O(\log(k))$ and 
$\Im{z_k} = 2k\pi/\alpha$, so that the series for $\kappa_1$,
$\kappa_2$ and $\kappa_3$ converge
at rates $O(\log(k)k^{-1})$, $O(k^{-1})$, and $O(\log(k)k^{-3})$ respectively.
The curious result is that low order cumulants are trickier to compute than higher order cumulants.
Assuming the existence of a closed form function $h$ as given
in equations \ref{helperFunctionEquation} and \ref{helperFunctionEquation2} one can easily compute cumulants $\kappa_1$ and $\kappa_2$ 
as shown below.
   
Let $g(\theta) = \log(h(\theta)\theta^{m_p})$
so that
\begin{equation}
g^\prime(0) = -\alpha/2 - \sum_{k = 0}^{n_p}u_k^{-1} -  
							\sum_{k = n_0}^\infty (w_k^{-1} + \overline{w}_k^{-1})
\end{equation}
and for $j > 1$
\begin{equation}
g^{(j)}(0) = -(j-1)!\left (\sum_{k = 0}^{n_p}u_k^{-j} +  
							\sum_{k = n_0}^\infty (w_k^{-j} + \overline{w}_k^{-j})\right )
\end{equation}
Using the binomial theorem on the asymptotic estimate $w_k = z_k + O(k^{-l})$ yields
\begin{equation}
w_k^{-j} = z_k^{-j} + O(k^{-j -l - 1})
\end{equation}
It then follows that
\begin{equation}
W_n^j = \sum_{k = n}^\infty (w_k^{-j} + \overline{w}_k^{-j}) = Z_n^j + O(n^{-j -l})
\end{equation}
where $Z_n^j$
\begin{equation}
Z_n^j = \sum_{k = n}^\infty (z_k^{-j} + \overline{z}_k^{-j}) 
\end{equation}
We can compute $W_n^j$ using 
\begin{equation}
W_n^j = -g^{(j)}(0)/(j-1)! - \alpha \delta_1(j)/2 -\sum_{k = 0}^{n_p} u_k^{-j} - \sum_{k = n_0}^{n - 1} (w_k^{-j} + \overline{w}_k^{-j})
\end{equation}
and approximate $\kappa_j$ with an error of $ O(n^{-j -l})$ by
\begin{equation}\
\kappa_j = -\alpha \delta_1(j)/2 + (j-1)!(-1)^j\left ({z_0}^{-j} + \sum_{k=1}^{n - 1}  
2\Re(z_k^j)|z_k|^{-2j} + W_n^j \right)  
\end{equation}
where $\delta_1(j)$ is $1$ for $j=1$ and $0$ otherwise.
\subsubsection{Moments Via Spectral Expansion}
The following expansion for the $\nu$-th waiting time moment $\mu _{\nu}$ may be derived by setting $\theta$ to zero in equation \ref{partFracDerivResult}.
\begin{equation}
		\mu _{\nu} = (-1)^{(\nu)}\phi^{(\nu)}(0) = \lim_{N \rightarrow \infty}
		      \sum_{n= -N}^N \sum_{j=1}^{k_n}  
		           a_{n,j} z_n^{-\nu} (\nu + j - 1)! / j!
                 \label{spectralMomentFormula}
\end{equation}

\subsection{Asymptotic Computation of the Spectral Coefficients}
\label{asymCompOfSpctrCoefSection}
In order to exploit equation \ref{waitingTailExpansion} it is a necessity to compute the spectral coefficients $a_{n,j}$ for $j \leq k_n$. For all but a finite number of poles $k_n = 1$ hence we restrict attention  to those $n$ for which $k_n = 1$ and simplify notation by setting $a_n = a_{n,1}$. Keeping in mind that $z_n a_n = A_n$ it suffices to compute $A_n$ via 
\ref{a_n_j_Formula} which simplifies to
\begin{equation}
A_n = \lim_{\theta\rightarrow z_n}\psi(\theta)(\theta -z_n)
\end{equation}

We may rewrite the previous equation using \ref{idlenessIdentity} to obtain
\begin{equation}
A_n = \lim_{\theta\rightarrow z_n}c_0(1 -I(-\theta))F(\theta)^{-1}(\theta -z_n) =
c_0(1-I(-z_n))/F^{\prime}(z_n)
\end{equation}

Since the idleness distribution is concentrated on the interval $[0, c]$ where $c = a_1 - b_0$ it is heuristically suggestive to develop $I(z_n)$ as
as follows (cf \ref{finiteExpansionOfLaplaceTransform} of Section \ref{transformAsympSection})

\begin{equation}
I(\theta) \approx I^{(N)}_0(\theta) - I^{(N)}_c(\theta)
\end{equation}

where $N$ is a small integer say no more than 3 and
\begin{equation}
I_u^{(N)}(\theta) = \exp(-\theta u)\sum_{j = 0}^N d_{j, u}\theta^{-j - 1}
\end{equation}

where $d_{j, 0} $ and $d_{j, c}$ are unknown constants to be determined by least squares fit relative to exactly computed coefficients.

\mysection{Helper Functions and Transform Asymptotics}
\label{transformAsympSection}
We present the asymptotics needed to justify Theorem \ref{prodWaitingThm}
and to obtain the closed form helper functions for telescoping the series and products which arise in the tail probability and moment computations. We also provide estimates of how well the helper functions' zeroes approximate those of $F$.

\subsection{Results for Smooth Densities on Finite Intervals}
We begin by developing an asymptotic formula for $P(\theta)$ the Laplace transform 
of a probability density function $p$ of a non-negative variable concentrated on $[p_0, p_1]$ having $N + 1$ continuous derivatives. By iteratively integrating by parts we find (cf Bleistein and Handelsman \cite{bleistein})

\begin{equation}
P(\theta) = {\bf P}_{p_0}^{(N)}(\theta) -  {\bf P}_{p_1}^{(N)}(\theta)  + 
					 {\bf \hat{P}}^{(N)}(\theta)
\end{equation}
where ${\bf P}_u^{(N)}$ and ${\bf \hat{P}}^{(N)}$ are defined as
\begin{equation}
\label{finiteExpansionOfLaplaceTransform}
{\bf P}_u^{(N)}(w) = \exp(-w u)\sum_{n = 0}^N p^{(n)}(u)w^{-n - 1}
\end{equation}

and
\begin{equation}
{\bf \hat{P}}^{(N)}(\theta) = \theta^{-N - 1}\int_{p_0}^{p_1}p^{(N + 1)}(u) \exp(-\theta u)du
\end{equation}

\begin{eqnarray}
{\bf \hat{P}}^{(N)}(\theta) &=&  \theta^{-N - 2} 
\left(p^{(N + 1)}(p_0)\exp(-p_0 \theta) - p^{(N + 1)}(p_1)\exp(-p_1 \theta) \right .
 \nonumber \\
 & & \left . + \int_{p_0}^{p_1}p^{(N + 2)}(u) \exp(-\theta u)du \right)
\end{eqnarray}

Moreover if $p$ has derivatives bounded by powers of $K$, ie, $|p^{(n)}(u)| \leq \kappa_0 K^n$ for $u \in[p_0, p_1]$
\begin{equation}
|{\bf \hat{P}}^{(N)}(\theta)| \leq \kappa_0 K^{N+1}   (2 + (p_1 - p_0) K) |\theta|^{-N-2}
            \exp(\max(-p_1 \theta_x, -p_0 \theta_x))
\end{equation}
thus for $|\theta| > K$ we have
$\lim_{N \rightarrow \infty} {\bf P}_{p_0}^{(N)}(\theta) - 
{\bf P}_{p_1}^{(N)}(\theta) = P(\theta)$.

\subsection{Decomposing $F$}
In the following section we are motivated by the case that both the service time distribution and the arrival time distributions have polynomial densities or are the result of a number of familiar operation on those densities such as convolution, order statistics or randomization. 

Moreover the analysis given below will still hold even if the arrival time distribution is the transform of a Cox phase distribution or one convolved with that of polynomial density.

For a large number of such cases it is possible to
decompose $F$ as the sum of $H$ and $G$ where $G$ is analytic except at zero and bounded in the left half plane and where
$H$ is an analytic except at zero and bounded in the right half plane with a limit of -1 as $x$ approaches infinity. 
In this section we present results that illustrate how $H$ can be decomposed at different levels facilitating its analysis and that of $\psi$. These decompositions allow us to apply Theorem \ref{prodWaitingThm} and the telescoping  methods of Section \ref{computationSchemaSection} to greatly accelerate the computation of average waiting time quantities of interest.

With the usual notation let $F(\theta) = A(-\theta)B(\theta) - 1$ where
$A(\theta)$ and $B(\theta)$ are Laplace transforms of $a$ and $b$, where
$a$ and $b$ are
probability densities with geometrically bounded derivatives concentrated on $[a_0, a_1]$ and $[b_0, b_1]$ 
respectively.
\begin{equation}
\label{RepresentationOfFAsInInfinityNeighborhood}
F(\theta) = -1 + \sum_{j = 0}^3 \exp(\alpha_j\theta)\Phi_j(\theta^{-1})\theta^{-k_j}
\end{equation}
where $\Phi_j$ is analytic around $0$ with $\Phi_j(0) \neq 0$, $k_j \geq 2$,
	and where $\alpha_3 = a_1 - b_0$, 
$\alpha_2 = \max(a_1 - b_1, a_0 - b_0)$, $\alpha_1 = \min(a_1 - b_1, a_0 - b_0)$, and
$\alpha_0 = a_0 - b_1$. $\alpha_0 < \alpha_1 \leq \alpha_2 < \alpha_3$. Observe $\alpha_3 > 0$ to satisfy the stability equation. Similarly $\alpha_0 < 0$ to prevent null waiting times. 
Define the partial sum approximations to $\Phi_j(w)$ by
\begin{equation}
\label{H_Approx_To_Psi}
\Phi_j^{(N)}(w) = \sum_{n = 0}^N \phi_n^{(j)} w^n
\end{equation}
and its tail by
\begin{equation}
\overline{\Phi}_j^{(N)}(w) = \Phi_j(w)
 - \Phi_j^{(N)}(w) = w^{N + 1} \sum_{n \geq 0} 
\phi_{n + N + 1}^{(j)} w^n = w^{N + 1} \hat{\Phi}_j^{(N)}(w)
\end{equation}
where
\begin{equation}
\label{PhiCaretTail}
\hat{\Phi}_j(N, w) =
\sum_{n \geq 0}\phi_{n + N + 1}^{(j)}w^n
\end{equation} 

Let $j_p$ be defined by $\alpha_{j_p} \leq 0 < \alpha_{j_p + 1}$. For technical reasons we need to add the assumption that $k_0 \le k_j$ for $j \le j_p$.

With this notation define the helping functions $H$ and $h$ by
\begin{equation}
\label{H_representation}
H(\theta) = -1 + \sum_{j = 0}^{j_p} \exp(\alpha_j\theta)
						\Phi_j^{(\kappa_p)}(\theta^{-1})\theta^{-k_j}
\end{equation}
where $\kappa_p = \min(k_{n_p + 1} ,.. k_3)$
and $h(\theta) = H(\theta) \theta^{m_p}$ where $m_p$ is the order of $H$'s pole at 0.

Observe that the core term of $H$ is given by
\begin{equation}
\label{TheCoreOf_F_EquationFor_T}
T(\theta) = c_0 \exp(\alpha_0\theta)\theta^{-k_0} -1
\end{equation}
 where $c_0 = \phi_0^{(0)}$ and that $H$ may be put into the form 
\begin{equation}
\label{H_ksiRepresentation}
H(\theta) = (T(\theta) + 1)(1 + \xi(\theta)) - 1
\end{equation}
where
$\xi(\theta) = \sum_{j_0\le j \le j_p} \exp(\beta_j)\xi_j(\theta^{-1})\theta^{-m_j}$, $m_j = k_j - k_0$, $\beta_j = \alpha_j - \alpha_0$, and
if $\Phi_0^{(\kappa_p)}$ is not a constant $j_0 = 0$  and $m_0 = 1$ otherwise $j_0 = 1$. 
$\xi_j(\theta) =\Phi_j^{(\kappa_p)}(\theta) c_0^{-1}$ for $j > 0$ and $\xi_0(\theta) =(\Phi_0^{(\kappa_p)}(\theta) -c_0)/c_0^{-1}\theta^{-1}$ for $j = 0$.

\begin{thm}
\label{HProductRepresentationThm}
\begin{equation}
\label{hFunctionProdRepresentation}
H(\theta) = \exp(\alpha\theta/2)q(\theta) 
						\prod_{n = n_0}^\infty(1-\theta/w_n)(1 - \theta/\overline{w}_n)
\end{equation}
where $q = q_0/q_1$ and $q_i$ are polynomials such that 
$\deg (q_0) - \deg (q_1) = 2 n_0 - 1$, where $\alpha = \alpha_0$, and where
$w_k = u_k + O(k^{-\lambda})$ and the $u_k$ are the complex roots of $T$ and $\lambda = \min_{j_0\le j \le j_p}(m_j + (\alpha_0 - \alpha_j)/\alpha_0)$ .
\end{thm}
{\bf Proof}
If $j_p = 0$ and $\Phi_0^{(\kappa_p)}$ is constant then $T = H$ and the theorem follows from Corollary
\ref{T_RepresCorollary}.
Observe that $H$'s representation given by equation \ref{H_ksiRepresentation} satisfies the conditions of Lemma \ref{BoldT_ExpansionLemma} from which the theorem now follows.

\begin{thm}
\label{psiRepThmRelatedToT}
Let $F$ and $T$ be respectively defined by equations \ref{RepresentationOfFAsInInfinityNeighborhood} and \ref{TheCoreOf_F_EquationFor_T} then the waiting time Laplace transform  $\psi$ has a product representation given by 
\begin{equation}
		\psi (\theta) =\exp\{-\alpha\theta/2\}
		\prod_{z\in {\mathcal Z}_0}(1-\theta/z)^{-1}
		\prod_{z\in {\mathcal Z}_1} (1-\frac{2\Re z_k\theta}{z\overline{z}}
		+\frac{\theta^{2}}{z\overline{z}})^{-1}
\end{equation}
where ${\mathcal Z}_0$ and ${\mathcal Z}_1$ are defined in Lemma \ref{psisZeroesRelatedToThoseOfT}.

Moreover $|u_n - z_n| = O(n^{-\lambda})$ where $\lambda > 0$ and the $u_n$ are zeroes appearing in equation \ref{littleTRepresentation}.
\end{thm}
{\bf Proof} The theorem will follow by proving that $F$ satisfies the conditions of Theorem \ref{prodWaitingThm}. Lemma \ref{psisZeroesRelatedToThoseOfT} conclusions cover all the premises of Theorem \ref{prodWaitingThm} except for the premise that $F$ is of exponential order $-\alpha$ at $-\infty$, and the set of zeroes in the strict left plane have not been mustered in the requisite form.

Since $F$ and $T$ satisfy the assertions of Lemma
\ref{psisZeroesRelatedToThoseOfT} we only need prove the missing assertions. 

Observe that as $T(z) +1$ is of exponential order $-\alpha$ at $-\infty$ and that 
\linebreak
$lim_{x\rightarrow -\infty} F(x)/(T(x) + 1) = 1$. These two limits imply that $F$ is of exponential order $-\alpha$ at $-\infty$.

Lemma \ref{psisZeroesRelatedToThoseOfT} shows that 
${\mathcal Z}$ the set of $F$'s zeroes in the strict left hand plane may be written as the disjoint union
of ${\mathcal Z}_0$ and ${\mathcal Z}_1$ where 
${\mathcal Z}_0 =  \{\zeta_n,  |n| < n_0\}$ where
$\zeta_n = \overline{\zeta}_{-n}$, and 
${\mathcal Z}_1 = \{z_n, \overline{z}_n, n\geq n_1\}$ where $ \Im{z_n} = 2m\pi\alpha^{-1}n +O(1)$ and $\Re{z_n} = O(n^b)$ for $b < 1$.

From ${\mathcal Z}_0$ and ${\mathcal Z}_1$ we will define a new sequence of zeroes that satisfies the premises of Theorem \ref{prodWaitingThm}.

Set $m = n_1 - n_0$, and for $0 \leq j < n_0$ define $v_j = \zeta_j$ while for
$j \geq n_0$ define $v_j = z_{j+m}$, and for $j< 0$ define $v_j = \overline{v}_{-j}$. 

Note that the asymptotic estimates for $v_j$ still hold true since the estimates are invariant under finite index shifts.

\begin{lem}
\label{psisZeroesRelatedToThoseOfT}
Let $F$ and $T$ be respectively defined by equations \ref{RepresentationOfFAsInInfinityNeighborhood} and \ref{TheCoreOf_F_EquationFor_T} then ${\mathcal Z}$, the zeroes of $F$ in the strict left plane, may be expressed as disjoint union of multisets
${\mathcal Z}_0$ and ${\mathcal Z}_1$, 
where ${\mathcal Z}_0 = \{\zeta_n, \overline{\zeta}_n, n < n_0\}$ and   
${\mathcal Z}_1 = \{z_n, \overline{z}_n, n\geq n_1\}$ such that $ \Im{z_n} = 2m\pi\alpha^{-1}n +O(1)$, $\Re{z_n} = O(n^b)$ for $b < 1$ and $|z_n -u_n| = O(n^{-\lambda})$ where $u_n$ are the zeroes of $T$ and $\lambda > 0$. Moreover let ${\mathcal E}_n(a)$ be the semi-circular curves as defined in Lemma \ref{LemmaForEstimatingRootsOf_fByT}, then $n \ge n_1$ implies  $|F| > 0.2$ on $ {\mathcal E}_n(0)$. Furthermore the radii $\{r_n\}$ of ${\mathcal E}_n(0)$ satisfy the asymptotic relationship $r_{n+1} = r_n + 2\pi/|\alpha| + o(1)$
\end{lem}
{\bf Proof}
We begin by showing that we may choose $a_1$, $a_2$, $N_1$ and $\lambda$ that satisfy the premises of Lemma \ref{LemmaForEstimatingRootsOf_fByT}.

Set $\xi(z) = (F(z) - T(z))/(T(z) + 1)$ then by identifying  $\Psi$ with  $\xi$ and $\gamma$ with $-m/\alpha$ in Lemma \ref{ExtendedFactorEstimatingLemma} we may infer that $\xi$ satisfies equation  
\ref{boundOnXiAroundTsZeroes} and the existence of an $a_1$ for which
$\Re{z} \leq a_1$ implies $|\xi(z)| < 0.1$.

Using easy asymptotic arguments one may show there exists an $N_1$ such that  $|z|\ge N_1$ and $a_1 \le \Re{z} \leq 0$ jointly imply $|F(z)| > 0.2$.

Identifying $F$ with $\Upsilon$ and setting $a_2 = 0$, we now see that 
$F$, $T$, $\xi$, $a_1$, $a_2$ and $N_1$ satisfy the premises of Lemma \ref{LemmaForEstimatingRootsOf_fByT}. From this we may deduce the existence of ${\mathcal Z}_1$ with the specified properties and the existence of semi-circular curves ${\mathcal E}_n(a_1)$ for $n \geq n_0$ such that $|F| > 0.2$ on  ${\mathcal E}_n(a_1)$ and that 
$z \in {\mathcal E}_n(a_1)$ implies $|z| \geq N_1$. From this we may deduce that $|F| > 0.2$ on  ${\mathcal E}_n(0)$. Also note that the asymptotic estimates for $u_n$ given in Corollary  \ref{horseshoeRegionForT} can be used to show $r_{n+1} = r_n + |2\pi/\alpha| + o(1)$.

Define ${\mathcal B}^- = \{z \ni |z| < r_{n_1}, \Re{z} < 0\}$. Since $F$ is an entire function it may vanish only a finite number of times on the bounded domain
${\mathcal B}^-$. Set ${\mathcal Z}_0$ to the zeroes of $F$
in  ${\mathcal B}^-$. Observe that since ${\mathcal B}^-$ is
symmetric relative to imaginary axis and that
$F(\overline{z}) = \overline{F(z)}$ it follows that
${\mathcal Z}_0$ is closed under conjugation.  Note that
$F$ is convex on the negative half-line with one negative real root $z_0$. Since ${\mathcal Z}_1$ has only complex zeroes, $z_0$ must belong to 
${\mathcal Z}_0$. Because ${\mathcal Z}_0$ is finite, closed under conjugation and has only one real root it must have an odd number of elements.

\begin{lem}
\label{sublinearityForHelpingFunH}
Let $H$ be defined by \ref{H_representation} then there exist $n_H$, $\delta_H > 0$, $\nu_H > 0$ and $\mu_H  > 0$
such that if $n \ge n_H$, $|\eta| \le \delta_H$ then
$\mu_H |\eta| \leq |H(w_n + \eta) - H(w_n)| \le \nu_H |\eta|$.
\end{lem}
{\bf Proof} Set ${\mathcal H} = H-T$ where $T$ is defined by \ref{TheCoreOf_F_EquationFor_T}. Since ${\mathcal H}$ satisfies the conditions of Corollary \ref{BoundingDeltaSumsOfTsNonDominantTerms} we may deduce the existence of positive constants $d_{\mathcal H}, \delta_{\mathcal H}, \nu_{\mathcal H},$ and $M_{\mathcal H}$ for which $|z| > M_{\mathcal H}$, $\Re{z} -  k_0\log(|z|)/\alpha_0 > 3|\log|c||$ and $|\eta| < \delta_{\mathcal H}$ jointly imply 
\begin{equation}
|{\mathcal H}(z+\eta)-{\mathcal H}(z)| < \nu_{\mathcal H}|\eta z^{-d_{\mathcal H}}|
\end{equation}

Observe by Corollary \ref{uniformityOfSemiLinearityOfT}
if $n > n_0$,
$|w_n - u_n| < \delta_1(\alpha_0)$ and
$|\eta| <\delta(\alpha_0)$ then
\begin{equation}
\label{extendedLinearIneqForTSect8}
0.75 \mu(\alpha_0) |\eta| \le |T(w_n + \eta) - T(w_n)| \leq 1.25 \nu(\alpha_0)|\eta|
\end{equation}
where  $n_0, \delta(\alpha_0), \delta_1(\alpha_0), \mu(\alpha_0),$ and $ \nu(\alpha_0)$ are positive constants.
 
By Prop \ref{continuityOfTsRootsRealPartEqn} we see that if
$|u_n| > 3\max(1/|\alpha_0|, 1) $ and 
$|u_n - w_n| < |\beta| $ where $\beta = -\log|c|$ then
$|\Re{w_n} -k_0\log|w_n|/\alpha_0  -\beta| < 2 |\beta|$ or that $\Re{w_n} -\log|w_n|/\alpha_0 > -3|\beta|$.

If  $\Re{w_n} - k_0\log(|w_n|)/\alpha_0 > 3|\beta|$ and $|w_n^{-d_{\mathcal H}}| < 0.25\mu(\alpha_0)/\nu_{\mathcal H}$
then $|{\mathcal H}(w_n+\eta)-{\mathcal H}(w_n)| < 0.25\mu(\alpha_0) |\eta|$.
From the last inequality and  \ref{extendedLinearIneqForTSect8} we deduce by the triangle inequality
\begin{equation}
0.50\mu(\alpha_0)|\eta| \le |H(w_n + \eta) - H(w_n)| \le 1.5 \nu(\alpha_0)|\eta|
\end{equation}
We can complete the proof by setting $\nu_H = 1.25\nu(\alpha_0)$, $\mu_H = 0.75\mu(\alpha_0)$,
$\delta_H = \delta(\alpha_0)$ and 
$n_H = \max(n_0, n_1)$ where $n_1$ guarantees that $n \ge n_1$ implies 
$|w_n| > \max(M_{\mathcal H}, 3, 3|\alpha_0|^{-1},
(0.25 \mu(\alpha_0)\nu(\alpha_0)^{-1})^{1/d_{\mathcal H}})$.

\bigskip
The next theorem in conjunction with Theorems \ref{psiRepThmRelatedToT} and \ref{hFunctionProdRepresentation} show why we can approximate statistics of $\psi$ distribution with related quantities of $H$. 
\begin{thm}
Let $F$ and $H$ be respectively defined by equations \ref{RepresentationOfFAsInInfinityNeighborhood} and \ref{H_representation}. Then we have the following asymptotic relationship between the strict left hand zeroes of $F$ and the zeroes of $H$: $|z_n - w_n| = O(n^{-\kappa_p - \lambda_p})$ where $\kappa_p = \min(k_{n_p + 1} ,.. k_3)$ and $\lambda_p = \min(1, -\alpha_{j_p + 1}k_0/\alpha_0)$.
\end{thm}
{\bf Proof} The proof proceeds by deploying Proposition \ref{bigOhTransitiveProp} to show that 
$|F-H|/|H + 1|$ is of order $O(|u_n|^{-\kappa_p - \lambda_p})$ in fixed size neighborhoods of the $u_n$. First we prove that 
$|H(\theta) -T(\theta)|/|T(\theta) + 1| = o(1)$, and then we show that $|H(\theta) -T(\theta)|/|T(\theta) + 1| = O(|\theta^{-\kappa_p - \lambda_p})$.

Employing Lemma \ref{ExtendedFactorEstimatingLemma} 
where we identify $\Psi(z)$ with $(H(z) -T(z))/(T(z) + 1)$ and 
$\gamma$ with $ -k_0/\alpha$ we may deduce  
that $|\theta| \rightarrow \infty$ in such a way that $x(\theta) = k_0\log|\theta|/\alpha_0 + O(1)$ implies
\begin{equation}
\label{HTDiffOverTPlus1}
|(H(\theta) -T(\theta))/(T(\theta) + 1)| = o(1)
\end{equation}
Set $F - H$ to $H_1 + G_1$ and define $H_1$ and $G_1$ by
\begin{equation}
H_1(\theta)= \sum_{j=0}^{j_p}\exp(\alpha_j\theta)\hat{\Phi}_j(\kappa_p, \theta^{-1}) \theta^{-k_j -\kappa_p -1}
\end{equation}
and by
\begin{equation}
G_1(\theta)= \sum\limits_{j=j_p + 1}^{3}\exp(\alpha_j \theta)\Phi_j(\theta^{-1})\theta^{-k_j}
\end{equation}
and where $\hat{\Phi}_j(\kappa_p, w)$ and
$\Phi_j(\theta)$ are respectively defined by \ref{PhiCaretTail} and \ref{H_Approx_To_Psi}.

\begin{equation}
\label{H1OverTPlus1Equation}
H_1(\theta)/(1 + T(\theta))= \sum_{j=0}^{j_p}\exp(\beta_j\theta)\hat{\Phi}_j(\kappa_p, \theta^{-1}) \theta^{-k_j + k_0 -\kappa_p -1}
\end{equation}

By applying Lemma \ref{ExtendedFactorEstimatingLemma} to \ref{H1OverTPlus1Equation} and noting that
$k_j \geq k_0$ and $\beta_j \geq 0$ we deduce
$|\theta|$ tending to $ \infty$ in such a way that
$x(\theta) = -\gamma\log|\theta| + O(1)$ implies
\begin{equation}
\label{H1_TPlus1RatioBigOEstimate}
 |H_1(\theta)/(1 + T(\theta))| = 
O(\theta^{-\lambda(\hat{\beta}, \hat{k})}) =
O(\theta^{ -\kappa_p - 1})
\end{equation} 
where 
$\lambda(\hat{\beta}, \hat{k}) = 
\min_{0\leq j\leq j_p}
\{\gamma\beta_j + k_j - k_0 + \kappa_p + 1\} \geq \kappa_p +1$.

\begin{equation}
\label{G1OverTPlus1Equation}
G_1(\theta)/(1 + T(\theta))= \sum\limits_{j=j_p + 1}^{3}\exp(\beta_j \theta)\Phi_j(\theta^{-1})\theta^{-k_j + k_0}
\end{equation}

By applying Lemma \ref{ExtendedFactorEstimatingLemma} to \ref{G1OverTPlus1Equation} and noting that
$k_j \geq k_0$ and $\beta_j > -\alpha_0$ we can infer that 
$|\theta| \rightarrow \infty$ in such a way that
$x(\theta) = -\gamma\log|\theta| + O(1)$ implies
\begin{equation}
\label{G1_TPlus1RatioBigOEstimate}
 |G_1(\theta)/(1 + T(\theta))| = 
O(\theta^{-\mu(\hat{\beta}, \hat{k})}) =
O(\theta^{-\kappa_p + \alpha_{j_p + 1}k_0/\alpha_0})
\end{equation} 
where 
$\mu(\hat{\beta}, \hat{k}) = 
\min_{0\leq j\leq j_p}
\{\gamma\beta_j + k_j - k_0\} 
\geq \kappa_p -k_0 + (\alpha_{j_p +1} -\alpha_0)\gamma $. Noting that $\gamma\alpha_0 = -k_0$ we see that $\mu(\hat{\beta}, \hat{k}) \geq \kappa_p - \alpha_{j_p + 1}k_0/\alpha_0 > \kappa_p$.

From equations \ref{H1_TPlus1RatioBigOEstimate} and
\ref{G1_TPlus1RatioBigOEstimate} we see that
\begin{equation}
\label{FHDiffOverTPlus1}
|(F(\theta) - H(\theta))/(T(\theta) + 1)| 
= O(\theta^{-\kappa_p -\lambda_p})
\end{equation}
 where 
$\lambda_p = \min(1, -\alpha_{j_p + 1}k_0/\alpha_0)$.

Proposition \ref{bigOhTransitiveProp} shows that equations \ref{HTDiffOverTPlus1} and \ref{FHDiffOverTPlus1} jointly imply that
\begin{equation}
|(F(\theta)-H(\theta))/(H(\theta) +1)| = O(|\theta|^{-\kappa_p -\lambda_p})
\end{equation}

Lemma \ref{sublinearityForHelpingFunH} indicates that
$H$ and $F-H$ satisfy the conditions of Lemma 
\ref{zeroApproxLemma} which implies the conclusion of this theorem.

\subsection{An illustration of the function $H$}
Let the interarrival time be the constant $1/2$ and service time distribution defined on $[0,1]$ with density $b(x) = (22 - 6x - 12x^2)/15$. In this case

$B(\theta) = 22\theta^{-1} -6\theta^{-2} - 24\theta^{-3}
-\exp(-\theta)\left(4\theta^{-1} -30\theta^{-2} - 24 \theta^{-3} \right)$.

Then $F(\theta) = H(\theta) + G(\theta)$ where

$H(\theta) = \exp(-\theta/2)(-4\theta^{-1} +30\theta^{-2} + 24 \theta^{-3} ) -1$
and

$G(\theta) = \exp(\theta/2)( 22\theta^{-1} -6\theta^{-2} - 24\theta^{-3}) $.

As $\Re{\theta} \rightarrow -\infty $ it is easily shown that $G(\theta) \rightarrow 0$. Moreover we can also see that the $|4\exp(-\theta/2)\theta^{-1}|$ dominates $H(\theta)$ as  $\Re{\theta} \rightarrow -\infty $ which leads to $T(\theta) = -4\exp(-\theta/2)\theta^{-1} -1$.

\mysection{Comparison Functions}
In this section we study the function $\sigma$, defined by \ref{tDefEqn} and three other closely related functions $S(\theta) = \sigma(\theta)\theta^{-m}$, $t(\theta) = c\exp(\alpha \theta)  - \theta^m$, and $T(\theta) = t(\theta)\theta^{-m}$. The latter function $T$ provides the estimates needed to prove
the expansion given by \ref{partFracResult}. Note that $T$ is linearly equivalent to $F$ for $D/E_k/1$ and $E_k/D/1$ queues (cf Proposition \ref{E_D_1_ToTeeTransform}), and forms the dominant term for a number of other $F$. Moreover, through $T$ one may compute the zeroes of $F$ corresponding to $D/E_k/1$ and $E_k/D/1$ queues via straightforward real valued computations (cf \ref{hEquationForT}).

\begin{equation}
\sigma(\theta) = \exp(\theta + \beta) -  \theta^m
\label{tDefEqn}
\end{equation}
where $\beta$ is real and $m$ is a positive integer. 

\begin{prop}
Define $f(\theta) = \exp(c\theta)(1 + \theta/b)^{-m} - 1$
where $b, c$ are real numbers such that $bc > 0$.
Then $f(\theta) = \sigma(z)z^{-m}$
where $z =  c(b + \theta)$ and $\beta = m \ln(bc) -bc$ .
\label{E_D_1_ToTeeTransform}
\end{prop}
\begin{thm}
\label{CompareFuncExp}
\begin{equation}
\label{tauExpEqn}
\sigma(\theta) = \exp(\theta/2)p(\theta)
				\prod_{n=1}^\infty (1-\theta/z_n)(1-\theta/\overline{z}_n)
\end{equation}
where $p$ is a polynomial of degree $m + 1$ whose roots occur in complex conjugates  except for at most 3 real roots, and where
$\Im{z_n} = (2n + m/2)\pi + o(n^{\eta -1})$ and
$\Re{z_n} = m\ln(2n\pi) + 4^{-1}m^2n^{-1} -\beta + o(n^{\eta -2})$. Moreover
\begin{equation}
\label{znArgRelationship}
\Im{z_n} = 2n\pi + m\arg(z_n)
\end{equation}
\end{thm}

{\bf Proof} One may show that $|\sigma(\theta)| \leq C \exp(|\theta|)$, where
$C = \exp(\beta) + m^m\exp(-m)$. It is a straightforward computation to show that $\sigma$ has exponential order 0 at $-\infty$. Since the zeroes of $\sigma$ satisfy the hypotheses of Theorem \ref{inf_prod_characterization} as shown by Corollary \ref{zeroesOfTauLemma},
Theorem \ref{inf_prod_characterization}  implies that $\sigma$ has a representation given equation \ref{tauExpEqn}. The asymptotic estimates for $z_n$ and
\ref{znArgRelationship} are implied by Lemma \ref{strongLemRhoSubnAsymp}.

\begin{cor}
\label{T_RepresCorollary}
Let $t(z) = c\exp(\alpha z)  - z^m$ where $\alpha < 0$, $c\alpha^m> 0$ and $m$ is a positive integer. Then
\begin{equation}
\label{littleTRepresentation}
t(z) = \exp(\alpha z/2)\xi(z)\prod_{n = 1}^\infty(1 - z/u_n)(1-z/\overline{u}_n)
\end{equation}
where $\xi$ is a polynomial of degree $m + 1$ whose roots occur in complex conjugates  except for at most 3 real roots, and where $\Im{u_n} = (2n\pi + m\arg(-u_n))/\alpha = (2n + m/2)\pi\alpha^{-1} + o(n^{\eta -1})$,
$\Re{u_n} = (m\ln(2n\pi) + 4^{-1}m^2n^{-1} -\beta)\alpha^{-1} + o(n^{\eta -2})$, and
$\beta = m\ln|\alpha| +\ln|c|$. 

Furthermore $u$ is a root of $t$ if and only if there exists a unique integer $n_u$ such that
\begin{equation}
\label{characterizationOfImagOftz}
\Im{u}   = (2n_u\pi + m\arg(u) -\arg(c))/\alpha
\end{equation}
and
\begin{equation}
\label{characterizationOfRealOftz}
\Re{u}  =( m\log|u| -\log|c|) /\alpha
\end{equation}
where $\arg(z)$ is the principal value of the argument of $z$.
\end{cor}
{\bf Proof} Equations \ref{characterizationOfImagOftz} and \ref{characterizationOfRealOftz} may be proved by observing that if $t(u)=0$ if and only if 
\begin{equation}
\label{tEquationForm}
c\exp(\alpha u) = u^m
\end{equation}
Taking the argument of each side of \ref{tEquationForm} yields $\arg(c) + \alpha\Im{u} = m\arg(u) +2k\pi$, for some integer $k$, but since $u$, $\alpha$ and $c$ are all fixed $k$ must be unique. Computing the logarithm of the moduli of each side of \ref{tEquationForm} and solving for $\Re{u}$ completes the proof.

Substituting $\theta = \alpha z$ and $\beta = m\ln|\alpha| +\ln|c|$
into $\sigma(\theta)$ shows that $\sigma(\theta) = \alpha^m t(z)$. The rest of the corollary follows from Theorem \ref{CompareFuncExp}.
\begin{cor}
Let $\tau(z) = c\exp(\alpha z)  +  z^m$ where $\alpha < 0$, $c\alpha^m > 0$ and $m$ is a positive integer. Then
\begin{equation}
\tau(z) = \exp(\alpha z/2)\zeta(z)\prod_{n = 1}^\infty(1 - z/w_n)(1-z/\overline{w}_n) =
\exp(\alpha z/2)\zeta(z)\gamma_0(z)
\end{equation}
where $\zeta$ is a polynomial of degree $m$ whose roots occur in complex conjugates  except for at most 1 real root, and where $\Im{w_n} = (2n + 1 + m/2)\pi\alpha^{-1} + o(n^{\eta -1})$,
$\Re{w_n} = (m\ln((2n - 1)\pi) + 2^{-1}m^2 (2n - 1)^{-1} -\beta)\alpha^{-1} + o(n^{\eta -2})$, and
$\beta = m\ln|\alpha| +\ln|c|$. 
\end{cor}
{\bf Proof} One could recapitulate the detailed arguments for proving the expansion of $t(z)$ in slightly modified form to derive the above expansion for $\tau$ but a cleaner way is to leverage the known expansion for $t$ with a bit of algebra.

Observe
\begin{equation}
t_{c^2,2\alpha,2m}(z) = t_{c,\alpha, m}(z)\tau(z)
\end{equation}
where $t_{d,a,j}(z) = d\exp(a z)  -  z^j$. 

As per Corollary \ref{T_RepresCorollary} we introduce the following expansions
\begin{equation}
t_{c^2,2\alpha,2m}(z)\exp(-\alpha z) = \xi_2(z)\prod_{n = 1}^\infty(1 - z/u_{2,n})(1-z/\overline{u}_{2,n}) =
\xi_2(z)\gamma_2(z)
\end{equation}

\begin{equation}
t_{c,\alpha,m}(z) = \exp(\alpha z/2)\xi_1(z)\prod_{n = 1}^\infty(1 - z/u_{1,n})(1-z/\overline{u}_{1,n})
\end{equation}

As per equations \ref{characterizationOfImagOftz} and \ref{characterizationOfRealOftz} we see that if $u$ is a root of $t_{c^2,2\alpha,2m}$ it follows that there exists $n_u$ for which
\begin{equation}
\Im{u} = \left (2n_u\pi + 2m\arg(u) -\arg(c^2)\right)/2\alpha = \left(n_u\pi + m\arg(u)  -\arg(c) \right)/\alpha
\label{argEqnForZeroesOfTSquare}
\end{equation}
and that
\begin{equation}
\Re{u} = \left (2m\log|u| - \log(c^2)\right)/2\alpha = 
\left (m\log|u| - \log|c|\right)/\alpha 
\end{equation}
The two previous equations imply $u$ is a root of  $ t_{c,\alpha, m}$ for $n_u$ even, and for odd $n_u$   simple algebra shows that $u$ is a root of $\tau$.

In particular if we select $u_{2,n}$ for $u$ we find that 
$\Im u_{2,n} = (n\pi + m\arg(u_{2,n}) - \arg(c))/\alpha$. According to Corollary \ref{T_RepresCorollary} if $n$ is even then $u_{2,n}= u_{1, n/2}$ otherwise we have $u_{2,n}$ is a zero of $\tau$. If we define $w_n$ to be
$u_{2,2n - 1}$ and we see that $\gamma_2(z) = \gamma_1(z)
\gamma_0(z)$. By cancellation we can deduce that $\zeta(z)= \xi_2(z)/\xi_1(z)$ which proves that $\zeta$ is a polynomial of degree $m$.

The asymptotic estimates for $w_n$ follow from those of $u_{2,n}$ given by Corollary \ref{T_RepresCorollary} and some algebraic simplifications.

The claim that $\tau$ has at most one real root may be inferred by applying $c\alpha^m > 0$ to equation \ref{argEqnForZeroesOfTSquare}.

\begin{prop}
Let $c$ and $\alpha$ be non-zero reals and $m$ be integer. Suppose $z$ is non-zero complex number such that
\begin{equation}
\label{genericEquationForZxForTeeNTau}
\alpha\Re{z} = m \log|z| - \log|c|
\end{equation}
\begin{equation}
\label{genericEquationForZyForTeeNTau}
\alpha\Im{z} = m \arg(z) + k\pi
\end{equation}
where $k$ is integer.
If $\arg(c) + k\pi = 0 \mod{2\pi}$ then
$c\exp(\alpha z) -z^m = 0$ else
$c\exp(\alpha z) +z^m = 0$.
\end{prop}
{\bf Proof}
The conclusion of this proposition is equivalent to $|c\exp(\alpha z)| = |z|^m$ and $\arg(c\exp(\alpha z)z^{-m}) = \arg(1)$ for  $\arg(c) + k\pi = 0 \mod{2\pi}$ and  otherwise $\arg(c\exp(\alpha z)z^{-m}) = \arg(-1)$.

Notice that the first equation implies that 
$|c\exp(\alpha z)| = |z|^m$
while the second implies $\arg(c\exp(\alpha z)z^{-m}) = \arg(c) + k\pi \mod{2\pi}$. If $\arg(c) + k\pi = 0 \mod{2\pi}$ we may deduce $\arg(c\exp(\alpha z)z^{-m}) = \arg(1)$ and otherwise $\arg(c) + k\pi = \pi \mod{2\pi}$ implying that
 $\arg(c\exp(\alpha z)z^{-m}) = \arg(-1)$. 

\subsection{Analysis of a $T$ Extension}
We introduce ${\bf T}$ an extended form of $T$ that figures prominently in providing sharp approximations.

\noindent
${\bf T}(z) = c\exp(\alpha z)z^{-m}(1 + \xi(z)) -1$ where $m > 0$, $\xi(z) = \sum_{0\le j \le j_p}\exp(\beta_i z)\xi_j(z^{-1})z^{-m_j}$ and $0 \leq \beta_j \leq -\alpha$, $m_j - \beta_j /\alpha > 0$, and  $c\alpha^m > 0$, and the $\xi_j$ are polynomials nonvanishing at zero.

\begin{lem}
\label{BoldT_ExpansionLemma}
Let ${\bf T}$ be defined as above then the following expansion holds
\begin{equation}
{\bf T}(z) = \exp(\alpha z/2){\bf q}(z)\prod_{n = n_0}^{\infty} (1 - z/w_n)(1 - z/\bar{w}_n)
\end{equation}
\noindent
where ${\bf q}$ is a rational function whose numerator is a polynomial of degree 
$2n_0 - 1 +  M$ and denominator is a monomial of degree $M$, where $M$ is the order of the pole of {\bf T} at 0,
and where $z_n = w_n + O(n^{-\lambda})$ with  $\lambda = \min_{0\le j \le j_p}(m_j - \beta_j/\alpha) > 0$.
\end{lem}
{\bf Proof}
Set ${\bf t(z)} = {\bf T}(z)z^M$. We begin by showing that ${\bf t}$ satisfies the conditions of 
Corollary \ref{inf_prod_characterization_dual}. $M$ has been selected so that ${\bf t}(0) \neq 0$.
Observe that ${\bf t}$ is an analytic function for which there exists a constant $C$ such that $|{\bf t}(z)| < C \exp(2\alpha |z|)$. Moreover it may be shown that {\bf t}'s exponential limit at $\infty$ is $0$.
 
From Lemma \ref{tripartiteTboldTRegions} we deduce that ${\bf t}$'s zeroes may be partitioned into the requisite disjoint sets $W_b$ and
$W_d$ where $W_b$ is finite and $W_d =\{n \geq n_0, w_n, \overline{w}_n\}$ with  $\Im{w_n} = 2\pi\alpha^{-1}n +O(1)$ and
        $\Re{w_n} = O(n^a)$ for $a < 1$.
Thus ${\bf t}$ satisfies the conditions of Corollary \ref{inf_prod_characterization_dual} implying that {\bf t} has the following representation.
\begin{equation}
{\bf t}(z) = \exp(\alpha z/2) {\bf t}_1(z) 
\gamma_{\bf t}(z)
\end{equation}
where ${\bf t}_1(z) = {\bf t}(0)\prod_{w \in W_b} (1 - z/w)$ and

$\gamma_{\bf t}(z) = \prod_{n \geq n_0} (1 - z/w_n)(1 - z/\bar{w}_n)$.

The desired representation for ${\bf T}$ obtained by dividing the last equation by $z^m$ and setting $q(z) = {\bf t}_1(z)z^{-M}$

Observe by leveraging the expansion for $T(z)z^m$ given by equation \ref{littleTRepresentation} we may write $T$ as
\begin{equation}
T(z) = \exp(\alpha z/2)
	   q_t(z)\prod_{n \geq n_0}
	   (1 - z/u_n)(1 - z/\bar{u}_n)
\end{equation}
where

$q_t(z) = z^{-m}\xi(z)\xi_1(z)$
where $\xi_1(z) = \prod_{n=1}^{n = n_0 -1}(1 - z/u_n)(1 - z/\bar{u}_n)$.

where $\xi$ is polynomial of degree $m + 1$

We may now apply Corollary \ref{degCompForInfProd} to the representations of $T$ and ${\bf T}$ to 
deduce 
\begin{equation}
  \deg ({\bf t}_1(z) z^m) = \deg(\xi(z)\xi_1(z)z^M)
\end{equation}
Since $\deg(\xi(z) \xi_1(z)z^M) = 2n_0 - 1 + M + m$ it follows the last equation that $\deg({\bf t}_1(z)) = 2n_0 - 1 + M$.
\vspace{0.20in}

The following lemma facilitates the analysis of ${\bf T}$'s zeroes by dividing the plane into three regions, defined by
${\mathcal D}(a,N) = \{z: \Re{z}\leq a$ 
and $|z| \geq N \} $,
 ${\mathcal B}(N) = \{z: |z| < N \}$ and
 ${\mathcal C}(a, b,N) = \{z: a \le \Re{z} < b$ and $|z| \geq N \}$.
\begin{lem}
\label{tripartiteTboldTRegions}
Let $T$, ${\bf T}$ and $M$ be as defined in Lemma
\ref{BoldT_ExpansionLemma}, then we may choose $a$ and $N$ such that ${\bf T}$ is non-vanishing on 
${\mathcal C}(a,\infty, N)$, 
${\bf T}$ vanishes only on a finite number of points of
${\mathcal B}(N)$, and on  
 ${\mathcal D}(a,N)$  the zeroes of $T$ and {\bf T} may be put into a bijective relationship. In particular if we denote the zeroes of $T$ and those of {\bf T} on  ${\mathcal D}(a,N)$ by $U_d$ and $W_d$ respectively, then $U_d$ and $W_d$ may be expressed as $U_d =\{u_n, \overline{u}_n, n \ge n_1\}$ and $W_d =\{w_n, \overline{w}_n, n \geq n_1\}$ where the $u_n$ are the zeroes of $T$ as given by Corollary \ref{T_RepresCorollary}. Moreover $|u_n - w_n| = O(n^{-\lambda})$ with $\lambda$ as defined by Lemma \ref{BoldT_ExpansionLemma}. Finally the $\Im{w_n} = 2\pi\alpha^{-1}n +O(1)$ and
        $\Re{w_n} = O(n^b)$ for $b < 1$
\end{lem}
{\bf Proof} We begin by observing that since ${\bf T}(z)z^M$ is an entire function it vanishes at most a finite number of times on the bounded 
domain defined by ${\mathcal B}(R)$ for any finite $R$. 

The rest of proof consists of showing that we can find $a_1$, $a_2$, $N_1$ and $\lambda$ that satisfy the premises of Lemma \ref{LemmaForEstimatingRootsOf_fByT}.

Applying Lemma \ref{ExtendedFactorEstimatingLemma} with $\gamma$ set to $m/\alpha$ shows that 
$|\xi(z)|$ satisfies equation \ref{boundOnXiAroundTsZeroes}, and that we may choose $a_1$ for which $\Re{z} \le a_1$ implies $|\xi(z)| < 0.1$ satisfying equation \ref{uniformRightBoundOnXi}. 

Observe if  $ \Re{z}\geq a_1$ and $|z| \geq 1$ then $|{\bf T}(z) + 1|\leq K(a_1)|z^{-1}|$ where $K(a_1) = |c|\exp(|\alpha a_1|)(1 + K^\prime)$ and $K^\prime$ is equal to the sum of the absolute values of each of the coefficients of $\xi_0$, $\xi_1$
and $\xi_2$. Set $N = 4 K(a_1)$ so that $|z| \ge N$ and $\Re{z} \geq a_1$ jointly imply $|T(z) + 1| \leq 0.25$. We see by the triangle inequality that for $z \in {\mathcal C}(a_1,\infty, N)$ that
 $|{\bf T}(z)| \geq 0.75$.

The remaining conclusion of this lemma now follow from those of Lemma \ref{LemmaForEstimatingRootsOf_fByT}.

\begin{lem}
\label{LemmaForEstimatingRootsOf_fByT}
Let $T$ be as defined in Lemma
\ref{BoldT_ExpansionLemma} and let $\Upsilon(z)z^M$ be an entire function, where $M$ 
is a nonnegative integer. Define $\xi(z)$ as $(\Upsilon(z) - T(z))\exp(-\alpha z)z^m$. 
Suppose there exist $a_1$,  $a_2$, $\lambda > 0$ and $N_1$ for which the following claims hold.

If $\Re{z} \le a_1$  then
\begin{equation}
\label{uniformRightBoundOnXi}
|\xi(z)| < 0.1
\end{equation}
If $|z|$ tends to $\infty$ in such a way that
$x(z) = m\log|z|/\alpha + O(1)$ then
\begin{equation}
\label{boundOnXiAroundTsZeroes}
|\xi(z)| = O(|z|^{-\lambda|})
\end{equation}
and if $a_1 \leq \Re{z} < a_2$ and $|z| \geq N_1$ then $|F(z)| > 0.2$.

With these assumptions the following conclusions hold.
We may choose  $N \geq N_1$ such that on
${\mathcal D}(a_1,N) = \{z: \Re{z}\leq a_1$ 
and $|z| \geq N \} $  the zeroes of $T$ and $\Upsilon$ may be put into a bijective relationship. In particular if we denote the zeroes of $T$ and those of $\Upsilon$ on  ${\mathcal D}(a_1,N)$ by $U_d$ and $W_d$ respectively, then $U_d$ and $W_d$ may be expressed as $U_d =\{u_n, \overline{u}_n, n \ge n_1\}$ and $W_d =\{w_n, \overline{w}_n, n \geq n_1\}$ where the $u_n$ are the zeroes of $T$ as given by Corollary \ref{horseshoeRegionForT}.  
Define the semi-circular curves ${\mathcal C}_n^\prime(a_1) = \{|z| = r_n$ and $ \Re{z}\leq a_1\}$ where $r_n = (|u_n|+ |u_{n+1})/2$. Then there exists $n_1$ for which $n \geq n_1$ implies $\Upsilon$ is bounded below by 0.20 on  ${\mathcal C}_n^\prime(a_1)$.
Moreover the following asymptotic results hold $|u_n - w_n| = O(n^{-\lambda})$, $\Im{w_n} = 2m\pi\alpha^{-1}n +O(1)$ and
        $\Re{w_n} = O(n^b)$ for $b < 1$
\end{lem}
{\bf Proof}
   Corollary \ref{horseshoeRegionForT} shows the  existence of $n_1$ for which  $r_{n_1}\geq N_1$ and a sequence of horseshoe shaped regions ${\mathcal D}_n^\prime(a_1)$ defined for $n \geq n_1$ such that $|T|$ is bounded below by $1/2$ on perimeter of ${\mathcal D}_n^\prime(a_1)$ and $T$ has exactly two zeroes inside of ${\mathcal D}_n^\prime(a_1)$. Moreover by applying Proposition \ref{roucheIneqProp} with $\delta$ set to 0.1 shows that $\Upsilon > 0.2$ on ${\mathcal D}_n^\prime(a_1)$. Since ${\mathcal C}_n^\prime(a_1)
   \subset {\mathcal D}_n^\prime(a_1)$ we have $\Upsilon > 0.2$ on ${\mathcal C}_n^\prime(a_1)$ as well.
    
   Applying Lemma \ref{rouchForTFunctions} to $T$  and  $\Upsilon$ on the region defined by ${\mathcal D}_n^\prime(a_1)$ and its interior we find that $\Upsilon$ must also have exactly two roots in the interior of this region. It is convenient to introduce the half open, half closed disjoint regions
   ${\mathcal E}_n(a_1) = \{r_n \leq |z| < r_{n+1}$ and $\Re{z} \leq a_1\}$ and their partial unions
    ${\mathcal E}_{n, m}(a_1) =  
    \{r_n \leq |z| < r_{n +m+1}$ 
      and $\Re{z} \leq a_1\}$.
   Observing that
$\lim_{m\rightarrow \infty}{\mathcal E}_{n_1,m}(a_1) = {\mathcal D}(a_1,r_{n_1})$  provides a means for inductively constructing a bijective mapping on  ${\mathcal D}(a,N)$ of the zeroes of $T$ to the zeroes of $\Upsilon$ where 
$N = r_{n_1}$.

Observe if $|w_n - u_n| = O(n^{-\lambda})$ then  the asymptotic estimates for $\Re{w_n}$ and $\Im{w_n}$ follow from those for $\Re{u_n}$ and $\Im{u_n}$ given in Corollary \ref{horseshoeRegionForT}.

Corollary \ref{semiLinearityOfT} implies that for $n_1$ sufficiently large there exist $\mu > 0$, $\nu > 0$ and $\delta > 0$ such that if $n \ge n_1$ and $|\eta| < \delta$ we have 
\begin{equation}
\mu|\eta| < |T(u_n + \eta)| < \nu|\eta|
\end{equation}

Corollary \ref{T_RepresCorollary} yields
$u_n = m\ln(n)/\alpha + O(1)$ and thus by equation \ref{boundOnXiAroundTsZeroes} we may deduce
$|\xi(u_n + \eta)| = O(n^{-\lambda})$ for $|\eta| <  \delta$. 

Lemma \ref{zeroApproxLemma} may be applied to $T$ and $u_n$ to deduce $|w_n - u_n| = O(n^{-\lambda})$.

\subsection {Analysis of $\sigma$'s Roots}
To simplify the analysis of the roots we introduce the multisets where the multiplicity of each element $b$ is equal to $b$'s multiplicity as a root of $\sigma$. Define
${\mathcal Z}$ as $\{z \ni \sigma(z) = 0\}$  and define
${\mathcal U}$ as $\{z \ni \Im{z} \geq 0$ and $\sigma(z) = 0\}$.
Since $\sigma(\overline{z}) = \overline{\sigma(z)}$ and similarly for all of $\sigma$'s derivatives it follows that a root $b$ has the same multiplicity as its conjugate, hence ${\mathcal Z} = {\mathcal U} \cup \overline{{\mathcal U}} $.
Thus it suffices to characterize the elements of ${\mathcal U}$. Moreover according to Lemma \ref{tizeroes} there is at most one zero of $\sigma$ with a multiplicity of two and none higher.

The equation $\sigma(z) = 0$ with $\Im{z} \geq 0$ is equivalent to the following set of simultaneous equations in the real variable $\rho$ which is the main focus of the rest of this section.

\begin{eqnarray}
x(\rho) &=& m \ln(\rho) - \beta \\
\label{xRhoEquationForT}
|x(\rho)|&\leq& \rho \\
\omega(\rho) &=& \arccos(x(\rho)/\rho) \\
y(\rho) &=& (\rho^2 - x(\rho)^2)^{1/2} \\
h(\rho) &=& 2n\pi
\label{hEquationForT}
\end{eqnarray}
where $h(\rho) = y(\rho) -m\omega(\rho)$,  $n$ is an integer, and  
the square root and arccos are taken to be the principal branches.

The constants $r_0, r_1, r_2$ and set $A_{m, \beta}$ are defined in Lemma
\ref{xRhoCharacterization}, and are used in the remainder of this subsection.
\begin{lem}
The function $\sigma$ has at most one zero of order 2 and none higher. Moreover
$\sigma$ has zero of order 2 at $\theta = m$
only if $\sigma(m) = 0$.
\label{tizeroes}
\end{lem}
{\bf Proof} Differentiate $\sigma$ to deduce the following equations.
\begin{equation}
\sigma^\prime(\theta) = \sigma(\theta) - \theta^{m-1}(m-\theta)
\label{tiprime}
\end{equation}
\begin{equation}
\sigma^{\prime\prime}(\theta) = \sigma(\theta) - \theta^{m-2}(m^2 -m -\theta^2)
\label{tiprime2}
\end{equation}
Equation \ref{tiprime} indicates that $\theta = m$
is the only possible solution to the simultaneous system of
equations $\sigma(\theta) = 0$ and $\sigma^\prime(\theta) = 0 $.
Equation \ref{tiprime2} shows that if $\sigma(m) = 0$ then $\sigma^{\prime\prime}(m) =  m^{m-1}\neq 0$.

\begin{lem} 
The function $h$ maps $A_{m, \beta}$ onto $[-m\pi, \infty)$.
If $m -x(m) > 0$ then the mapping is strictly increasing, otherwise
it is strictly increasing on $[r_0, r_1]$ and on $[r_2, \infty)$ with $h(r_1) = h(r_2) = 0$.
\label{hRangeLemma}
\end{lem}
{\bf Proof} The lemma is jointly implied by Lemmas \ref{hMonoLemma} and  \ref{xRhoCharacterization} .
\begin{cor}
\label{hOfOmegaMonotonicity}
For every integer $n \geq -m/2$ and $n \neq 0$ there exists a unique $\rho_n > 0$ such that
$h(\rho_n) = 2n\pi$. If $m \geq x(m)$ then there exists a unique $\rho_0 > 0$ such that $h(\rho_0) = 0$, otherwise there exist $\rho_0 < \hat{\rho}_0 $ such that 
$h(\rho_0) = h(\hat{\rho}_0) = 0$. Moreover the $\rho_i$ are strictly monotonically increasing in $i$, and if $m < x(m)$ then $\hat{\rho_0} < \rho_1$.
\end{cor}

\begin{cor}
\label{zeroesOfTauLemma}
The multiset of $\sigma$'s zeroes can be partitioned into  
multisets $ = {\mathcal  Z}_0$ and ${\mathcal  Z}_1$ that satisfy the hypotheses of Theorem
\ref{inf_prod_characterization}, and the elements of ${\mathcal  Z}_1$ satisfy equation \ref{znArgRelationship} 
\end{cor}
{\bf Proof} Setting $z(\rho) = x(\rho) + iy(\rho)$ we may define
 ${\mathcal  Z}_1 = \{z(\rho_n), \overline{z(\rho_n)}, n > 0\}$. Notice that the definition of $\rho_n$ is equivalent to equation \ref{znArgRelationship} since $\arg(z_n) = \arccos(\Re{z_n}/|z_n|) = \arccos(x(\rho_n)/\rho_n)$.
 
Let ${\mathcal  Z}_0 = \{w_0, w_1, ... w_m\}$, where $w_i$ are defined as follows:
Define $w_0 = z(\rho_0)$.
If $x(m) < m$ define  $w_1 = z(\hat{\rho}_0)$
else set $w_1 =  \overline{z(\rho_0)}$.
For $0 < j < m/2$
define $w_{2j} = z(\rho_{-j})$ and $w_{2j + 1} = \overline{z(\rho_{-j})}$.
For $m$ even define $w_m = x(\rho_{-m / 2})$.
Observe that ${\mathcal  Z}_0$ and ${\mathcal  Z}_1$ are closed under conjugation and that ${\mathcal  Z}_0 \cup {\mathcal  Z}_1 \supset {\mathcal  U}$. Therefore by the opening remark of this subsection it follows that ${\mathcal  Z} ={\mathcal  Z}_0 \cup {\mathcal  Z}_1$

Lemma \ref{tizeroes} shows that $\sigma$'s zeroes of  are simple when
$x(m) \neq m$, and otherwise if $x(m) = m$ then  $w_0 = w_1$ and all of $\sigma$'s remaining zeroes are simple.
Lemma \ref{strong_y_rho_n_asymp} shows that for any $\eta > 0$
\begin{equation}
y_n = (2n + m/2)\pi + o(n^{\eta - 1}) 
\end{equation}
\begin{equation}
x_n = m \ln(2n\pi) + 4^{-1}m^2 n^{-1} -\beta + o(n^{\eta - 2}) = O(n^\eta)
\end{equation}
\begin{lem}
The function $h^\prime$ is strictly positive on the interior of $A_{m, \beta}$.
\label{hMonoLemma}
\end{lem}
{\bf Proof}
\begin{equation}
\omega^\prime(\rho) = \arccos(x(\rho)\rho^{-1})^\prime = (x(\rho) -m)\rho^{-1} y(\rho)^{-1}
\end{equation}
\begin{equation}
y^\prime(\rho) = (\rho^2 -mx(\rho))\rho^{-1} y(\rho)^{-1}
\end{equation}
\begin{equation}
h^\prime(\rho) = y^\prime(\rho) - m\omega^\prime(\rho) = (\rho^2 -2mx(\rho) + m^2)(\rho y(\rho))^{-1}
\end{equation}
Note that $h^\prime$ is positive since $|x(\rho)| < \rho$ on the interior of  $A_{m, \beta}$
and $A_{m, \beta}^{(1)}$.

\begin{lem}
There exists a set $A_{m, \beta}$ such that $|x(\rho)| \leq \rho$ if and only if
$\rho \in A_{m, \beta}$. Moreover $|x(\rho)| < \rho$ whenever $\rho \in A_{m, \beta}^o$.
\label{xRhoCharacterization}
\end{lem}
{\bf Proof}
Since the function $x(\rho) + \rho$ maps $(0, \infty)$ continuously and injectively onto 
$(-\infty, \infty)$, define by $r_0$ the positive real number for which $x(r_0) = -r_0$.
Note $x(\rho) < -\rho$ whenever $0 < \rho < r_0$ and $x(\rho) > -\rho$ whenever $ \rho > r_0$.

Setting $\eta(\rho) = \rho - x(\rho)$ we see that $\eta(\rho)$ is convex on $(0, \infty)$ since its second derivative is positive with a unique minimum at $\rho = m$. 

If $\eta(m) > 0$ then $x(\rho) < \rho$ 
on $(0, \infty)$ so that $|x(\rho)| \leq \rho$ on
$[r_0, \infty)$. In this case set $r_1$ and $r_2$ equal to $r_0$.
Otherwise there exist $r_1, r_2$ such that $r_0 \leq r_1 \leq r_2$ and 
$\eta(r_1) = \eta(r_2) = 0$. Moreover $x(\rho) > \rho$ for $\rho \in (r_1, r_2)$
and $|x(\rho)| < \rho$ whenever $\rho \in (r_0, r_1) \bigcup (r_2, \infty)$.

If we define $A_{m, \beta} = [r_0, r_1] \bigcup \, [r_2, \infty)$ we see that 
$|x(\rho)| \leq \rho$ on $A_{m, \beta}$, $|x(\rho)| > \rho$ on 
$(0, \infty) \bigcap A_{m, \beta}^c$, and $|x(\rho)| < \rho$ for $\rho \in A_{m, \beta}^o$

\begin{prop}
\label{hxyDeltaBoundsProp}
Let $x_\delta(r) = m\ln(r) - \beta + \delta$, $y_\delta(r) = (r^2 - x_\delta(r)^2)^{1/2}$,
$v(x,y) = y - m\arctan(x,y)$ and $|\delta| \leq \delta_0$ then uniformly for $\delta \in [0, \delta_0]$
\begin{equation}
\label{estOfHxyOfr}
|v(x_\delta(r), y_\delta(r)) - r + m\pi/2| < x_\delta(r)^2r^{-1}
\end{equation}
for $r$ sufficiently large.
\end{prop}
{\bf Proof}
Simple asymptotic arguments show that there exists an $R$ for which if $r > R$ then
$2m < (\ln(r))^{1/2} < x_\delta(r) < r^{1/4}$, $y_\delta(r) > r - r^{-1/2}$, and thus also $x_{\delta}(r)y_{\delta}(r) > mr$.
If $r > R$ it follows that 
\begin{equation}
\label{arcTanAngleBndWithDelta}
\pi/2 - m^{-1}x_\delta(r)^2r^{-1} < \pi/2 - x_\delta(r)y_\delta(r)^{-1} < \arctan(x_\delta(r), y_\delta(r)) < \pi/2
\end{equation}
and 
\begin{equation}
\label{yrBndWithDelta}
r - x_\delta(r)^2r^{-1} < y_\delta(r) < r
\end{equation}
Inequalities \ref{arcTanAngleBndWithDelta} and \ref{yrBndWithDelta} jointly imply
\begin{equation}
r - x_\delta(r)^2r^{-1} - m\pi/2 < v(x_\delta(r), y_\delta(r)) < r - m\pi/2 + x_\delta(r)^2r^{-1}
\end{equation}
Subtracting $r - m\pi/2$ from the previous inequality yields \ref{estOfHxyOfr}.
\begin{cor}
\label{coshOfrSubnDeltaEst}
Let $\delta_0$ be a fixed positive number and suppose $\{r_n, \epsilon_n, n > 0\}$ where 
$0 < r_n = (2n + 1 + m/2)\pi + \epsilon_n$, and $\epsilon_n = o(1)$. Then uniformly for $\delta \in [0, \delta_0]$
there exists $n_0$ for which 
$\cos(v(x_\delta(r_n), y_\delta(r_n)) < 0$  whenever $n > n_0$.
\end{cor}
{\bf Proof}
By Proposition \ref{hxyDeltaBoundsProp} we may choose $n_1$ sufficiently large so that
$|v(x_\delta(r_n), y_\delta(r_n)) - r_n + m\pi/2| < x_\delta(r_n)^2 r_n^{-1}$ 
whenever $n_1 < n$. Substituting for $r_n$ we find
$|v(x_\delta(r_n), y_\delta(r_n)) -(2n + 1)\pi-\epsilon_n| < x_\delta(r_n)^2 r_n^{-1}$.
Thus
$|v(x_\delta(r_n), y_\delta(r_n)) -(2n + 1)\pi| < |\epsilon_n| + x_\delta(r_n)^2 r_n^{-1}$. We may choose $n_0 \geq n_1$ such that 
$x_\delta(r_n)^2 r_n^{-1} + |\epsilon_n| < \pi/2$. 
Thus if $n > n_0$ we have $\pi /2 < v(x_\delta(r_n), y_\delta(r_n)) -2n\pi < 3\pi/2$.
The last inequality implies $\cos(v(x_\delta(r_n), y_\delta(r_n))) < 0$ whenever $n > n_0$.
\begin{lem}
\label{weakLemRhoSubnAsymp} Let $\rho_n$ be as defined in Lemma \ref{hOfOmegaMonotonicity}.
For any $\eta > 0 $ the sequences $\rho_n$, $x(\rho_n)^2\rho^{-1}$, $x(\rho_n)\rho_n^{-1}$ and $y(\rho_n)$ satisfy the following set of asymptotic relations
\end{lem}
\begin{equation}
\label{x_rho_div_rho_n_asymp}
0 < x(\rho_n)\rho_n^{-1} = o(n^{\eta - 1})
\end{equation}
\begin{equation}
\label{x_rho_sqr_div_rho_n_asymp}
x(\rho_n)^2\rho_n^{-1} =  o(n^{\eta - 1})
\end{equation}
\begin{equation}
\label{weak_rho_n_estimate}
\rho_n \geq  h(\rho_n) = 2 n\pi
\end{equation}
\begin{equation}
\label{weak_y_rho_n_estimate}
0 < \rho_n - y(\rho_n) = o(n^{\eta - 1})
\end{equation}
\begin{equation}
\label{weak_omega_rho_n_estimate}
\pi/2 - m^{-1}x(\rho_n)^2\rho_n^{-1} < \omega(\rho_n)  < \pi /2 
\end{equation}

{\bf Proof}
Basic asymptotics show $0 < x(\rho) = o(\rho^{\eta/2})$ as $\rho \rightarrow \infty$ for any $\eta > 0 $. This in turn implies that $x(\rho)^2\rho^{-1}$ and
$x(\rho)\rho^{-1}$ are $ o(\rho^{\eta - 1})$. To prove \ref{x_rho_div_rho_n_asymp} and \ref{x_rho_sqr_div_rho_n_asymp} we need only establish \ref{weak_rho_n_estimate}.
Since $y(\rho) \leq \rho$ and $\omega(\rho) \geq  0$ on $A_{m,\beta}$ it follows definitionally that $\rho_n \geq h(\rho_n) = 2n\pi$. 
Thus $x(\rho_n)^2\rho_n^{-1}$ and $x(\rho_n)\rho_n^{-1}$ are $ o(n^{\eta - 1})$ as claimed.

To establish \ref{weak_y_rho_n_estimate} note whenever $|b| < a$,
 $a - b^2a^{-1} <(a^2- b^2)^{1/2} < a$. Applying this inequality to the definition of $y(\rho)$ we find
\begin{equation}
\rho_n -x(\rho_n)^2\rho_n^{-1} < y(\rho_n) < \rho_n 
\end{equation}
To complete the proof of the lemma we now justify \ref{weak_omega_rho_n_estimate}.
Noting that if we set $\delta$ to zero in  inequality \ref{arcTanAngleBndWithDelta} we get
\begin{equation}
\pi/2 - m^{-1}x(r)^2r^{-1}  < \arctan(x(r), y(r)) < \pi/2
\end{equation}

\begin{lem}
\label{strongLemRhoSubnAsymp}
The sequences $\rho_n$, $x(\rho_n)$ and $y(\rho_n)$ satisfy the follow set of asymptotic equations as 
$n \rightarrow \infty$.
\begin{equation}
\label{strong_rho_n_asym}
\rho_n = (2n + m/2)\pi + o(n^{\eta - 1})
\end{equation}
\begin{equation}
\label{strong_y_rho_n_asymp}
y(\rho_n) = (2n + m/2)\pi + o(n^{\eta - 1}) 
\end{equation}
\begin{equation}
\label{strong_x_rho_n_asymp}
x(\rho_n) = m \ln(2n\pi) + 4^{-1}m^2n^{-1} -\beta + o(n^{\eta - 2})
\end{equation}

\end{lem}

{\bf Proof}
Note \ref{strong_rho_n_asym} and Lemma \ref{weakLemRhoSubnAsymp} in combination show \ref{strong_y_rho_n_asymp}.
Lemma \ref{weakLemRhoSubnAsymp} and $h(\rho)$'s definition jointly imply
\begin{equation}
\label{hRhoInequality}
\rho_n - x(\rho_n)^2\rho_n^{-1} - m\pi/2  < h(\rho_n) < 
\rho_n - m\pi/2 + x(\rho_n)y(\rho_n)^{-1}
\end{equation}
Using the identity $h(\rho_n) = 2n\pi$ in \ref{hRhoInequality}
yields 
\begin{equation}- x(\rho_n)y(\rho_n)^{-1} < \rho_n -(2n + m/2)\pi < 
x(\rho_n)^2\rho_n^{-1} 
\end{equation}
The last inequality and Lemma \ref{weakLemRhoSubnAsymp} jointly yield \ref{strong_rho_n_asym}.

Observe 
\begin{equation}
\ln(\rho_n) = \ln((2n + m/2)\pi + o(n^{\eta - 1})) = 
\ln(2n\pi) + m4^{-1}n^{-1} + o(n^{\eta - 2})
\end{equation}
The last equality and the definition of $x(\rho)$ together imply \ref{strong_x_rho_n_asymp}.
\begin{prop}
\label{estOfExpAroundZeroProp}
Set $\sigma_0(w) = \exp(w) - 1$ and $w = x + iy$, then $\sigma_0$  satisfies the following inequalities

\begin{equation}
|\sigma_0(w)| \geq \min(1/2, |x|/2)
\label{sigma0LowerBoundEst}
\end{equation}
\begin{equation}
0.5|w| \leq |\sigma_0(w)| \leq 1.5|w|
\end{equation}
for $|w| < 1/2$,
\begin{equation}
|\sigma_0(w)| > 1 
\end{equation}
for $\cos(y) < 0$,
\begin{equation}
\label{shiftedHyperbolicSineIneq}
|\sigma_0(w)|^2 \geq 1 - \cos(y)^2 
\end{equation}
\end{prop}
{\bf Proof}
To verify the first inequality  note that for $u \geq 0$, $|\exp(-u) - 1| \geq \min(1/2, |u|/2)$ which may be proved by a Taylor series expansion with remainder argument. Substituting this inequality into
$|\sigma_0(w)| \geq ||\exp(w)| - 1| \geq |\exp(-|x|) - 1| \geq  \min(1/2, |x|/2)$ yields the first inequality. 

The second inequality follows from 
\begin{equation}
0.5|w| \leq |w| - |w|^2(1- |w|)^{-1} \leq |\sigma_0(w)| \leq |w|+ |w|^2(1- |w|)^{-1} < 1.5|w|
\end{equation}
for $|w| < 1/2$. All of which may be justified by substituting each term with its Taylor series expansion.

The third inequality is proved by observing
\\
 $|\sigma_0(w)| \geq |\Re{w}| = 1 + |\cos(y)|\exp(x) > 1$.

To prove the last inequality note that $\exp(2x) + 1 -2\cos(y) \exp(x) \geq \exp(2x) + 1 -2|\cos(y)|\exp(x) $ which reaches an absolute minimum of $1 - \cos(y)^2 $ at $x = -\ln|\cos(y)|$.
\begin{lem}
Let $h(\omega, \eta, \xi, m) = \exp(\omega\eta) (1 + \eta/\xi)^{-m} -1$ for $\omega \neq 0$ and $m$ a positive integer then
\label{tQuasiLinearityLemma}
\begin{equation}
\label{linearBoundsOnhEtaXi}
\mu(\omega)|\eta| \leq |h(\omega, \eta, \xi,m)| \leq \nu(\omega)|\eta |
\end{equation}
whenever  $|\eta| < \delta(\omega)$
and $|\xi| > M_0(m, \omega)$, where $\delta(\omega) = \min(0.5, |1.5\omega|^{-1})$,
\linebreak 
$M_0(m,\omega) = (8|\omega|^{-1} + 1)m$,  $\mu(\omega) = 4^{-1}|\omega|$ and $\nu(\omega) = 7|\omega|/4$.
\end{lem}
{\bf Proof}
\begin{equation}
h(\omega, \eta, \xi, m) = \exp(\omega \eta)(1 + \eta/\xi)^{-m} -1 = g(\eta,\xi)+  \sigma_0(\omega\eta) 
\end{equation} 
where $g(\eta, \xi) = \exp(\omega\eta)(1 + \eta/\xi)^{-m} - \exp(\omega\eta)$ and where $\sigma_0$ as defined in Proposition \ref{estOfExpAroundZeroProp}. 

Proposition \ref{estOfExpAroundZeroProp} implies 
\begin{equation}
\label{tauLinearBounds}
0.5|\omega\eta| \leq |\sigma_0(\omega\eta)| \leq 1.5  |\omega\eta|
\end{equation}
It follows from inequality \ref{tauLinearBounds} that $|\exp(\omega\eta)| < 2$ 
whenever $|\eta| < \delta(\omega)$.
By setting $k = 8 |\omega^{-1}|$ and $ v = \eta$ in Proposition \ref{negBinIneqProp} we may deduce that for $|\xi| > M_0(m, \omega)$
\begin{equation}
\label{hEtaUnBounds}
|g(\eta, \xi)| \leq 2 |\eta|(8 |\omega^{-1}|)^{-1} = 4^{-1}|\omega\eta|
\end{equation}
Inequality \ref{linearBoundsOnhEtaXi} follows jointly from
\begin{equation}
|\sigma_0(\omega\eta)| - |g(\eta, \xi)| \leq |h(\omega, \eta, \xi, m)| 
\leq |\sigma_0(\omega\eta)| + |g(\eta, \xi)|
\end{equation}
and inequalities \ref{tauLinearBounds} and \ref{hEtaUnBounds} 
\begin{cor}
\label{semiLinearityOfT}
Let $T(z) = c\exp(\alpha z)z^{-m} -1$ where $m$ is a positive integer, $\alpha < 0$ and $c \alpha^m > 0$. Let $\{u_n, n > 0\}$ be the zeroes of $T(z)z^m$ as per Corollary 
\ref{T_RepresCorollary}. Then there exist $n_0, \delta(\alpha) > 0, \mu(\alpha) > 0, \nu(\alpha) > 0$ such that
if $n > n_0$ and $|\eta| < \delta(\alpha)$ then 
\begin{equation}
\label{linearBoundsOnTandU_n}
\mu(\alpha) |\eta| \leq |T(u_n + \eta)| \leq \nu(\alpha) |\eta|
\end{equation}
\end{cor}
{\bf Proof}
By simple algebra we see that
\begin{equation}
T(u_n + \eta) = \exp(\alpha \eta)(1 + \eta/u_n)^{-m} -1
\end{equation} 
Choose  $n_0$ so large that $n > n_0$ we have $|u_n| > (8|\alpha|^{-1} + 1)m$, and apply Lemma \ref{tQuasiLinearityLemma} where $\xi = u_n$ and $\omega = \alpha$ to complete the proof.

Our next result shows that the linear growth bounds on $T$ at a zero holds in neighborhood of the zero.
\begin{cor}
\label{uniformityOfSemiLinearityOfT}
Let $\nu(\alpha)$, $\delta(\alpha)$, $n_0$ and $u_n$ be those of Lemma \ref{semiLinearityOfT}. Set
$\delta_1(\alpha) =\min(\delta(\alpha), 1/4\nu(\alpha))$ and suppose $n > n_0$,
$|w - u_n| < \delta_1(\alpha)$ and
$|\eta| <\delta(\alpha)$ then
\begin{equation}
\label{extendedLinearInequalityForT}
0.75 \mu(\alpha) |\eta| \le |T(w + \eta) - T(w)| \leq 1.25 
\nu(\alpha)|\eta|
\end{equation}
\end{cor}
{\bf Proof} We have from Corollary \ref{semiLinearityOfT} that
\begin{equation}
|T(w)-T(u_n)| 
\le \nu(\alpha) |w-u_n| \le 0.25
\end{equation}
By the triangle inequality we deduce that
\begin{equation}
\label{boundOnTPlus1NearAZero}
0.75 < |c\exp(\alpha w)w^{-m}| < 1.25
\end{equation}
Observe that 
$|T(w + \eta) - T(w)| = |c\exp(\alpha w)w^{-m}
h(\alpha, \eta, w, m)| $
which by Lemma \ref{tQuasiLinearityLemma} and inequality \ref{boundOnTPlus1NearAZero} implies  \ref{extendedLinearInequalityForT}.

\begin{cor}
\label{BoundingTsLinearlyForGrowingZ}
Let $\nu(\omega)$, $\delta(\omega)$ and $M_0(m, \omega)$ be those of Lemma \ref{tQuasiLinearityLemma}, and let $\beta_1$ be a real number and let $H(z, \omega,m,j) = c\exp(\omega z)z^{-m-j}$ where $\omega < 0$, $j$ is a non-negative integer, and $\alpha \le \omega$. If 
$ \omega -\alpha + j > 0$ then there exists $d > 0$ for which $|z| > M_0(m +j, \omega)$, $\Re{z} > m\log(|z|)/\alpha + \beta_1$ and $|\eta| < \delta(\omega)$ jointly imply $|H(z+\eta, \omega,m,j)-H(z, \omega,m,j)| < \nu(\omega)|c\eta z^{-d}|\exp(\omega\beta_1)$.
\end{cor}
{\bf Proof}
We may rewrite $|H(z + \eta, \omega,m,j) - H(z, \omega,m,j)|$ as
$|c\exp(\omega z)z^{-m-j}||h(\eta, \omega, m+j)|$ where $h(\eta, \omega, m + j)$ is that defined in Lemma \ref{tQuasiLinearityLemma}.
Note that $|\exp(\omega z)z^{-m -j}| = \exp(\omega \Re{z}) |z|^{-m-j}| < \exp(\omega\beta_1 + \omega m \log(|z|)/\alpha)|z|^{-m-j} = \exp(\omega\beta_1)|z|^{-m -j +m \omega/\alpha}$. If we set $d =  m(1 - \omega/\alpha) +j$ the corollary now follows from Lemma \ref{tQuasiLinearityLemma}.

\begin{cor}
\label{BoundingDeltaSumsOfTsNonDominantTerms}
Let $\beta_1$ be a fixed, real number. Set ${\mathcal H}(z) = \sum_{1\leq k \le k_0}c_k\exp(\omega_k z)z^{-m-j_k}$ where $\omega_k < 0$, $j_k$ is a non-negative integer, and $\omega_0 \le \omega_k$. If 
$ \min_{1 \le k \le k_0} (\omega_k -\omega_0 + j_k) > 0$ 
then there exist $d_{\mathcal H} > 0, \nu_{\mathcal H} > 0, \delta_{\mathcal H} > 0$ and $M_{\mathcal H} > 0$  for which $|z| > M_{\mathcal H}$, $\Re{z} > \beta_1 + m\log(|z|)/\omega_0$ and $|\eta| < \delta_{\mathcal H}$ jointly imply $|{\mathcal H}(z+\eta)-{\mathcal H}(z)| < \nu_{\mathcal H}|\eta z^{-d_{\mathcal H}}|$.
\end{cor}
{\bf Proof}
Define $\nu_{\mathcal H} = \sum_{1 \leq k \le k_0}\nu(\omega_k)|c_k|\exp(\omega_k\beta_1)$,
$\delta_{\mathcal H} = \min_{1\leq k \le k_0}\delta(\omega_k)$,

$d_{\mathcal H} = \min_{1\leq k \le k_0}d(m,j_k,\alpha, \omega_k )$ and $M_{\mathcal H} = \max_{1\le k \le k_0}M_0(m+j_k, \omega_k)$, 

where $\delta(), \nu(), M_0(),d()$ are as defined in Corollary \ref{BoundingTsLinearlyForGrowingZ}.
\begin{equation}
|{\mathcal H}(z+\eta)-{\mathcal H}(z)| \leq \sum_{1\le k \le k_0}
|H(z + \eta, \omega_k,m,j_k) -H(z, \omega_k,m,j_k)|
\end{equation}
applying Corollary \ref{BoundingTsLinearlyForGrowingZ} to each individual term on the right hand side we get
for $|z| > M_0$ and $|\eta| <\delta_{\mathcal H}$
\begin{equation}
|{\mathcal H}(z+\eta)-{\mathcal H}(z)| < \sum_{1\le k \le k_0}
\nu(\omega_k)\exp(\omega_k\beta_1)|c_k\eta z^{-d(m,j_k,\alpha, \omega_k)}|< \nu_{\mathcal H}|\eta z^{-d_{\mathcal H}}|
\end{equation}
\begin{prop}
\label{negBinIneqProp}
If $0 \leq |v| \leq 1 $ and $ (k + 1)m \leq |u|$ where $m \geq 1$ and $k > 0$, we have the following inequality
\begin{equation}
|(1 + vu^{-1})^{-m} - 1| \leq |v|/k
\end{equation}
\end{prop}
{\bf Proof}
Applying the binomial theorem in conjunction with triangle inequality yields
\begin{equation}
\label{negBinThmSeries}
|(1 + vu^{-1})^{-m} - 1| \leq \sum_{j = 1}^\infty |v|^j |u|^{-j}\binom {m + j - 1} {j}
\end{equation}
Using the inequality $m + n - 1 \leq mn $ for $n \geq 1$ one may prove
\begin{equation}
\label{binCoeffEst}
\binom {m + j - 1} {j} m^{-j} \leq 1
\end{equation}
Substituting \ref{binCoeffEst} into \ref{negBinThmSeries} we have
\begin{equation}
|(1 + vu^{-1})^{-m} - 1| \leq  \sum_{j = 1}^\infty  |v|^j(k + 1)^{-j} = |v|(k + 1 - |v|)^{-1} \leq |v| /k
\end{equation}
\begin{lem}
\label{horseshoeRegionForRouche}
Let $S(z) = \exp(z + \beta)z^{-m} - 1$ and $a$ be a real fixed number then there exist an integer $n_0(a)$ and horseshoe shaped curves  ${\mathcal D}_n(a)$ for $n \geq n_0(a)$ such that if $n\geq n_0(a)$ and $z \in {\mathcal D}_n(a)$ then $|S(z)| \geq 1/2$ where
${\mathcal D}_n(a)$ is defined by ${\mathcal C}_n(a) \bigcup {\mathcal A}_n(a) \bigcup {\mathcal C}_{n +1}(a)$
and where
${\mathcal C}_k(a) = \{z: \Re{z} \geq a $  and  $|z| = r_k\}$ and
${\mathcal A}_k(a)$ is  $\{z: \Re{z} = a $  and  $|z| \in [r_k, r_{k + 1}]\}$ where 
$r_k = (|z_k| + |z_{k + 1}|)/2$.
Moreover $S(z)$'s roots inside of ${\mathcal D}_n$
are exactly $z_{n+1}$ and $\overline{z}_{n+1}$ with
$\Im{z_n} = (2n + m/2)\pi + o(n^{\eta -1})$ and
$\Re{z_n} = m\ln(2n\pi) + 4^{-1}m^2n^{-1} -\beta + o(n^{\eta -2})$
\end{lem}
{\bf Proof}
We may assume w.l.o.g. that $\Im{z} \geq 0$ 
which follows from $S(\bar{z}) =\overline{S(z)}$, next we may rewrite $S$ as $|S(z)| = |\exp(u(x,y) + iv(x,y)) - 1|$ where $z = x + iy$,
$u(x,y) = x - m\ln|z| -\beta$ and $v(x,y) = y - m\arctan(x, y)$.

Since by Lemma \ref{strongLemRhoSubnAsymp} we have $r_j = (2j + 1 + m/2)\pi + o(1)$, we 
can choose an $n_1$ such that
$a < m\ln(r_j) + \beta - 1$ whenever $j > n_1$, hence by definition
$|u(x,y)| \geq 1$ when $z \in {\mathcal A}_n(a)$ and $n > n_1$. We may now apply inequality \ref{sigma0LowerBoundEst} of Lemma \ref{estOfExpAroundZeroProp} to complete the proof for ${\mathcal A}_n(a)$.

For $z \in {\mathcal C}_n(a) \bigcup {\mathcal C}_{n +1}(a)$ it suffices by Proposition \ref{estOfExpAroundZeroProp} to prove that if $|u(x,y)| < 1$ then
$\cos(v(x,y)) < 0)$. We may assume w.l.o.g. that $z \in {\mathcal C}_n(a)$. By Corollary \ref{coshOfrSubnDeltaEst} if we set $\delta_0 = 1$ we may select $n_0(a) \geq n_1$ for which 
$\cos(v(x_\delta(r_n), y_\delta(r_n)) < 0$  
whenever $n > n_0(a)$ and $\delta < 1$. If we set $\delta = |u(x,y)|$ we find by definition that 
$x_\delta(r_n) = m \ln(r_n) -\beta + u(x, y)  =  x$ which in turn implies $y_\delta(r_n) = |y| = y$ and that $\cos(v(x,y)) < 0 $

Lastly the assertion for the roots $z_n$ follows from Theorem \ref{CompareFuncExp} and the fact that for $n_0(a)$ sufficiently large guarantees that 
$r_n < |z_{n + 1}| < r_{n+1}$.
\begin{cor}
\label{horseshoeRegionForT}
Let $T(\theta) = c\exp(\alpha\theta)\theta^{-m} - 1$
where $c \alpha^m > 0$. For every fixed real number 
$a$ there exist an integer $n^\prime_0(a)$ and horseshoe shaped curves ${\mathcal D}_n^\prime(a)$ for $n \geq n^\prime_0(a)$ such that if $n\geq n^\prime_0(a)$ and $\theta \in {\mathcal D^\prime}_n(a)$ then $|T(\theta)| \geq 1/2$ where
${\mathcal D}_n^\prime(a)$ is defined by ${\mathcal C}_n^\prime(a) \bigcup {\mathcal A}_n^\prime(a) \bigcup 
{\mathcal C}_{n +1}^\prime(a)$
and where
${\mathcal C}_k^\prime(a) = \{z: \Re{z} < a $  and  $|z| = r_k^\prime\}$ and
${\mathcal A}_k^\prime(a)$ is  $\{z: \Re{z} = a $  and  $|z| \in [r_k^\prime, r_{k + 1}^\prime]$ for 
 $\{r_n^\prime \}$ a sequence that strictly increases to infinity.
 
 Moreover $T(z)$'s roots inside of ${\mathcal D}_n^\prime$
 are exactly $u_{n+1}$ and $\overline{u}_{n+1}$ with
 $\Im{u_n} = \alpha^{-1}(2n + m/2)\pi + o(n^{\eta -1})$ and
 $\Re{u_n} = \alpha^{-1}m\ln(2n\pi) + 4^{-1}m^2n^{-1} -\beta/\alpha + o(n^{\eta -2})$.
\end{cor}
{\bf Proof}
Observe that $T(z) = S(\alpha z)$ for $S$ defined by 
$S(z) = \exp(z + \beta)z^{-m} - 1$ where 
$\beta = m\log(|\alpha|) + \log(|c|)$. Let 
${\mathcal C}_n(\alpha a)$, ${\mathcal A}_n(\alpha a)$,
$ r_n $ and $n_0(a)$ be those defined in Lemma \ref{horseshoeRegionForRouche}. If we set
$r_n^\prime =  r_n / |\alpha|$ and $n^\prime_0(a) = n_0(\alpha a)$
it is easy to show that if
$ z \in {\mathcal C}_n^\prime(a)$ and $n \ge n^\prime_0(a)$ then $\alpha z \in {\mathcal C}_n(a\alpha)$, and if 
$z \in {\mathcal A}_n^\prime(a)$ then $\alpha z \in {\mathcal A}_n(a\alpha)$.
The proof now follows from Lemma \ref{horseshoeRegionForRouche} by setting $u_k = z_k/\alpha$ and applying the identity $T(z) = S(\alpha z)$.

\subsection {Possible Generalizations}
As we will see in the next section there are G/G/1 systems for which $\sigma$ as defined by \ref{tDefEqn} doesn't dominate $F(\theta)$ as $\Re{\theta} \rightarrow -\infty $. 

The function $\sigma$ may be varied by adding $\theta^n$ instead of subtracting it:
$\hat{\sigma}(\theta) = \exp(\theta + \beta) + \theta^n$.
One may use methods detailed above to prove   $\hat{\sigma}$ has canonical form very similar to that of $\sigma$ given by Theorem \ref{CompareFuncExp}.

In some instances the dominant term of $F$ will have multiple exponents such as $T(\theta)=\beta_0 \beta_1\exp((\alpha_0 +\alpha_1)\theta)\theta^{-m_0 -m_1} + \beta_1\exp( \alpha_1\theta)\theta^{-m_1} -1$ where $\alpha_i < 0$ and the $m_i$ are positive integers satisfying $(m_0 +m_1)/(\alpha_0 +\alpha_1) <  m_1/\alpha_1 $. 

One may express $T(\theta)$ as $T_0(\theta)T_1(\theta) + \beta_0 \exp(\alpha_0\theta)$ where $T_0(\theta) = \beta_0 \exp(\alpha_0\theta)\theta^{-m_0} +1$ and $T_1(\theta) = \beta_1 \exp(\alpha_1\theta)\theta^{-m_1} -1$ The inequality on $\alpha_i$ and $m_i$ simplifies to $m_0/\alpha_0 < m_1/\alpha_1$. This last inequality guarantees that when $T_0(\theta)$ is close to $0$ for large $\theta$ that $T_1(\theta)$ grows large with $\theta$. Also for $T_1(\theta)$ close to zero for $\theta$ large implies that $\beta_0\exp(\alpha_0\theta)\theta^{-m_0}$ will be decreasingly small compared to $\beta_1(\exp(\alpha_1\theta))\theta^{-m_1}$ allowing one to use Lemma \ref{rouchForTFunctions}

We will now illustrate these generalizations with the following example.

\subsubsection{A queue whose dominant term has multiple exponentials}
\label{ExploreOfMultiplyDominatedQueue}
Let the interarrival time be the constant $1/2$ and service time distribution defined on $[0,1]$ as the fifty-fifty mixture of a Unif[0, 7/8] and ${\mathcal B}(1,2)$ where ${\mathcal B}(1,2)$ is the beta distribution defined on $[0,1]$ by the density $2(1-x)$. 
The service time transform is given by 

$B(\theta) = \exp(-\theta)\theta ^{-2} - \theta^{-2} + \theta^{-1} + 4\left (1 -\exp(-7\theta/8)\right)
(7\theta)^{-1}$

Setting  $F(\theta) = H(\theta) + G(\theta)$

where

$H(\theta) =  \exp(-\theta/2)\theta ^{-2} -4\exp(-3\theta/8) (7\theta)^{-1} -1 $

and

$G(\theta) = \exp(\theta/2)(4(7\theta)^{-1} -\theta^{-2})$

$T_0(\theta) = 7\exp(-\theta/8)(4\theta)^{-1} - 1$

$T_1(\theta) = 4\exp(-3\theta/8)(7\theta)^{-1} + 1$

Computational evidence and the heuristics given above suggest that $H(\theta)= q(\theta)\hat{T}_1(\theta)\hat{T}_0(\theta)$
where $q(\theta)$ is the ratio of low degree polynomials, and the $\hat{T}_i$ have the same form as $T_i$ with asymptotically convergent zeroes.

For example solving the simultaneous equations \ref{genericEquationForZxForTeeNTau} and \ref{genericEquationForZyForTeeNTau}  for $k = 10, 100$ when the target function is $T_0$ and $T_1$ respectively and then solving for $H$'s and $F$'s associated zeroes we find that

\[
\begin{tabular}{|c|c|r|r|}\hline
\multicolumn{4}{|c|}{Sample Root Comparison for $T_0$, $T_1$ with $H$}\\ \hline
\multicolumn{1}{|c|}{k}&
\multicolumn{1}{|c|}{Function}& 
\multicolumn{1}{c|}{$Z_k$} &
\multicolumn{1}{l|}{$ \Im{Z_k}/\alpha\pi$} \\ \hline
10&$T_0$&-45.879622 + 539.675461j& 21.473004 \\ \hline
10&$ H$&-45.879369 + 539.675421j& 21.473003 \\ \hline
10&$F$&-45.879369 + 539.675421j& 21.473003 \\ \hline
10&$T_1$&-15.221727 + 171.504339j& 20.471823 \\ \hline
10&$ H$&-15.132358 + 171.351343j& 20.453560 \\ \hline
10&$F$&-15.132361 + 171.351346j& 20.453560 \\ \hline
100&$T_0$&-63.763235 + 5064.146634j&201.495992 \\ \hline
100&$ H$&-63.763232 + 5064.146634j&201.495992 \\ \hline
100&$F$&-63.763232 + 5064.146634j&201.495992 \\ \hline
100&$T_1$&-21.296132 + 1679.671064j&200.495964 \\ \hline
100&$ H$&-21.276224 + 1679.636887j&200.491885 \\ \hline
100&$F$&-21.276224 + 1679.636887j&200.491885 \\ \hline
\end{tabular}
\]

\mysection {Exponential Asymptotics}
The notion of exponential order introduced in this section 
is the key to obtaining an exact mathematical
formula for $\psi$, the waiting time Laplace transform. 
An exponential order computation plays a key role in pinning 
down the exponential term in $\psi$'s Hadamard product representation.

\begin{defin}
\label{expOrdDef}
	A real function $f$ is said to be of exponential
	order $a$ at $b$ if for $\forall \epsilon >0$
		\begin{equation}
			\lim_{x \rightarrow b} f(x)\exp(-(a+\epsilon)\mid x\mid)
			= \lim_{x \rightarrow b} \exp((a-\epsilon)\mid x\mid)/f(x) =0
		\end{equation}
\end{defin}
Note that we omit $b$ if the exponential order is the same for $\pm\infty$.
It is easy to prove that all non-null polynomials are of order 0 for $b$ = $\pm\infty$.

\begin{lem}
	Suppose that $\{a_{n},\ n=1,2,...\}$ 
	is sequence for which
	\begin{equation}
		a_{n} = O(n^{-1-\eta}) \ where \ 0 < \eta < 1
                \label{aSubNEquation}
	\end{equation}

	and $\{b_{n},\ n=1,2,..\}$ is a sequence of positive numbers for which
	\begin{equation}
		0 < b_{n}^{-1} = n^{2} \ + \ O(n)
		\label{bSubNEquation}
	\end{equation}
	Then $\gamma(x)$ defined by
	\begin{equation}
		\gamma(x)\equiv \prod_{n=1}^{\infty}\ (1 +a_{n}x +b_{n}x^{2})
	\end{equation}
	is of exponential order $\pi$.
	\label{prod_order}
\end{lem}
{\bf Proof}
	Define
	$\zeta(x) = \prod_{n=1}^{\xi(x)} (1 +a_{n}x +b_{n}x^{2})$,
	where
	$\xi(x) =\lceil |x|^{\alpha}\rceil$ with
	$\alpha = 1+\eta$, and $\Gamma(x) = \gamma(x)/\zeta(x)$.
	We will now show that $\zeta$ is of exponential order $\pi$.
	Let $0<\epsilon <1$ be given.
	If $1\leq n \leq \xi(x)$ then
	$|a_{n}/(b_{n}x)| = O(\xi(x)^{1-\eta})/|x|=O(|x|^{-\eta^2})$.
	It follows that by making $|x|$ sufficiently
	large we can ensure that $|a_n / (b_n x)| < \epsilon$ for $1\leq n \leq \xi(x)$, and thus

		\begin{equation}
			\prod_{n=1}^{\xi(x)} (1 +b_{n}(1-\epsilon)x^{2})
			\leq
			\prod_{n=1}^{\xi(x)} (1 +b_{n}x^{2}(1+a_{n}/(b_{n}x)))
			\leq
			\prod_{n=1}^{\xi(x)} (1 +b_{n}(1+\epsilon)x^{2})
		\end{equation}
	The last equation and Corollary \ref{funct_asym} jointly imply that
	$\zeta$ is of exponential order $\pi$.
	Proposition \ref{xi-tail-is-zero} shows that
	$\Gamma$ is of exponential order 0.
	Applying the multiplicative law (cf. Proposition \ref{prod_exp_order}) for exponential orders to the product
	of $\zeta(x)\Gamma(x)$ completes the proof.

\begin{cor}
\label{expBehaviorEstCor}
If the $b_n$ in Lemma \ref{prod_order} satisfy $0\leq b_{n}^{-1} = 4\pi^2 \alpha^{-2}n^{2} \ + \ O(n)$
where $\alpha > 0$
then the exponential order of $\gamma(x)$ is $\alpha/2 $.
\end{cor}
\begin{cor}
Let $a_n$, $b_n$ and $\gamma(x)$ be as in Corollary \ref{expBehaviorEstCor} then for every $\epsilon > 0$ a constant
$C(\epsilon)$ such that
     \begin{equation}
     |\gamma(\theta)| < C(\epsilon)\exp(|\theta|(\alpha/2 + \epsilon))
     \end{equation}
\label{prodExpoBound}
\end{cor}
{\bf Proof}
Without loss of generality we may assume that $a_n > 0$
     \begin{equation}
     |1 + a_n \theta + b_n \theta^2| \leq 1 + a_n |\theta| + b_n |\theta|^ 2
     \label{gammaTermIneq}
     \end{equation}
Forming the product of \ref{gammaTermIneq} over $n$ yields
$|\gamma(\theta)| \leq \gamma(|\theta|)$. By Corollary
\ref{expBehaviorEstCor} there exists $M > 0$ such that if $x > M$ then
$|\gamma(x)| < \exp(x(\alpha/2 + \epsilon))$.
Choosing $C(\epsilon)$ to be the maximum of $|\gamma(\theta)|$ on $|\theta| \leq M$
completes the proof.
\begin{prop}
	Suppose that $\{b_{n},\ n=1,2,..\}$.
	is a sequence for which
		\begin{equation}
			0\leq b_{n}^{-1} =  n^{2}(1 \ + \ O(n^{-1}))
			\label {b-asym}
		\end{equation}
	Then $\gamma(x)$ defined by
		\begin{equation}
			\gamma(x)\equiv \prod_{n=1}^{\infty}\ (1 + b_{n}x^{2})
		\end{equation}
	is of exponential order $\pi$.
	\label{sinhLikeBehavior}
\end{prop}
{\bf Proof}
	Note that $\sinh(\pi x) = (\exp(\pi x) - \exp(-\pi x)) / 2$ is of exponential order $\pi$ and
	has product form
		\begin{equation}
			\sinh(\pi x)=\pi x\prod_{n=1}^{\infty}\ (1\ + \ x^{2} n^{-2}).
		\end{equation}
	Observe by virtue of \ref{b-asym}
	there exists an $n_{0}$ for which $n\geq n_{0}$
		\begin{equation}
			\label{bnInequality}
			(n-n_0)^{2}\leq b_{n}^{-1} \leq (n+ n_0)^{2}
		\end{equation}
        Thus for $n > n_0$
		\begin{equation}
  			 1\ + x^{2} (n + n_{0})^{-2}\leq  1 + b_{n}x^{2}
			 \leq 1\ + \ x^{2} (n - n_{0})^{-2}
                 \label{bnPolyInequality}
		\end{equation}
        Since the leftmost term in \ref{bnPolyInequality} is positive we may form the
        product over $n_0 < n < \infty$ to get
  		\begin {equation}
			\sinh(\pi x)/p_{2n_{0}}(x) \leq \gamma(x)/\gamma_{n_{0}}(x)
			\leq \sinh(\pi x)  / (\pi x)
			\label{sine-sand}
		\end{equation}
        Where $\gamma_{n_{0}}(x) = \prod_{n=1}^{n_{0}}\ (1 +b_n x^{2})$
	and $p_{n}(x)= \pi x \prod_{k=1}^{n}\ (1 +x^{2}k^{-2})$.
        Since $\sinh(\pi x)/ (\pi x)$ and $\sinh(\pi x)/p_{n}(x)$ have order $\pi$ it
        follows from \ref{sine-sand} and Proposition \ref{prod_exp_order}
        that $\gamma(x)/\gamma_{n_{0}}(x)$ has
        exponential order $\pi$. Because $\gamma_{n_{0}}(x)$ is of order 0, $\gamma(x)$
        must be of exponential order $\pi$.

\begin{prop}
	Suppose $\{a_{n}\}$ and $\{b_n\}$ satisfy \ref{aSubNEquation} and \ref{bSubNEquation} respectively
	and set $Q(x)= \lceil |x|^{1 + \eta}\rceil + 1$, then
	$\Gamma(x)$ defined by
		\begin{equation}
			\Gamma(x)\equiv \prod_{n=Q(x)}^{\infty} (1 + a_{n}x + b_{n}x^{2})
		\end{equation}
	is of exponential order 0.
	\label{xi-tail-is-zero}
\end{prop}
{\bf Proof}
Set $A_n = \sum_{k=n}^{\infty} |a_k|$, $B_n = \sum_{k=n}^{\infty} |b_k|$ and
$h(x) = A_{Q(x)}|x| + B_{Q(x)}x^2$.
From 4.2.37 of \cite{abram} we have $\exp(-2|c|) \leq 1 + c \leq \exp(|c|)$ when $|c| \leq 1/2$ .
Take $|x|$ sufficiently large so that if $n \geq Q(x)$ implies $|a_{n}x| + |b_{n}|x^{2} < 1/2$ and thus

	\begin{equation}
         \exp(-2(|a_{n}x| + |b_{n}|x^{2})) \leq
	(1 + a_{n}x + b_{n}x^{2}) \leq \exp(|a_{n}x| + |b_{n}|x^{2})
	\label{expLinBounds}
	\end{equation}
Applying \ref{expLinBounds} to the factors of $\Gamma(x)$ yields
		\begin{equation}
			\exp(-2h(x)) \leq
			\Gamma(x)
			\leq  \exp(h(x))
    \end{equation}
	The estimates $A_n = O(n^{-\eta})$ and $B_n = O(n^{-1})$ jointly imply that
	$h(x) = O(|x|^{1-\eta})$. Consequently $\exp(-2h(x))$ and $\exp(h(x))$ are both of exponential
	order 0. The result now follows from Proposition \ref{sandwichProp}.
\begin{prop}
        If $f$ and $g$ are of exponential order $a$ at $b$ and $|f| \leq |h| \leq |g|$ for some neighborhood of $b$,
        then $h$ is also of exponential order $a$ at $b$.
        \label{sandwichProp}
\end{prop}

\begin{prop}
	If $f$ and $g$ are of exponential order $c$, $d$ respectively then $fg$ is
	of exponential order $c+d$
	\label{prod_exp_order}
\end{prop}
	Propositions \ref{sinhLikeBehavior}, \ref{xi-tail-is-zero} and \ref{prod_exp_order}
	jointly imply the next corollary.

\begin{cor}
	Let $\{b_{n}\}$ satisfy \ref{bSubNEquation}, and let $\xi(x)= \lceil |x|^{1 + \eta}\rceil$
	where $0< \eta < 1$, then
	$\zeta(x)$ defined by
		\begin{equation}
			\zeta(x)\equiv \prod_{n=1}^{\xi(x)} (1 +b_{n}x^{2})
		\end{equation}
	is of exponential order $\pi$.
	\label{funct_asym}
\end{cor}
\mysection {Auxiliary Results from Complex Analysis}
The following theorem is a simplification of the more general form given Ahlfors \cite{ahlfors}, Chapter 5, Section 1.5, equation 19.
\begin{thm}
\label{HadamardProdTheorem}
Hadamard Factorization Theorem. Let $H$ be analytic function satisfying
	\begin{equation}
	|H(\theta)| \leq C\exp(M\mid\theta\mid)
	\label{hExpoBound}
	\end{equation}
	and $H(0) \neq 0$ then
	\begin{equation}
	H(\theta) = H(0)\exp(\mu \theta) \prod_{n = 0}^\infty (1 -\theta/u_n)\exp(\theta/u_n)
	\end{equation}
	and
	\begin{equation}
	\sum_{n = 0}^\infty |u_n|^{-2} < \infty
	\end{equation}
	where $\{u_n, n \geq 0\}$ is an enumeration of the zeroes of $H$ in with their full multiplicity, and $\mu$ is a complex constant.
\end{thm}
\begin{thm}
	Let $G$ be an entire function, nonvanishing at zero, which satisfies 	$\overline{G(\theta)} = G(\overline{\theta})$ and inequality \ref{hExpoBound}
	with zeroes $ {\mathcal Z}$ is equal to the disjoint union of ${\mathcal Z}_0$ and $ {\mathcal Z}_1 $
	where 
	${\mathcal Z}_1 = \{z_n , \overline{z}_n, 0 \leq n_1 < n < \infty \}$ such that    
	$\Im{z_n} = 2\pi|\alpha|^{-1}n +O(1)$ and
       $\Re{z_n} = O(n^\eta)$ where $\alpha < 0$ and $\eta < 1$,
  and where
  ${\mathcal Z}_0 =\{w_n, 0 \leq n < n_0 \}$.
	In addition suppose either $G$ is of exponential order $0$ as $x \rightarrow -\infty$, {\bf or}  $M = |\alpha|$ in inequality \ref{hExpoBound} and $G(\theta)$ is of exponential order is less than or equal to zero as $x \rightarrow -\infty$.

 Then
	$G(\theta)= G(0)p(\theta)\exp\{-\alpha\theta/2)\}
	\gamma(\theta)$, 
	
	where
	$\gamma(\theta) = 	\prod_{n=n_1}^{\infty}  (1- \frac{2\Re z_n\theta}{z_n\overline{z}_n}
		+\frac{\theta^{2}}{z_n\overline{z}_n})$ and
	 $p(\theta) = \prod_{n=0}^{n_0} (1 - \theta/u_n)$
	\label{inf_prod_characterization}
\end{thm}
{\bf Proof}
	By virtue of Theorem \ref{HadamardProdTheorem} there exists a constant $\mu$ such that
	\begin{equation}
	G(\theta) = G(0)\exp(\mu \theta) 
			\prod_{z \in {\mathcal Z}} (1 -\theta/z)\exp(\theta/z)
	\end{equation}
	Setting $\lambda = \mu +  \sum_{n=1}^\infty(z_n^{-1} + \overline{z}_n^{-1})
	+ \sum_{n=n_1}^{n_0} u_n^{-1}$ we may simplify the above equation to
		\begin{equation}
			G(\theta) = G(0)p(\theta)\exp\{\lambda\theta\}
			\prod_{n=n_1}^{\infty} (1- 
			\frac{2\Re z_n\theta}
			{z_n\overline{z}_n}
			+\frac{\theta^{2}}{z_n\overline{z}_n}) .
                \label{hHadRep}
		\end{equation}
		
	The assumption that $\overline{G(\theta)} = G(\overline{\theta})$ implies that $ {\mathcal Z}_0 = \overline{{\mathcal Z}_0}$ and in turn that $p(\overline{\theta}) = \overline{p(\theta)}$, and from that it may be inferred $\lambda$ is real valued.
		 
	Corollary \ref{expBehaviorEstCor} indicates that $\gamma$ is of
	exponential order $|\alpha|/2$ as $|\theta| \rightarrow \infty$. 
	We conclude by showing that $\lambda = |\alpha|/2 = -\alpha/2$.
	If $G$ has exponential order $0$ at $-\infty$ then Proposition \ref{prod_exp_order} may be applied to show that $\lambda = |\alpha|/2$. Otherwise we may use  Proposition \ref{prod_exp_order} in conjunction with hypothesis that the exponential order of $G(\theta)$ at $-\infty$ is no greater than 0 and the exponential order of $\gamma$ at $-\infty$ equals $|\alpha|/2$ to show
		 $\lambda \geq  |\alpha|/2$. Conversely we may show that $\lambda \leq |\alpha|/2$ by
	applying Proposition \ref{prod_exp_order} in conjunction with the fact that the exponential order of $\gamma$ at $\infty$ is $|\alpha|/2$ and the hypothesis that 
	$|G(\theta)| \leq C\exp(|\alpha\theta|)$. 

\begin{cor}
	Let $H$ be an entire function, nonvanishing at zero, which satisfies 	$\overline{H(\theta)} = H(\overline{\theta})$ and inequality \ref{hExpoBound} and
	whose multiset of zeroes is the union of the disjoint multisets ${\mathcal W}_0$ and ${\mathcal W}_1 $
	where 
	${\mathcal W}_1 = \{w_n , \overline{w}_n, 0 \leq n_1 < n < \infty \}$ such that    
	$\Im{w_n} = 2\pi\alpha^{-1}n +O(1)$ and
       $\Re{w_n} = O(n^\eta)$ where $\alpha < 0$ and $\eta < 1$,
  and where
  ${\mathcal W}_0 =\{v_n, 0 \leq n < n_0 \}$.
	In addition suppose either $H$ is of exponential order $0$ as $x \rightarrow \infty$, {\bf or}  $M = |\alpha|$ in inequality \ref{hExpoBound} and $H(\theta)$ is of exponential order is less than or equal to zero as $x \rightarrow \infty$.

 Then
	$H(\theta)= H(0)\hat p(\theta)\exp\{\alpha\theta/2)\}
	\hat\gamma(\theta)$, 
	
	where
	$\hat\gamma(\theta) = 	\prod_{n=n_1}^{\infty}  (1- 
	\frac{2\Re w_n\theta}{w_n\overline{w}_n}
		+\frac{\theta^{2}}{w_n\overline{w}_n})$ and
	 $\hat p(\theta) = \prod_{n=0}^{n_0} (1 - \theta/u_n)$
	\label{inf_prod_characterization_dual}
\end{cor}
{\bf Proof}
Set $G(\theta) = H(-\theta)$, $z_n = -w_n$ and $u_n = -v_n$. 
It follows that $G$ satisfies the hypotheses of Theorem
\ref{inf_prod_characterization}. Thus
$H(\theta)=G(-\theta) = G(-0)\exp(\alpha\theta)
p(-\theta)\gamma(-\theta)$.
To conclude the proof notice that simple algebra shows $H(0) = G(-0)$, $p(-\theta) = \hat p(\theta)$ and 
$\gamma(-\theta) = \hat \gamma(\theta)$.

\begin{defin}
A non-decreasing sequence $\{r_n , n > 0\}$ is {\bf $\eta$-sublinear} if there exists $n_0$ for which
$ r_{n+1} - r_n < \eta$  whenever $n \geq n_0$.
\label{subLinDef}
\end{defin}
\begin{thm}
	Let $f$ be a globally meromorphic non-vanishing function satisfying the following conditions.
	\begin{enumerate}
		\item There exists a sequence of circles $\rho_n$ centered around $0$ of radius $r_n$ such that
	        $r_n \uparrow \infty$ and $\limsup r_{n+1} - r_n < \eta$ and for which
	        $|f(z)|$ is bounded by $C\exp(\lambda|z|)$ on $\rho_n$
		\item All poles occur in conjugate pairs except for one non-zero real pole.
		\item The poles $\{z_n\}$ of the upper plane have asymptotic form
			\begin{equation}
			\Re{z_n} = o(n^{1/2})
			\end{equation}
			\begin{equation}
			\Im{z_n} = -2\pi\alpha^{-1} n  + \beta + o(1)
			\end{equation}
			where $\alpha < 0 $.
                \item $f$ is of exponential order $-|\alpha|$ at $-\infty$
                \item $f(0) = 1$
	\end{enumerate}
	Then 
	\begin{equation}
		f(\theta) = \exp(-\alpha \theta/2) \gamma(\theta)^{-1}
	\end{equation}
	where
	\begin{equation}
		\gamma(\theta) = (1 - \theta/z_0)
			\prod_{k=1}^{\infty} 
			(1- \frac{2\Re z_k\theta}
			{z_{k}\overline{z}_{k}}
			+\frac{\theta^{2}}{z_{k}\overline{z}_{k}} )
	\end{equation}
        \label{onThemeOfWeierstrass}

\end{thm}
{\bf Proof}
	Observe that $F = f \gamma$ is an analytic non-vanishing function. By Corollary \ref{prodExpoBound}
	\begin{equation}
		|\gamma(\theta)| < C_2\exp(|\theta|(|\alpha| /2+ \epsilon))
	\end{equation}
	Thus on $\rho_n$
	\begin{equation}
		|F(\theta)| < C_3\exp(r_nM)
	\end{equation}
	where $M = \alpha /2 + \epsilon + \lambda$ and $C_3 = C C_2$.
	The last inequality and Lemma \ref{maxPrinLem} jointly imply that the existence of
	$C_4$ such that $|F(z)| \leq C_4\exp(M|z|)$. Hadamard's product theorem in turn implies
	$F(z) = F(0)\exp(\kappa z)$. Notice that $F(0) = f(0)\gamma(0) = 1$.
	Since $f$ is of exponential order $\alpha$ at $-\infty$ and $\gamma$ is of exponential
	order $-\alpha/2$ at $-\infty$ it follows from Proposition \ref{prod_exp_order}
	that $\kappa = -\alpha/2$, completing the proof.

\begin{lem}
	\label{maxPrinLem}
	Let $h$ be a holomorphic function bounded by $C\exp(\lambda|z|)$ with $\lambda > 0$
	on $\rho_n$ where $\rho_n$ are concentric circles
        centered around $0$ of radius $r_n$ such that
	$r_n \uparrow \infty$ and $\limsup r_{n+1} - r_n < \eta$. Then $h$ is bounded by
	$C_1 \exp(\lambda|z|)$
	on the complex plane.
\end{lem}
{\bf Proof}
         The limit supremum condition implies the existence of $n_0$
         such that $r_{n+1} - r_n < \eta$ whenever $n \geq n_0$.
	 Set $C_1 = C\exp(\lambda \max(\eta, r_{n_0}))$.
	 Let $z$ be given. If $|z| \leq r_{n_0}$ then by the maximum principle
	 $|h(z)| \leq C_1 \leq C_1 \exp(\lambda |z|)$.
	 
	 For $|z| > r_{n_0}$ there exists $m > n_0$ for which
	 $r_m \geq |z| > r_{m - 1}$. Applying the maximum principle to the circle $\rho_m$ yields
	 \begin{equation}
	 |h(z)| \leq C \exp(\lambda (r_m - |z|) +
	                  \lambda |z|) \leq C\exp(\lambda\eta)\exp(\lambda|z|) \leq C_1 \exp(\lambda |z|)
         \end{equation}

\begin{defin}
	$f$ is said to be {\em globally meromorphic} if all its poles 
	$\{\xi_{n},\ -\infty < n < \infty \}$ are of finite order
	$k_n$ and are isolated, ie, there exists $r_0 > 0$ such that
	$f$ may be represented around $\xi_{n}$ by
		\begin{equation}
			f(z)=  \sum_{j= - k_n}^{\infty} (z - \xi_{n})^j A_{j}^{(n)}
		\end{equation}
	\label{mero}
	where $|z - \xi_n| < r_n$.
\end{defin}

\begin{defin}
    $f$ is said to satisfy {\em Mittag-Leffler Integrability} if there exist 
	a sequence
	$\{r_{n},\ n=1,2,...\}$
	increasing to $\infty$ and a sequence
	$\{\epsilon_{n},\ n=1,2,...\}$ decreasing to zero for which

	\begin{equation}
		|z| = r_{n} \Rightarrow |f(z)| < \epsilon_n |z|
		\label{linBound}
	\end{equation}
	\label{integProp}
\end{defin}

\begin{thm}
	Let $f$ satisfy definitions \ref{mero} and \ref{integProp} and analytic at $0$
	then uniformly for all bounded sets
	\begin{equation}
		f(z) = f(0) + \sum_{-\infty < n < \infty} \sum_{j=1}^{k_n}  A_{j}^{(n)}((z - \xi_{n})^{-j} - (-\xi_n)^{-j})
	\end{equation}
	\label{mittagTheorem}
	where $z$ is an analytic point of $f$.
\end{thm}

{\bf Proof}  Define
	$\rho_n = \{\zeta: |\zeta| = r_n\}$ and
	$g(\zeta, z) = f(\zeta)((\zeta - z) \zeta)^{-1}$.
	
	 The conditions on $f$ imply that
        for $\zeta \in \rho_n$, $g(\zeta, z) = o({r_n}^{-1})$ uniformly for $|z| \leq r_n /2$, and in turn implies that $\eta_n(z) = o(1)$ uniformly for $|z| \leq r_n/2$ where $\eta_n(z)$ is defined as
    \begin{equation}
    \eta_n(z) = \int_{\rho_n}g(\zeta,z)d\zeta = o(1)
    \end{equation}
	Applying Cauchy's residue theorem to $g$ on $\rho_n$
	shows that
	uniformly for $|z| < r_n / 2$ and $z$ an analytic point of $f$
		\begin{equation}
		    \lim_{n \rightarrow \infty}	f(z)/z - f(0)/z + \lim_{n \rightarrow \infty}
				\sum_{\xi_j \in {\rho_n} }\mbox {\bf res} \left(f(\zeta)((\zeta - z) 
				\zeta)^{-1} \right) \left|_{\zeta = \xi_j} \right. = 0
		\end{equation}
		
	Thus uniformly for analytic points $z$ of a bounded domain we have
		\begin{equation}
			f(z) = f(0) - 
			\lim_{n \rightarrow \infty}
			z \sum_{\xi_j \in {\rho_n} }\mbox {\bf res}
					  \left(f(\zeta)((\zeta - z) \zeta)^{-1}  \right) \left |_{\zeta = \xi_j} \right.
			\label{sumOfRes}
		\end{equation}
	The proof of the theorem follows jointly from \ref{sumOfRes} and Lemma \ref{resComputation}

\begin{lem}
	\label{resComputation}

	Let $f$ have a pole of degree $k$ at $\xi$ then
	\begin{equation}
		-z\mbox {\bf res}(f(\zeta)\left((\zeta - z)\zeta\right)^{-1} \left|_{\zeta = \xi } \right. = \sum_{j = -1}^{-k}
		c_j((z - \xi)^j - (-\xi)^j)
		\label{gResidue}
	\end{equation}
	where locally
	\begin{equation}
		f(\zeta)=  \sum_{j= - k}^{\infty} c_j (\zeta -  \xi)^j
		\label{fLaurent}
	\end{equation}
\end{lem}
{\bf Proof} Define $h(\zeta) =  z(\zeta - z)^{-1} \zeta^{-1}$.
Straightforward algebra shows that
          \begin{equation}
         h(\zeta) = (\zeta - z)^{-1} - \zeta^{-1}
         \end{equation}
and
          \begin{equation}
                          (\zeta -  \xi + b)^{-1}  = -\sum_{j = 0}^{\infty}(\zeta -  \xi)^j (-b)^{-j-1}
          \end{equation}
Combining the previous identities yields
      \begin{equation}
               -h(\zeta) =
               \sum_{j = 0}^{\infty}(\zeta -  \xi)^j  ((z -  \xi)^{-j -1} - (- \xi)^{-j -1})
               \label{hExpansion}
      \end{equation}
The residue of $-hf$ at $\zeta =  \xi$ is the coefficient of $(\zeta -  \xi)^{-1}$ term in the
Laurent expansion of $-h f$. Multiplying out the right-hand sides of \ref{fLaurent} and \ref{hExpansion}
yields \ref{gResidue}.
\begin{lem} Let $t \neq 0$ and $z \neq 0$ be fixed then
\begin{equation}
	\label{erlangInvFormula}
	{\rm Res}((\theta - z)^{-j} \theta^{-1} \exp(\theta t) , z) = 
	\sum_{k = 0}^{j - 1} z^{k-j}(-1)^{j-k-1} t^k\exp(zt)/k!
\end{equation}
where ${\rm Res}(f(\theta), z))$ denotes $ {\bf res}(f(\theta))\left|_{\theta = z } \right. $.
\end{lem}
{\bf Proof} If $h$ is analytic in a neighborhood $z$, then
${\rm Res}((\theta - z)^{-j} h(\theta), z) = h^{(j-1)}(z)/(j-1)!$ where $h^{(k)}$ denotes the 
$k$th derivative of $h$. Applying Leibniz's identity for products to $h(\theta) = \theta^{-1}\exp(\theta t)$ shows the r.h.s.
of \ref{erlangInvFormula} equals $h^{(j-1)}(z)/(j-1)!$ . 
\begin{lem}
	\label{LaplaceInversionTheorem}
	Let $F$ be the distribution function of a nonnegative random variable with Laplace transform $\psi$. Then for $c > 0 $ and $t > 0$

	\begin{equation}
		\lim_{T \rightarrow \infty} 
		     \int_{c - iT}^{c + iT} \psi(\theta)\theta^{-1}\exp(\theta t)d\theta =
		(F(t{\scriptscriptstyle +}) + F(t{\scriptscriptstyle -}))\pi i
		\label{LaplaceInversionFormula}
	\end{equation}
For a proof consult Widder \cite{widder}, Theorem 7.6a, page 69.
\end{lem}

\begin{cor}
\label{laplaceInversionCorollary}
Suppose there exist a constant $M$ and a sequence of radii $r_n$ increasing to infinity such that $|\psi(z)| \leq M$ whenever $|z| = r_n$. Then under the conditions of Lemma \ref{LaplaceInversionTheorem}
	\begin{equation}
	\label{lapIntegralInvEqn}
		\lim_{n \rightarrow \infty} 
		     \int_{L_n} \psi(\theta)\theta^{-1}\exp(\theta t)d\theta =
		(F(t{\scriptscriptstyle +}) + F(t{\scriptscriptstyle -}))\pi i
	\end{equation}
	where $L_n = A_n \bigcup D_n$
	with $A_n = \{z: |z| = r_n \: and \:\Re {z} < c \}$ and
	$D_n = \{z = c + iu: u \in [-T_n, T_n] \}$ with $T_n = (r_n^2 - c^2)^{1/2}$

\end{cor}
{\bf Proof} By virtue of Lemma \ref{LaplaceInversionTheorem} it suffices to show that the integral on ${\mathcal A}_n$ vanishes as $n$ grows large. 
\begin{equation}
\int_{A_n} \psi(\theta)\theta^{-1}\exp(\theta t)d\theta =
\int_{\Omega_n} i\psi(r_n\exp(i\omega))\exp(r_n\exp(i\omega) t)d\omega
\end{equation}
where $\theta = r_n\exp(i\omega)$ and where $\Omega_n =  [\omega_0 , \omega_3 ]$
and $\omega_0 = \arccos(c/r_n) $ and $\omega_3 = 2\pi - \omega_0$.
\begin{equation}
\left | \int_{\Omega_n} i\psi(r_n\exp(i\omega))\exp(r_n\exp(i\omega) t)d\omega \right |
\leq \int_{\Omega_n} M \exp(r_n\cos(\omega)t) d\omega
\end{equation}
Defining ${\mathcal B}_n =  [\omega_1 , \omega_2 ] $
where $\omega_1 = \arccos(-r_n^{-1/2})$ and
$\omega_2 = 2\pi - \omega_1$, we find that $\cos(\omega)$ is dominated on ${\mathcal B}_n$  by $\cos(\omega_1)$
and $\Omega_n - {\mathcal B}_n$ by
$\cos(\omega_0)$ respectively.
Noting that $\cos(\omega_1) = -r_n^{-1/2}$ and that
$\cos(\omega_0) = cr_n^{-1}$ we
 deduce the following inequality.
\begin{equation}
\int_{\Omega_n} M \exp(r_n\cos(\omega)t) d\omega \leq
\int_{{\mathcal B}_n} M\exp(-r_n^{1/2} t)d\omega
+ \int_{\Omega_n - {\mathcal B}_n} M \exp(ct) d\omega
\end{equation}
Noting the concavity of  $\sin$ on the interval $[0, \pi/2]$ one may show that $\arcsin(x) < x\pi/2$ for $x \in (0, \pi/2)$. From this we see that the arc length of
$\Omega_n - {\mathcal B}_n$ is bounded by $r_n^{-1/2} + cr_n^{-1}$.

\begin{flushleft}
Combining the bounds on the integrands and arc lengths shows the integral over ${\mathcal B}_n$ is no more than $\pi M \exp(-r_n^{1/2} t)$,
and the integral over $\Omega_n - {\mathcal B}_n$ is dominated by 
$\pi M\exp(ct)  (r_n^{-1/2} +  cr_n^{-1})$. Hence 
$\int_{A_n} \psi(\theta)\theta^{-1}\exp(\theta t)d\theta \rightarrow 0$ as $n \rightarrow \infty$.
\end{flushleft}

\begin{lem}
	\label{WaitingTimeExpFromPsi}
    Let $F$ be the distribution function of a nonnegative random variable with Laplace transform    
    $\psi$. Suppose $\psi$ satisfies the conditions of Lemma \ref{LaplaceInversionTheorem} and
    Corollary \ref{laplaceInversionCorollary}. Additionally suppose that $\psi$ has a partial fraction expansion given by Equation \ref{psiPartFracExpansionEqn}. Let $t > 0$ be a point of continuity of $F$ then

	\begin{equation}
	\label{WaitingTimeExpFromPsiEqn}
		F(t) = 1 -	\lim_{N \rightarrow \infty} \sum_{n= -N}^N \sum_{j=1}^{k_n} a_{n,j}
		\sum_{l = 0}^{j - 1} \exp(z_nt)(-tz_n)^l / l! 
	\end{equation}

\end{lem}
where $a_{n,j} = A_{n,j}(-z_n)^{-j}$

{\bf Proof} Equation \ref{lapIntegralInvEqn} and the continuity of $F$ at $t$ imply
\begin{equation}
\label{lapIntegralInvResEqn}
F(t) = 
\lim_{n \rightarrow \infty} (2\pi i)^{-1}
		     \int_{L_n} \psi(\theta)\theta^{-1}\exp(\theta t)d\theta
\end{equation}
where $L_n$ is the curve defined in Lemma \ref{laplaceInversionCorollary}.

Applying Cauchy's residue theorem to the right hand side of \ref{lapIntegralInvResEqn} yields
\begin{equation}
\label{lapLimInvResEqn}
F(t) = 1 + \lim_{m \rightarrow \infty}\sum_{z \in \Lambda_m} {\rm Res}(\psi(\theta) \theta^{-1} \exp(\theta t) , z) 
\end{equation}
where $\Lambda_m$ is the interior of the simply closed curve $L_m$.

Equations \ref{psiPartFracExpansionEqn}, \ref{erlangInvFormula} and \ref{lapLimInvResEqn} jointly imply
\begin{equation}
\label{AsubnFormOfCumFunct}
F(t) = 1 + \lim_{m \rightarrow \infty}\sum_{z_n \in \Lambda_m} 
\sum_{j=1}^{k_n} A_{n,j} \sum_{l = 0}^{j - 1} z_n^{l-j}(-1)^{j-l-1} t^l\exp(z_nt)/l!
\end{equation}
Substituting $a_{n,j} = A_{n,j}(-z_n)^{-j}$ into \ref{AsubnFormOfCumFunct} and simplifying yields
\begin{equation}
F(t) = 1 - \lim_{m \rightarrow \infty}\sum_{z_n \in \Lambda_m} 
\sum_{j=1}^{k_n} a_{n,j} \sum_{l = 0}^{j - 1} \exp(z_nt)(-tz_n )^l /l!
\end{equation}
Finally, by supposition $\{z_n \in \Lambda_m\} = \{z_n: n\in [-m, m]\}$ 
reduces the last equation to the r.h.s. of \ref{WaitingTimeExpFromPsiEqn}.

\begin{lem}
\label{ahlforsRationalFunctionLemma}
Let $R(z) = p(z)/q(z)$ where $p$ and $q$ are polynomials. If 
$\lim_{x \rightarrow \infty} R(x) = a$ where $a$ is finite and non-zero, then $\deg(p) = \deg(q)$.
\end{lem}
The limit condition implies that $R$ has neither a pole or zero at $\infty$. 
The lemma then follows from a comment of Ahlfors \cite{ahlfors} in section 1.3, page 44.

\begin{cor}
Let $H(z) = R(z)A(z)$ where $R$ is the ratio of polynomials $p$ to $q$. If 
$\lim_{x \rightarrow \infty} H(x) = b$ and
$\lim_{x \rightarrow \infty} A(x) = c$ where $b$ and $c$ are non-zero and finite, then $\deg(p) = \deg(q)$.
\end{cor}
{\bf Proof} Since the limit of a ratio is equal to the ratio of the limits whenever the denominator limit is non-zero, we have
$\lim_{x \rightarrow \infty} R(x) = a $ where $a = b/c$. Since $b$ and $c$ are non-zero finite so is $a$.

\begin{cor}
\label{degCompForInfProd}
Let $H = \exp(-\alpha\theta)h_1(\theta)h_2^{-1}(\theta) h_3(\theta)$ and 
$G = \exp(-\alpha\theta)g_1(\theta)g_2^{-1}(\theta) g_3(\theta)$ where
$h_1, h_2, g_1, g_2$ are all polynomials, and 
\begin{equation}
h_3(\theta) = \prod_{n=n_0}^\infty (1-\theta/z_n)(1-\theta/\overline{z}_n)
\end{equation}
\begin{equation}
g_3(\theta) = \prod_{n=n_0}^\infty (1-\theta/w_n)(1-\theta/\overline{w}_n)
\end{equation}
such that $w_n = z_n + O(n^{-\eta})$ for $\eta > 0$ and $|z_n| > Mn > 0$ and $\Re z_n < 0$ for $n$ sufficiently large. If $\lim_{x \rightarrow \infty} |H(x)/G(x)| = b \in (0, \infty)$ then
$\deg(h_1 g_2) = \deg(h_2 g_1)$.
\end{cor}
{\bf Proof} Notice $H(x)/G(x)$ can be expressed as $R(x)A(x)$ where $R(x) = h_1(x)g2(x)h_2^{-1}(x)g_1^{-1}(x)$ and $A(x)=h_3(x)g_3^{-1}(x)$. From Lemma \ref{limitOfMobiusRatios} we deduce that $A(x)$ has a finite, positive limit at $\infty$, and consequently Lemma \ref{ahlforsRationalFunctionLemma} may be applied to complete the proof.
\begin{lem}
\label{limitOfMobiusRatios}
Let $h_3$ and $g_3$ be as in Corollary \ref{degCompForInfProd}
then $\lim_{x \rightarrow \infty} h_3(x)/g_3(x)$ is $ \prod_{k = n_0}^\infty |z_k|^2
|w_k|^{-2})$ which is finite and positive.
\end{lem}
{\bf Proof}
\begin{equation}
 g_3(x)/h_3(x) = \prod_{k = n_0}^\infty \frac{z_k\overline{z}_k (w_k - x)(\overline{w}_k -x)}
{w_k\overline{w}_k (z_k - x)(\overline{z}_k - x) }
\end{equation}
Observe that 
\begin{equation}
\label{zkToWkRatioEst}
z_k  w_k^{-1} = 1 + (w_k - z_k)w_k^{-1} = 1 + O(k^{-1-\eta})
\end{equation} 
This proves that
\begin{equation}
\prod_{k = n_0}^\infty |z_k|^2 |w_k|^{-2}
\end{equation}
is properly convergent.
We will now show that
\begin{equation}
\label{infiniteRatio_w_z_limit}
\lim_{x \rightarrow \infty}\prod_{k = n_0}^\infty \frac{(w_k - x)(\overline{w}_k -x)}
{(z_k - x)(\overline{z}_k - x) } = 1
\end{equation}
By hypothesis there exists $k_0 \geq n_0$ such 
$\Re{z_k} < 0$ and $|z_k| > M k$ whenever $k \geq k_0$. If we set 
\begin{equation}
x_0 = 2\max_{n_0\leq k < k_0}\{|\Re{z_k}|\}
\end{equation}
then
$\max(x/2,|z_k|) < |z_k - x|$ for $x > x_0$. Whence it follows for $x > x_0$
\begin{equation} 
\label{wk_xToZk_xRatioEst}
1 - |(w_k - z_k)z_k^{-1}| \leq |(w_k - x)(z_k - x)^{-1}| \leq 1 + |(w_k - z_k)z_k^{-1}|
\end{equation}
Thus for $x > x_0$ and $k > k_0$

\begin{equation} 
\label{wkzkTermSandwich}
1 - O(k^{-1-\eta}) \leq |(w_k - x)(z_k - x)^{-1}| \leq 1 + O(k^{-1-\eta})
\end{equation} 
Note for fixed $k$ we have
\begin{equation}
\label{wkzkxTermWiseLimit}
\lim_{x \rightarrow \infty} (w_k - x)(z_k - x)^{-1} = 1 
\end{equation}
By the dominated convergence theorem for products (cf. Walker\cite{walker}) we may deduce \ref{infiniteRatio_w_z_limit}
from \ref{wkzkTermSandwich}  and \ref{wkzkxTermWiseLimit}.

\begin{lem}
\label{zeroApproxLemma}
Let $f = f_0 - 1$ and $g = f + f_0g_0$ be analytic on a domain ${\mathcal D}$. Suppose $f$ has a sequence of simple zeroes $\{z_n \in {\mathcal D}, n > 0\}$ for which
there exist positive constants
$\delta > 0$, $\mu > 0$, $\nu > 0$, $\lambda > 0$, $M > 0$ for which $|\eta| < \delta$ implies $z_n + \eta \in {\mathcal D}$, 
$\mu|\eta| < |f(z_n + \eta)| < \nu|\eta|$, and
$|g_0(z_n + \eta)| < Mn^{-\lambda}$. Then there exist $n_0$ and a sequence $\{w_n \in {\mathcal D}, n > n_0\}$ such that $g(w_n) = 0$ and $|w_n - z_n| < 2M\mu^{-1}n^{-\lambda}$.
\end{lem}
{\bf Proof}
Let $n_0$ be so large that if $n > n_0$ then $\delta_n < \delta$  and
$\delta_n \nu < 1$, where
$\delta_n = 2M\mu^{-1}n^{-\lambda}$. We claim that if $|\eta| = \delta_n$ then  $|f_0(z_n + \eta)g_0(z_n + \eta)| < \mu\delta_n$. By hypothesis $|\eta| < \delta$ then
$|f_0(z_n + \eta)| < 1 + \nu|\eta|$ , and 
$|g_0(z_n + \eta)| < Mn^{-\lambda}$. Thus
$|f_0(z_n + \eta)g_0(z_n + \eta)| < (1 + \nu|\eta|)Mn^{-\lambda} < 2Mn^{-\lambda} = \delta_n \mu$.

Since $f$, and $f g_0$ are analytic in ${\mathcal D}$, and $|f(z_n + \eta)| > \mu\delta_n$ we may deduce from Rouche's Theorem the existence of 
$w_n$ such that $ |z_n - w_n| < \delta_n$. QED

The following lemma complements Lemma \ref{zeroApproxLemma}

\begin{lem}
\label{rouchForTFunctions}
Let $f = f_0 - 1$ and $g = f + f_0g_0$ be analytic on a simply connected domain ${\mathcal D}$ containing a Jordan curve ${\mathcal C}$. Suppose $|f(z)| > 4\delta$ and $|g_0(z)| < \delta \leq 0.1$ for $z \in  {\mathcal C}$, then  $f$ and $g$ have the same number of zeroes inside of ${\mathcal C}$. 
\end{lem}
{\bf Proof}
By Rouche's theorem it suffices to prove that $|f| > |f - g|$ on $ {\mathcal C}$ which is implied by Prop \ref{roucheIneqProp}. 

\begin{prop}
\label{roucheIneqProp}
Let $f_0$ and $g_0$ be functions defined on a set ${\mathcal A}$ satisfying the following inequalities:
$|f_0 - 1| > 4 \delta > 0$ and $|g_0| < \delta \leq 0.1$, then $|f_0 - 1| - |f_0g_0| > 2\delta$ on ${\mathcal A}$.
\end{prop}
{\bf Proof} Let $z \in {\mathcal A}$ such that $|f_0(z)| \leq 2$ then
$|f_0(z)g_0(z)| \leq 2 |g_0(z)| < 2 \delta < |f_0(z) - 1| - 2\delta$.
Conversely suppose $|f_0(z)| > 2$ then $|f_0(z)| (1 - 4\delta) > 1$. Thus
$|f_0(z) - 1| \geq |f_0(z)| - 1 > 4\delta |f_0(z)| > 2\delta |f_0(z)| > |g_0(z)f_0(z)|$ Thus
$|f_0(z) - 1| - |g_0(z)f_0(z)| > 4\delta |f_0(z)| - 2\delta |f_0(z)| > 2\delta$.

\begin{prop}
\label{bigOhTransitiveProp}
Let $f(\theta)$, $g(\theta)$ and $h(\theta)$ be functions on a domain  ${\mathcal A}$ satisfying the following inequalities:
$g(\theta)/f(\theta) = O(|\theta|^{-c})$ where $c > 0$, $|h(\theta)/f(\theta)| = o(1)$ then
$|g(\theta)/(f(\theta) + h(\theta))| =  O(|\theta|^{-c})$
\end{prop}
{\bf Proof}
\begin{equation}
g(\theta)/(f(\theta) + h(\theta)) = g(\theta)f(\theta)^{-1}(1 + h(\theta)/f(\theta))^{-1}
\end{equation}
Thus arguing somewhat informally we find
\begin{equation}
|g(\theta)/(f(\theta) + h(\theta))| =  O(|\theta|^{-c}) (1 + o(1))^{-1} = O(|\theta|^{-c})
\end{equation}
\begin{lem} 
\label{ExtendedFactorEstimatingLemma}
Suppose $\Phi_j$ for $j \in \{0, ... J-1\}$ are functions bounded in a neighborhood of zero and where $x(\theta) = -\gamma\log|\theta| +O(1)$ and $\gamma > 0$.
Define $\Psi(\theta)$ by
\begin{equation}
\Psi(\theta) = \sum_{0 \leq j < J} \exp(\beta_j\theta) \Phi_j(\theta ^{-1})\theta^{-n_j}
\end{equation}
Then 
\begin{equation}
|\Psi(\theta)| = O(|\theta|^{-m( \hat{\beta},\hat{ n})})
\end{equation}
 
where $\hat{n} = (n_0, ..., n_{J-1}) $ and $\hat{\beta} = (\beta_0, ..., \beta_{J-1}) $
and $m( \hat{\beta},\hat{ n}) = \min_{0 \leq j < J}\gamma\beta_j + n_j $. Moreover if 
$\beta_j \geq 0$, $n_j \geq 0$ 
and $\beta_j + n_j > 0$ for $j < J$ 
then for every $\epsilon > 0$ there exists an $a < $ such that if $\Re{\theta} \leq a$ then $|\Psi(\theta)| < \epsilon$.
\end{lem}
{\bf Proof}
From Proposition \ref{localEstimateOfHterms} we find that
\begin{equation}
\Psi(\theta)= \sum_{0\leq j < J} O(|\theta|^{-\gamma\beta_j -n_j }) = \max_{0\leq j < J} O(|\theta|^{-\gamma\beta_j -n_j})= O(|\theta|^{-m( \hat{\beta},\hat{ n})})
\end{equation}
\begin{prop}
\label{localEstimateOfHterms}
Let $\psi_{\beta, n}(\theta) = \exp(\beta\theta) \Phi(\theta ^{-1})\theta^{-n}$ 
where $\Phi$ is bounded in a neighborhood of zero and where $x(\theta) = -\gamma\log|\theta| +O(1)$ and $\gamma > 0$ 
then $|\psi_{\beta, n}(\theta)| = O(|\theta|^{-\beta\gamma -n})$. Moreover if $\beta$ and $n$ are nonnegative and their sum positive then for every $\epsilon > 0$ there exists an $a$ such that if $\Re{\theta}< a$ then $\Psi_{\beta,n}(\theta)<\epsilon$.
\end{prop}
{\bf Proof}
By hypothesis there exist $M_0$ and $\delta > 0$ such that if $|w| < \delta$ then $|\Phi(w)| < M_0$ hence if $|\theta| > \delta^{-1}$ then $|\Phi(\theta)| < M_0$.Also by hypothesis there exist $M_1$ and $C$ such that $|\theta| > C$  implies
$|x(\theta) + \gamma\log|\theta|| < M_1$. Putting it all together if $|\theta| > \max(C, \delta^{-1})$ we have
\begin{equation}
|\psi_{\beta, n}(\theta)| < \exp(-\beta\gamma\log|\theta| +\beta M_1) 
M_0 |\theta|^{-n} = |\theta|^{-\beta\gamma -n}D 
= O(|\theta|^{-\beta\gamma -n})
\end{equation}
where $D = M_0 \exp(\beta M_1)$.

Observe if $x(\theta) \leq a < 0$ and $|a| > \delta^{-1}$ then
\begin{equation}
|\psi_{\beta, n}(\theta)| < M_0\exp(\beta a)|a|^{-n} =  M_0 |a|^{\log_{|a|}(e)\beta a -n}
\label{psiBoundOnX}
\end{equation}
Since the last term has a negative exponent the term can be made as small as possible by letting $a$ tend to $-\infty$ thus proving the proposition's last claim.

\begin{prop}
\label{continuityOfTsRootsRealPartEqn}
Let $z = x(z) + i y(z)$ where $x(z) = -\gamma\log|z| +\beta$ where $\gamma > 0$ and $\beta$ is real then for every $\delta < 1/2$ and for $|z| > 3\max(\gamma, 1)/2$ 
then if $|\theta - z| < \delta$ we have $|x(\theta) +\gamma\log|\theta| -\beta| < 2\delta$
\end{prop}
{\bf Proof}
\begin{equation}
\log|\theta| = \log|z| + \log|1 + (\theta  -z)/z| = \log|z| + u(z,\theta) 
\end{equation}
where $|u(\theta, z)| < 3  |(\theta - z)/2z|$ by virtue of Proposition \ref{logAbsLocBehavior}.

Multiplying the previous equation by $-\gamma$ we find
\begin{equation}
-\gamma\log|\theta| = -\gamma\log|z| -\gamma u(z, \theta) = x(z) - \gamma u(z, \theta)
\end{equation}
Now subtracting $x(\theta)$ from both sides we get
\begin{equation}
-\gamma\log|\theta| - x(\theta)  =  x(z)- x(\theta) - \gamma u(z, \theta)
\end{equation}
Now taking absolute values and applying the triangle inequality we find
\begin{equation}
| x(\theta) +\gamma\log|\theta| -\beta|  \leq  |x(z)- x(\theta)| + \gamma |u(z, \theta)| < \delta + 3\delta\gamma/(2z)
\end{equation}
The desired inequality now follows from the assumption on $z$.

\begin{prop}
$|\log|1 + u|| \leq -\log(1-|u|) \leq 3|u|/2$ for $|u| < 1/2$
\label{logAbsLocBehavior}
\end{prop}
{\bf Proof}
If $|1 + u| < 1$ then by the monotonicity of the $\log$ function and the triangle inequality we find that
$\log(1-|-u|) \leq \log|1 + u| < 0$ or reversing signs we find $|\log|1 + u|| \leq - \log(1 - |u|)$. Conversely, if $|1 + u| \geq 1$ then also by monotonicity and the triangle inequality we have $0 < \log|1 + u| \leq \log(1 + |u|)$ Inequality (4.1.38) of \cite{abram} yields $\log(1 + |u|) < |\log(1 -|u|)$, and inequality (4.1.35) ibid
gives $|\log(1 -|u|) \leq 3|u|/2$ when $|u| < 1/2$.
  
\mysection{Utility Summation and Product Formulas}
The following formulas are due to Euler.
\begin{equation}
\omega(b) =\sum_{j = 0}^{\infty} (b^2 + j^2)^{-1} = ( b\pi\coth(b\pi) - 1) )
2^{-1}b^{-2}
\end{equation}
\begin{equation}
\label{eulerSum}
\Omega(b) =\sum_{j = 0}^{\infty} (b^2 + (2\pi j)^2)^{-1} = ( 2^{-1}b\coth(2^{-1}b) - 1) )
2^{-1}b^{-2}
\end{equation}
We may employ the last identity to telescope the following infinite product.
\begin{equation}
\label{eulerProdApprox}
\prod_{j = n}^{\infty}(1 +  2\lambda\mu(\lambda^2 + 4\pi^2 j^2)^{-1}) =
\exp\{2\lambda\mu( \Omega(\lambda)
- \sum_{j = 1}^{n-1}(\lambda^2 + 4\pi^2j^2)^{-1})\} + O(n^{-3})
\end{equation}

\begin{equation}
|1 - \prod_{j = 1}^n (1 + z_j)| \leq 2 \sum_{j = 1}^n |z_j| 
\end{equation}
whenever $\sum_{j = 1}^n |z_j| < 1$. 

\mysection{Conclusions and Remaining Challenges}
This paper furnishes a variety of computationally efficient means to determine the steady-state unfinished workload statistics for an uncommon set of G/G/1 models. (An alternative approach of great power is given in \cite{whittAbateChoudhury}.)
In particular the formulas derived for the the gated M/M/1's unfinished workload statistics should have interesting implications for investigations into queues with non-homogeneous Poisson arrival processes (cf \cite{closApproxForTimeDepQueue}).

This new class of solved G/G/1 queues should also be valuable to those studying sensitivity of waiting times to distributional assumptions.

Within the current scope there still remain some interesting challenges. For example, expansion given by equation \ref{waitingTailExpansion} converges very slowly for small values of $t$, so a way to accelerate this convergence or replace it with some other inversion formula would be useful. One potentially fruitful approach would be finding a more computationally efficient way to compute the spectral coefficients $a_{n,1}$ for $n$ large.

A promising idea is given in Section \ref{asymCompOfSpctrCoefSection} where it shown that the spectral coefficients depend on the asymptotic behavior of the idle distribution. In the case of bounded interarrival times the idleness transform should be amenable to asymptotic expansion since the idle times will also be bounded.

Another challenge is dealing with the case where $F$ is not dominated on the left plane by a term of the form $c\exp(\alpha\theta)\theta^{-m} - 1$ but rather by a product of such terms. The analysis of such $F$ was briefly considered in Section \ref{ExploreOfMultiplyDominatedQueue}


\begin{thebibliography} {dummyTitle}
\bibitem{abram}Milton Abramowitz and Irene A. Stegun
{\em Handbook of Mathematical Functions}, (Dover Publications, 1965).
\bibitem{ahlfors}Lars V. Ahlfors,
{\em Complex Analysis}, (MacGraw-Hill 1953).
\bibitem{apostol}Tom M. Apostol,
{\em Mathematical Analysis} (Addison Wesley, 1957)
\bibitem{bleistein}Norman Bleistein  and Richard A. Handelsman, {\em Asymptotic Expansions of Integrals} (Dover Press, 1986, Theorem 3.2, Chapter 3, p 78-79). 
\bibitem{cramer_cit}  Harald Cramer, {\em Mathematical Methods of Statistics}
(Princeton University Press, 1957)
\bibitem{que_con} D. Kennedy,
The continuity of the single-server queue,
J. Appl. Probab. 9, 370-381.
\bibitem{lev-red}N. Levinson and P. Redheffer,
{\em Complex Variables} (Holden-Day, 1970).
\bibitem{marshall}K.T. Marshall,
Some Relationships between the Distributions of Waiting Time, Idle Time, and
Interoutput Time in the GI/G/1 Queue,
SIAM Journal Applied Math, Vol 16,(1968) 324-327.
\bibitem{closApproxForTimeDepQueue}
Michael. H. Rothkopf and Shmuel. S. Oren, A Closure Approximation for the Nonstationary M/M/S Queue, Management Science, Vol 6 (1979) 522-534, No 11.
\bibitem{smith}W.L. Smith,
On the distribution of queueing times,
Proc. Cambridge Philosophical Society Vol 49 (1953) 449-461.
\bibitem{rootFinder}M. J. Schaefer, T. Bubeck,
A parallel complex zero finder, Journal of Reliable Computing,
Vol 1 (1995) 317-323, No 3.
\bibitem{takacs}Lajos Takacs,
{\em Combinatorial Methods in the Theory of Stochastic Processes} (John Wiley, 1967)
\bibitem{walker}Peter Leslie Walker, {\em Examples and Theorems in Analysis} (Springer, 2003, Theorem 6.76, Chapter 6, p 187).
Equation 20, Chapter 3, Section 19.
\bibitem{whittAbateChoudhury}Joseph Abate, Gagan L. Choudhury and Ward Whitt,
Calculation of the GI/G/1 Waiting-Time Distribution and its Cumulants from Pollaczek's Formula. Archiv für Elektronik und Übertragungstechnik, Vol 47 (1993) 311-321, Nos 5/6.
\bibitem{widder}D. Vernon Widder,{\em The Laplace Transform} (Princeton University Press, 1946).

\end{thebibliography}
\end{document}